\documentclass[11pt,a4paper]{article}

\usepackage[table]{xcolor}
\usepackage{slashed}
\usepackage{xkeyval}
\usepackage{url,amsmath,amssymb,latexsym,pstricks,mathrsfs,comment,amsthm,graphicx,tikz,tikz-cd,enumerate,accents,pgffor,cite,wrapfig,multicol,float}
\usepackage{bibspacing}
\newcommand{\nc}{\newcommand}
\nc{\rnc}{\renewcommand}

\nc{\OEIS}{}

\usepackage{geometry}  \rnc{\ss}{\smallskip} \nc{\ms}{\medskip}  \nc{\nss}{\vspace{-3mm}} 

\makeatletter
\newenvironment{myalign}{%
  \def\align@preamble{
    &\hfil
     \strut@
     \setboxz@h{\@lign$\m@th\displaystyle{####}$}%
     \ifmeasuring@\savefieldlength@\fi
     \set@field
     \hfil
     \tabskip\z@skip
    &\setboxz@h{\@lign$\m@th\displaystyle{{}####}$}%
     \ifmeasuring@\savefieldlength@\fi
     \set@field
     \hfil
     \tabskip\alignsep@
  }\align}{\endalign}

\DeclareMathSymbol{\widehatsym}{\mathord}{largesymbols}{"62}
\newcommand\lowerwidehatsym{%
  \text{\smash{\raisebox{-1.3ex}{%
    $\widehatsym$}}}}
\newcommand\fixwidehat[1]{%
  \mathchoice
    {\accentset{\displaystyle\lowerwidehatsym}{#1}}
    {\accentset{\textstyle\lowerwidehatsym}{#1}}
    {\accentset{\scriptstyle\lowerwidehatsym}{#1}}
    {\accentset{\scriptscriptstyle\lowerwidehatsym}{#1}}
}
\rnc{\widehat}{\fixwidehat}

\begin{document}


\nc{\PBnt}{R^\tau[\PB_n]}
\nc{\Mnt}{R^\tau[\M_n]}
\nc{\PBntC}{\bbC^\tau[\PB_n]}
\nc{\MntC}{\bbC^\tau[\M_n]}
\nc{\bbC}{\mathbb C}
\nc{\RSr}{R[\S_{\br}]}
\nc{\fs}{\mathfrak s}
\nc{\ft}{\mathfrak t}
\nc{\fu}{\mathfrak u}
\nc{\fv}{\mathfrak v}
\nc{\rad}{\operatorname{rad}}
\nc{\dominates}{\unrhd}

\nc{\ubluebox}[2]{\bluebox{#1}{1.7}{#2}2\udotted{#1}{#2}}
\nc{\lbluebox}[2]{\bluebox{#1}0{#2}{.3}\ldotted{#1}{#2}}
\nc{\ublueboxes}[1]{{
\foreach \x/\y in {#1}
{ \ubluebox{\x}{\y}}}
}
\nc{\lblueboxes}[1]{{
\foreach \x/\y in {#1}
{ \lbluebox{\x}{\y}}}
}

\nc{\bluebox}[4]{
\draw[color=blue!20, fill=blue!20] (#1,#2)--(#3,#2)--(#3,#4)--(#1,#4)--(#1,#2);
}
\nc{\redbox}[4]{
\draw[color=red!20, fill=red!20] (#1,#2)--(#3,#2)--(#3,#4)--(#1,#4)--(#1,#2);
}

\nc{\bluetrap}[8]{
\draw[color=blue!20, fill=blue!20] (#1,#2)--(#3,#4)--(#5,#6)--(#7,#8)--(#1,#2);
}
\nc{\redtrap}[8]{
\draw[color=red!20, fill=red!20] (#1,#2)--(#3,#4)--(#5,#6)--(#7,#8)--(#1,#2);
}

\usetikzlibrary{decorations.markings}
\usetikzlibrary{arrows,matrix}
\usepgflibrary{arrows}
\tikzset{->-/.style={decoration={
  markings,
  mark=at position #1 with {\arrow{>}}},postaction={decorate}}}
\tikzset{-<-/.style={decoration={
  markings,
  mark=at position #1 with {\arrow{<}}},postaction={decorate}}}
\nc{\Unode}[1]{\draw(#1,-2)node{$U$};}
\nc{\Dnode}[1]{\draw(#1,-2)node{$D$};}
\nc{\Fnode}[1]{\draw(#1,-2)node{$F$};}
\nc{\Cnode}[1]{\draw(#1-.1,-2)node{$\phantom{+0},$};}
\nc{\Unodes}[1]{\foreach \x in {#1}{ \Unode{\x} }}
\nc{\Dnodes}[1]{\foreach \x in {#1}{ \Dnode{\x} }}
\nc{\Fnodes}[1]{\foreach \x in {#1}{ \Fnode{\x} }}
\nc{\Cnodes}[1]{\foreach \x in {#1}{ \Cnode{\x} }}
\nc{\Uedge}[2]{\draw[->-=0.6,line width=.3mm](#1,#2-9)--(#1+1,#2+1-9); \vertsm{#1}{#2-9} \vertsm{#1+1}{#2+1-9}}
\nc{\Dedge}[2]{\draw[->-=0.6,line width=.3mm](#1,#2-9)--(#1+1,#2-1-9); \vertsm{#1}{#2-9} \vertsm{#1+1}{#2-1-9}}
\nc{\Fedge}[2]{\draw[->-=0.6,line width=.3mm](#1,#2-9)--(#1+1,#2-9); \vertsm{#1}{#2-9} \vertsm{#1+1}{#2-9}}
\nc{\Uedges}[1]{\foreach \x/\y in {#1}{\Uedge{\x}{\y}}}
\nc{\Dedges}[1]{\foreach \x/\y in {#1}{\Dedge{\x}{\y}}}
\nc{\Fedges}[1]{\foreach \x/\y in {#1}{\Fedge{\x}{\y}}}
\nc{\xvertlabel}[1]{\draw(#1,-10+.6)node{{\tiny $#1$}};}
\nc{\yvertlabel}[1]{\draw(0-.4,-9+#1)node{{\tiny $#1$}};}
\nc{\xvertlabels}[1]{\foreach \x in {#1}{ \xvertlabel{\x} }}
\nc{\yvertlabels}[1]{\foreach \x in {#1}{ \yvertlabel{\x} }}

\nc{\bbE}{\mathbb E}
\nc{\floorn}{\lfloor\tfrac n2\rfloor}
\rnc{\sp}{\supseteq}
\rnc{\arraystretch}{1.2}

\nc{\bn}{{[n]}} \nc{\bt}{{[t]}} \nc{\ba}{{[a]}} \nc{\bl}{{[l]}} \nc{\bm}{{[m]}} \nc{\bk}{{[k]}} \nc{\br}{{[r]}} \nc{\bs}{{[s]}} \nc{\bnf}{{[n-1]}}

\nc{\M}{\mathcal M}
\nc{\G}{\mathcal G}
\nc{\F}{\mathbb F}
\nc{\MnJ}{\mathcal M_n^J}
\nc{\EnJ}{\mathcal E_n^J}
\nc{\Mat}{\operatorname{Mat}}
\nc{\RegMnJ}{\Reg(\MnJ)}
\nc{\row}{\mathfrak r}
\nc{\col}{\mathfrak c}
\nc{\Row}{\operatorname{Row}}
\nc{\Col}{\operatorname{Col}}
\nc{\Span}{\operatorname{span}}
\nc{\mat}[4]{\left[\begin{matrix}#1&#2\\#3&#4\end{matrix}\right]}
\nc{\tmat}[4]{\left[\begin{smallmatrix}#1&#2\\#3&#4\end{smallmatrix}\right]}
\nc{\ttmat}[4]{{\tiny \left[\begin{smallmatrix}#1&#2\\#3&#4\end{smallmatrix}\right]}}
\nc{\tmatt}[9]{\left[\begin{smallmatrix}#1&#2&#3\\#4&#5&#6\\#7&#8&#9\end{smallmatrix}\right]}
\nc{\ttmatt}[9]{{\tiny \left[\begin{smallmatrix}#1&#2&#3\\#4&#5&#6\\#7&#8&#9\end{smallmatrix}\right]}}
\nc{\MnGn}{\M_n\sm\G_n}
\nc{\MrGr}{\M_r\sm\G_r}
\nc{\qbin}[2]{\left[\begin{matrix}#1\\#2\end{matrix}\right]_q}
\nc{\tqbin}[2]{\left[\begin{smallmatrix}#1\\#2\end{smallmatrix}\right]_q}
\nc{\qbinx}[3]{\left[\begin{matrix}#1\\#2\end{matrix}\right]_{#3}}
\nc{\tqbinx}[3]{\left[\begin{smallmatrix}#1\\#2\end{smallmatrix}\right]_{#3}}
\nc{\MNJ}{\M_nJ}
\nc{\JMN}{J\M_n}
\nc{\RegMNJ}{\Reg(\MNJ)}
\nc{\RegJMN}{\Reg(\JMN)}
\nc{\RegMMNJ}{\Reg(\MMNJ)}
\nc{\RegJMMN}{\Reg(\JMMN)}
\nc{\Wb}{\overline{W}}
\nc{\Xb}{\overline{X}}
\nc{\Yb}{\overline{Y}}
\nc{\Zb}{\overline{Z}}
\nc{\Sib}{\overline{\Si}}
\nc{\Om}{\Omega}
\nc{\Omb}{\overline{\Om}}
\nc{\Gab}{\overline{\Ga}}
\nc{\qfact}[1]{[#1]_q!}
\nc{\smat}[2]{\left[\begin{matrix}#1&#2\end{matrix}\right]}
\nc{\tsmat}[2]{\left[\begin{smallmatrix}#1&#2\end{smallmatrix}\right]}
\nc{\hmat}[2]{\left[\begin{matrix}#1\\#2\end{matrix}\right]}
\nc{\thmat}[2]{\left[\begin{smallmatrix}#1\\#2\end{smallmatrix}\right]}
\nc{\LVW}{\mathcal L(V,W)}
\nc{\KVW}{\mathcal K(V,W)}
\nc{\LV}{\mathcal L(V)}
\nc{\RegLVW}{\Reg(\LVW)}
\nc{\sM}{\mathscr M}
\nc{\sN}{\mathscr N}
\rnc{\iff}{\ \Leftrightarrow\ }
\nc{\Hom}{\operatorname{Hom}}
\nc{\End}{\operatorname{End}}
\nc{\Aut}{\operatorname{Aut}}
\nc{\Lin}{\mathcal L}
\nc{\Hommn}{\Hom(V_m,V_n)}
\nc{\Homnm}{\Hom(V_n,V_m)}
\nc{\Homnl}{\Hom(V_n,V_l)}
\nc{\Homkm}{\Hom(V_k,V_m)}
\nc{\Endm}{\End(V_m)}
\nc{\Endn}{\End(V_n)}
\nc{\Endr}{\End(V_r)}
\nc{\Autm}{\Aut(V_m)}
\nc{\Autn}{\Aut(V_n)}
\nc{\MmnJ}{\M_{mn}^J}
\nc{\MmnA}{\M_{mn}^A}
\nc{\MmnB}{\M_{mn}^B}
\nc{\Mmn}{\M_{mn}}
\nc{\Mkl}{\M_{kl}}
\nc{\Mnm}{\M_{nm}}
\nc{\EmnJ}{\mathcal E_{mn}^J}
\nc{\MmGm}{\M_m\sm\G_m}
\nc{\RegMmnJ}{\Reg(\MmnJ)}
\rnc{\implies}{\ \Rightarrow\ }
\nc{\DMmn}[1]{D_{#1}(\Mmn)}
\nc{\DMmnJ}[1]{D_{#1}(\MmnJ)}
\nc{\MMNJ}{\Mmn J}
\nc{\JMMN}{J\Mmn}
\nc{\JMMNJ}{J\Mmn J}
\nc{\Inr}{\mathcal I(V_n,W_r)}
\nc{\Lnr}{\mathcal L(V_n,W_r)}
\nc{\Knr}{\mathcal K(V_n,W_r)}
\nc{\Imr}{\mathcal I(V_m,W_r)}
\nc{\Kmr}{\mathcal K(V_m,W_r)}
\nc{\Lmr}{\mathcal L(V_m,W_r)}
\nc{\Kmmr}{\mathcal K(V_m,W_{m-r})}
\nc{\tr}{{\operatorname{T}}}
\nc{\MMN}{\MmnA(\F_1)}
\nc{\MKL}{\Mkl^B(\F_2)}
\nc{\RegMMN}{\Reg(\MmnA(\F_1))}
\nc{\RegMKL}{\Reg(\Mkl^B(\F_2))}
\nc{\gRhA}{\widehat{\mathscr R}^A}
\nc{\gRhB}{\widehat{\mathscr R}^B}
\nc{\gLhA}{\widehat{\mathscr L}^A}
\nc{\gLhB}{\widehat{\mathscr L}^B}
\nc{\timplies}{\Rightarrow}
\nc{\tiff}{\Leftrightarrow}
\nc{\Sija}{S_{ij}^a}
\nc{\dmat}[8]{\draw(#1*1.5,#2)node{$\left[\begin{smallmatrix}#3&#4&#5\\#6&#7&#8\end{smallmatrix}\right]$};}
\nc{\bdmat}[8]{\draw(#1*1.5,#2)node{${\mathbf{\left[\begin{smallmatrix}#3&#4&#5\\#6&#7&#8\end{smallmatrix}\right]}}$};}
\nc{\rdmat}[8]{\draw(#1*1.5,#2)node{\rotatebox{90}{$\left[\begin{smallmatrix}#3&#4&#5\\#6&#7&#8\end{smallmatrix}\right]$}};}
\nc{\rldmat}[8]{\draw(#1*1.5-0.375,#2)node{\rotatebox{90}{$\left[\begin{smallmatrix}#3&#4&#5\\#6&#7&#8\end{smallmatrix}\right]$}};}
\nc{\rrdmat}[8]{\draw(#1*1.5+.375,#2)node{\rotatebox{90}{$\left[\begin{smallmatrix}#3&#4&#5\\#6&#7&#8\end{smallmatrix}\right]$}};}
\nc{\rfldmat}[8]{\draw(#1*1.5-0.375+.15,#2)node{\rotatebox{90}{$\left[\begin{smallmatrix}#3&#4&#5\\#6&#7&#8\end{smallmatrix}\right]$}};}
\nc{\rfrdmat}[8]{\draw(#1*1.5+.375-.15,#2)node{\rotatebox{90}{$\left[\begin{smallmatrix}#3&#4&#5\\#6&#7&#8\end{smallmatrix}\right]$}};}
\nc{\xL}{[x]_{\! _\gL}}\nc{\yL}{[y]_{\! _\gL}}\nc{\xR}{[x]_{\! _\gR}}\nc{\yR}{[y]_{\! _\gR}}\nc{\xH}{[x]_{\! _\gH}}\nc{\yH}{[y]_{\! _\gH}}\nc{\XK}{[X]_{\! _\gK}}\nc{\xK}{[x]_{\! _\gK}}
\nc{\RegSija}{\Reg(\Sija)}
\nc{\MnmK}{\M_{nm}^K}
\nc{\cC}{\mathcal C}
\nc{\cR}{\mathcal R}
\nc{\Ckl}{\cC_k(l)}
\nc{\Rkl}{\cR_k(l)}
\nc{\Cmr}{\cC_m(r)}
\nc{\Rmr}{\cR_m(r)}
\nc{\Cnr}{\cC_n(r)}
\nc{\Rnr}{\cR_n(r)}
\nc{\Z}{\mathbb Z}

\nc{\Reg}{\operatorname{Reg}}
\nc{\RP}{\operatorname{RP}}
\nc{\TXa}{\T_X^a}
\nc{\TXA}{\T(X,A)}
\nc{\TXal}{\T(X,\al)}
\nc{\RegTXa}{\Reg(\TXa)}
\nc{\RegTXA}{\Reg(\TXA)}
\nc{\RegTXal}{\Reg(\TXal)}
\nc{\PalX}{\P_\al(X)}
\nc{\EAX}{\E_A(X)}
\nc{\Bb}{\overline{B}}
\nc{\bb}{\overline{\be}}
\nc{\bw}{{\bf w}}
\nc{\bz}{{\bf z}}
\nc{\TASA}{\T_A\sm\S_A}
\nc{\Ub}{\overline{U}}
\nc{\Vb}{\overline{V}}
\nc{\eb}{\overline{e}}
\nc{\EXa}{\E_X^a}
\nc{\oijr}{1\leq i<j\leq r}
\nc{\veb}{\overline{\ve}}
\nc{\bbT}{\mathbb T}
\nc{\Surj}{\operatorname{Surj}}
\nc{\Sone}{S^{(1)}}
\nc{\fillbox}[2]{\draw[fill=gray!30](#1,#2)--(#1+1,#2)--(#1+1,#2+1)--(#1,#2+1)--(#1,#2);}
\nc{\raa}{\rangle_J}
\nc{\raJ}{\rangle_J}
\nc{\Ea}{E_J}
\nc{\EJ}{E_J}
\nc{\ep}{\epsilon} \nc{\ve}{\varepsilon}
\nc{\IXa}{\I_X^a}
\nc{\RegIXa}{\Reg(\IXa)}
\nc{\JXa}{\J_X^a}
\nc{\RegJXa}{\Reg(\JXa)}
\nc{\IXA}{\I(X,A)}
\nc{\IAX}{\I(A,X)}
\nc{\RegIXA}{\Reg(\IXA)}
\nc{\RegIAX}{\Reg(\IAX)}
\nc{\trans}[2]{\left(\begin{smallmatrix} #1 \\ #2 \end{smallmatrix}\right)}
\nc{\bigtrans}[2]{\left(\begin{matrix} #1 \\ #2 \end{matrix}\right)}
\nc{\lmap}[1]{\mapstochar \xrightarrow {\ #1\ }}
\nc{\EaTXa}{E}

\nc{\gL}{\mathscr L}
\nc{\gR}{\mathscr R}
\nc{\gH}{\mathscr H}
\nc{\gJ}{\mathscr J}
\nc{\gD}{\mathscr D}
\nc{\gK}{\mathscr K}
\nc{\gLa}{\mathscr L^a}
\nc{\gRa}{\mathscr R^a}
\nc{\gHa}{\mathscr H^a}
\nc{\gJa}{\mathscr J^a}
\nc{\gDa}{\mathscr D^a}
\nc{\gKa}{\mathscr K^a}
\nc{\gLJ}{\mathscr L^J}
\nc{\gRJ}{\mathscr R^J}
\nc{\gHJ}{\mathscr H^J}
\nc{\gJJ}{\mathscr J^J}
\nc{\gDJ}{\mathscr D^J}
\nc{\gKJ}{\mathscr K^J}
\nc{\gLh}{\widehat{\mathscr L}^J}
\nc{\gRh}{\widehat{\mathscr R}^J}
\nc{\gHh}{\widehat{\mathscr H}^J}
\nc{\gJh}{\widehat{\mathscr J}^J}
\nc{\gDh}{\widehat{\mathscr D}^J}
\nc{\gKh}{\widehat{\mathscr K}^J}
\nc{\Lh}{\widehat{L}^J}
\nc{\Rh}{\widehat{R}^J}
\nc{\Hh}{\widehat{H}^J}
\nc{\Jh}{\widehat{J}^J}
\nc{\Dh}{\widehat{D}^J}
\nc{\Kh}{\widehat{K}^J}
\nc{\gLb}{\widehat{\mathscr L}}
\nc{\gRb}{\widehat{\mathscr R}}
\nc{\gHb}{\widehat{\mathscr H}}
\nc{\gJb}{\widehat{\mathscr J}}
\nc{\gDb}{\widehat{\mathscr D}}
\nc{\gKb}{\widehat{\mathscr K}}
\nc{\Lb}{\widehat{L}^J}
\nc{\Rb}{\widehat{R}^J}
\nc{\Hb}{\widehat{H}^J}
\nc{\Jb}{\widehat{J}^J}
\nc{\Db}{\overline{D}}
\nc{\Kb}{\widehat{K}}

\hyphenation{mon-oid mon-oids}

\nc{\itemit}[1]{\item[\emph{(#1)}]}
\nc{\E}{\mathcal E}
\nc{\TX}{\T(X)}
\nc{\TXP}{\T(X,\P)}
\nc{\EX}{\E(X)}
\nc{\EXP}{\E(X,\P)}
\nc{\SX}{\S(X)}
\nc{\SXP}{\S(X,\P)}
\nc{\Sing}{\operatorname{Sing}}
\nc{\idrank}{\operatorname{idrank}}
\nc{\SingXP}{\Sing(X,\P)}
\nc{\De}{\Delta}
\nc{\sgp}{\operatorname{sgp}}
\nc{\mon}{\operatorname{mon}}
\nc{\Dn}{\mathcal D_n}
\nc{\Dm}{\mathcal D_m}

\nc{\lline}[1]{\draw(3*#1,0)--(3*#1+2,0);}
\nc{\uline}[1]{\draw(3*#1,5)--(3*#1+2,5);}
\nc{\thickline}[2]{\draw(3*#1,5)--(3*#2,0); \draw(3*#1+2,5)--(3*#2+2,0) ;}
\nc{\thicklabel}[3]{\draw(3*#1+1+3*#2*0.15-3*#1*0.15,4.25)node{{\tiny $#3$}};}

\nc{\slline}[3]{\draw(3*#1+#3,0+#2)--(3*#1+2+#3,0+#2);}
\nc{\suline}[3]{\draw(3*#1+#3,5+#2)--(3*#1+2+#3,5+#2);}
\nc{\sthickline}[4]{\draw(3*#1+#4,5+#3)--(3*#2+#4,0+#3); \draw(3*#1+2+#4,5+#3)--(3*#2+2+#4,0+#3) ;}
\nc{\sthicklabel}[5]{\draw(3*#1+1+3*#2*0.15-3*#1*0.15+#5,4.25+#4)node{{\tiny $#3$}};}

\nc{\stll}[5]{\sthickline{#1}{#2}{#4}{#5} \sthicklabel{#1}{#2}{#3}{#4}{#5}}
\nc{\tll}[3]{\stll{#1}{#2}{#3}00}

\nc{\mfourpic}[9]{
\slline1{#9}0
\slline3{#9}0
\slline4{#9}0
\slline5{#9}0
\suline1{#9}0
\suline3{#9}0
\suline4{#9}0
\suline5{#9}0
\stll1{#1}{#5}{#9}{0}
\stll3{#2}{#6}{#9}{0}
\stll4{#3}{#7}{#9}{0}
\stll5{#4}{#8}{#9}{0}
\draw[dotted](6,0+#9)--(8,0+#9);
\draw[dotted](6,5+#9)--(8,5+#9);
}
\nc{\vdotted}[1]{
\draw[dotted](3*#1,10)--(3*#1,15);
\draw[dotted](3*#1+2,10)--(3*#1+2,15);
}

\nc{\Clab}[2]{
\sthicklabel{#1}{#1}{{}_{\phantom{#1}}C_{#1}}{1.25+5*#2}0
}
\nc{\sClab}[3]{
\sthicklabel{#1}{#1}{{}_{\phantom{#1}}C_{#1}}{1.25+5*#2}{#3}
}
\nc{\Clabl}[3]{
\sthicklabel{#1}{#1}{{}_{\phantom{#3}}C_{#3}}{1.25+5*#2}0
}
\nc{\sClabl}[4]{
\sthicklabel{#1}{#1}{{}_{\phantom{#4}}C_{#4}}{1.25+5*#2}{#3}
}
\nc{\Clabll}[3]{
\sthicklabel{#1}{#1}{C_{#3}}{1.25+5*#2}0
}
\nc{\sClabll}[4]{
\sthicklabel{#1}{#1}{C_{#3}}{1.25+5*#2}{#3}
}

\nc{\mtwopic}[6]{
\slline1{#6*5}{#5}
\slline2{#6*5}{#5}
\suline1{#6*5}{#5}
\suline2{#6*5}{#5}
\stll1{#1}{#3}{#6*5}{#5}
\stll2{#2}{#4}{#6*5}{#5}
}
\nc{\mtwopicl}[6]{
\slline1{#6*5}{#5}
\slline2{#6*5}{#5}
\suline1{#6*5}{#5}
\suline2{#6*5}{#5}
\stll1{#1}{#3}{#6*5}{#5}
\stll2{#2}{#4}{#6*5}{#5}
\sClabl1{#6}{#5}{i}
\sClabl2{#6}{#5}{j}
}

\nc{\keru}{\operatorname{ker}^\wedge} \nc{\kerl}{\operatorname{ker}_\vee}

\nc{\coker}{\operatorname{coker}}
\nc{\KER}{\ker}
\nc{\N}{\mathbb N}
\nc{\LaBn}{L_\al(\B_n)}
\nc{\RaBn}{R_\al(\B_n)}
\nc{\LaPBn}{L_\al(\PB_n)}
\nc{\RaPBn}{R_\al(\PB_n)}
\nc{\rhorBn}{\rho_r(\B_n)}
\nc{\DrBn}{D_r(\B_n)}
\nc{\DrPn}{D_r(\P_n)}
\nc{\DrPBn}{D_r(\PB_n)}
\nc{\DrKn}{D_r(\K_n)}
\nc{\alb}{\al_{\vee}}
\nc{\beb}{\be^{\wedge}}
\nc{\Bal}{\operatorname{Bal}}
\nc{\Red}{\operatorname{Red}}
\nc{\Pnxi}{\P_n^\xi}
\nc{\Bnxi}{\B_n^\xi}
\nc{\PBnxi}{\PB_n^\xi}
\nc{\Knxi}{\K_n^\xi}
\nc{\C}{\mathscr C}
\nc{\exi}{e^\xi}
\nc{\Exi}{E^\xi}
\nc{\eximu}{e^\xi_\mu}
\nc{\Eximu}{E^\xi_\mu}
\nc{\REF}{ {\red [Ref?]} }
\nc{\GL}{\operatorname{GL}}
\rnc{\O}{\mathcal O}

\nc{\vtx}[2]{\fill (#1,#2)circle(.2);}
\nc{\lvtx}[2]{\fill (#1,0)circle(.2);}
\nc{\uvtx}[2]{\fill (#1,1.5)circle(.2);}

\nc{\Eq}{\mathfrak{Eq}}
\nc{\Gau}{\Ga^\wedge} \nc{\Gal}{\Ga_\vee}
\nc{\Lamu}{\Lam^\wedge} \nc{\Laml}{\Lam_\vee}
\nc{\bX}{{\bf X}}
\nc{\bY}{{\bf Y}}
\nc{\ds}{\displaystyle}

\nc{\uuvert}[1]{\fill (#1,3)circle(.2);}
\nc{\uuuvert}[1]{\fill (#1,4.5)circle(.2);}
\nc{\overt}[1]{\fill (#1,0)circle(.1);}
\nc{\overtl}[3]{\node[vertex] (#3) at (#1,0) {  {\tiny $#2$} };}
\nc{\cv}[2]{\draw(#1,1.5) to [out=270,in=90] (#2,0);}
\nc{\cvs}[2]{\draw(#1,1.5) to [out=270+30,in=90+30] (#2,0);}
\nc{\ucv}[2]{\draw(#1,3) to [out=270,in=90] (#2,1.5);}
\nc{\uucv}[2]{\draw(#1,4.5) to [out=270,in=90] (#2,3);}
\nc{\textpartn}[1]{{\lower1.0 ex\hbox{\begin{tikzpicture}[xscale=.3,yscale=0.3] #1 \end{tikzpicture}}}}
\nc{\textpartnx}[2]{{\lower1.0 ex\hbox{\begin{tikzpicture}[xscale=.3,yscale=0.3] 
\foreach \x in {1,...,#1}
{ \uvert{\x} \lvert{\x} }
#2 \end{tikzpicture}}}}
\nc{\disppartnx}[2]{{\lower1.0 ex\hbox{\begin{tikzpicture}[scale=0.3] 
\foreach \x in {1,...,#1}
{ \uvert{\x} \lvert{\x} }
#2 \end{tikzpicture}}}}
\nc{\disppartnxd}[2]{{\lower2.1 ex\hbox{\begin{tikzpicture}[scale=0.3] 
\foreach \x in {1,...,#1}
{ \uuvert{\x} \uvert{\x} \lvert{\x} }
#2 \end{tikzpicture}}}}
\nc{\disppartnxdn}[2]{{\lower2.1 ex\hbox{\begin{tikzpicture}[scale=0.3] 
\foreach \x in {1,...,#1}
{ \uuvert{\x} \lvert{\x} }
#2 \end{tikzpicture}}}}
\nc{\disppartnxdd}[2]{{\lower3.6 ex\hbox{\begin{tikzpicture}[scale=0.3] 
\foreach \x in {1,...,#1}
{ \uuuvert{\x} \uuvert{\x} \uvert{\x} \lvert{\x} }
#2 \end{tikzpicture}}}}

\nc{\dispgax}[2]{{\lower0.0 ex\hbox{\begin{tikzpicture}[scale=0.3] 
#2
\foreach \x in {1,...,#1}
{\lvert{\x} }
 \end{tikzpicture}}}}
\nc{\textgax}[2]{{\lower0.4 ex\hbox{\begin{tikzpicture}[scale=0.3] 
#2
\foreach \x in {1,...,#1}
{\lvert{\x} }
 \end{tikzpicture}}}}
\nc{\textlinegraph}[2]{{\raise#1 ex\hbox{\begin{tikzpicture}[scale=0.8] 
#2
 \end{tikzpicture}}}}
\nc{\textlinegraphl}[2]{{\raise#1 ex\hbox{\begin{tikzpicture}[scale=0.8] 
\tikzstyle{vertex}=[circle,draw=black, fill=white, inner sep = 0.07cm]
#2
 \end{tikzpicture}}}}
\nc{\displinegraph}[1]{{\lower0.0 ex\hbox{\begin{tikzpicture}[scale=0.6] 
#1
 \end{tikzpicture}}}}
 
\nc{\disppartnthreeone}[1]{{\lower1.0 ex\hbox{\begin{tikzpicture}[scale=0.3] 
\foreach \x in {1,2,3,5,6}
{ \uvert{\x} }
\foreach \x in {1,2,4,5,6}
{ \lvert{\x} }
\draw[dotted] (3.5,1.5)--(4.5,1.5);
\draw[dotted] (2.5,0)--(3.5,0);
#1 \end{tikzpicture}}}}

\nc{\partn}[4]{\left( \begin{array}{c|c} 
#1 \ & \ #3 \ \ \\ \cline{2-2}
#2 \ & \ #4 \ \
\end{array} \!\!\! \right)}
\nc{\partnlong}[6]{\partn{#1}{#2}{#3,\ #4}{#5,\ #6}} 
\nc{\partnsh}[2]{\left( \begin{array}{c} 
#1 \\
#2 
\end{array} \right)}
\nc{\partncodefz}[3]{\partn{#1}{#2}{#3}{\emptyset}}
\nc{\partndefz}[3]{{\partn{#1}{#2}{\emptyset}{#3}}}
\nc{\partnlast}[2]{\left( \begin{array}{c|c}
#1 \ &  \ #2 \\
#1 \ &  \ #2
\end{array} \right)}

\nc{\uuarcx}[3]{\draw(#1,3)arc(180:270:#3) (#1+#3,3-#3)--(#2-#3,3-#3) (#2-#3,3-#3) arc(270:360:#3);}
\nc{\uuarc}[2]{\uuarcx{#1}{#2}{.4}}
\nc{\uuuarcx}[3]{\draw(#1,4.5)arc(180:270:#3) (#1+#3,4.5-#3)--(#2-#3,4.5-#3) (#2-#3,4.5-#3) arc(270:360:#3);}
\nc{\uuuarc}[2]{\uuuarcx{#1}{#2}{.4}}
\nc{\udarcx}[3]{\draw(#1,1.5)arc(180:90:#3) (#1+#3,1.5+#3)--(#2-#3,1.5+#3) (#2-#3,1.5+#3) arc(90:0:#3);}
\nc{\udarc}[2]{\udarcx{#1}{#2}{.4}}
\nc{\uudarcx}[3]{\draw(#1,3)arc(180:90:#3) (#1+#3,3+#3)--(#2-#3,3+#3) (#2-#3,3+#3) arc(90:0:#3);}
\nc{\uudarc}[2]{\uudarcx{#1}{#2}{.4}}
\nc{\darcxhalf}[3]{\draw(#1,0)arc(180:90:#3) (#1+#3,#3)--(#2,#3) ;}
\nc{\darchalf}[2]{\darcxhalf{#1}{#2}{.4}}
\nc{\uarcxhalf}[3]{\draw(#1,1.5)arc(180:270:#3) (#1+#3,1.5-#3)--(#2,1.5-#3) ;}
\nc{\uarchalf}[2]{\uarcxhalf{#1}{#2}{.4}}
\nc{\uarcxhalfr}[3]{\draw (#1+#3,1.5-#3)--(#2-#3,1.5-#3) (#2-#3,1.5-#3) arc(270:360:#3);}
\nc{\uarchalfr}[2]{\uarcxhalfr{#1}{#2}{.4}}

\nc{\bdarcx}[3]{\draw[blue](#1,0)arc(180:90:#3) (#1+#3,#3)--(#2-#3,#3) (#2-#3,#3) arc(90:0:#3);}
\nc{\bdarc}[2]{\darcx{#1}{#2}{.4}}
\nc{\rduarcx}[3]{\draw[red](#1,0)arc(180:270:#3) (#1+#3,0-#3)--(#2-#3,0-#3) (#2-#3,0-#3) arc(270:360:#3);}
\nc{\rduarc}[2]{\uarcx{#1}{#2}{.4}}
\nc{\duarcx}[3]{\draw(#1,0)arc(180:270:#3) (#1+#3,0-#3)--(#2-#3,0-#3) (#2-#3,0-#3) arc(270:360:#3);}
\nc{\duarc}[2]{\uarcx{#1}{#2}{.4}}

\nc{\uuv}[1]{\fill (#1,4)circle(.1);}
\nc{\uv}[1]{\fill (#1,2)circle(.1);}
\nc{\lv}[1]{\fill (#1,0)circle(.1);}
\nc{\uvred}[1]{\fill[red] (#1,2)circle(.1);}
\nc{\lvred}[1]{\fill[red] (#1,0)circle(.1);}
\nc{\lvwhite}[1]{\fill[white] (#1,0)circle(.1);}
\nc{\buv}[1]{\fill (#1,2)circle(.18);}
\nc{\blv}[1]{\fill (#1,0)circle(.18);}

\nc{\uvs}[1]{{
\foreach \x in {#1}
{ \uv{\x}}
}}
\nc{\uuvs}[1]{{
\foreach \x in {#1}
{ \uuv{\x}}
}}
\nc{\lvs}[1]{{
\foreach \x in {#1}
{ \lv{\x}}
}}

\nc{\buvs}[1]{{
\foreach \x in {#1}
{ \buv{\x}}
}}
\nc{\blvs}[1]{{
\foreach \x in {#1}
{ \blv{\x}}
}}

\nc{\uvreds}[1]{{
\foreach \x in {#1}
{ \uvred{\x}}
}}
\nc{\lvreds}[1]{{
\foreach \x in {#1}
{ \lvred{\x}}
}}

\nc{\uudotted}[2]{\draw [dotted] (#1,4)--(#2,4);}
\nc{\uudotteds}[1]{{
\foreach \x/\y in {#1}
{ \uudotted{\x}{\y}}
}}
\nc{\uudottedsm}[2]{\draw [dotted] (#1+.4,4)--(#2-.4,4);}
\nc{\uudottedsms}[1]{{
\foreach \x/\y in {#1}
{ \uudottedsm{\x}{\y}}
}}
\nc{\udottedsm}[2]{\draw [dotted] (#1+.4,2)--(#2-.4,2);}
\nc{\udottedsms}[1]{{
\foreach \x/\y in {#1}
{ \udottedsm{\x}{\y}}
}}
\nc{\udotted}[2]{\draw [dotted] (#1,2)--(#2,2);}
\nc{\udotteds}[1]{{
\foreach \x/\y in {#1}
{ \udotted{\x}{\y}}
}}
\nc{\ldotted}[2]{\draw [dotted] (#1,0)--(#2,0);}
\nc{\ldotteds}[1]{{
\foreach \x/\y in {#1}
{ \ldotted{\x}{\y}}
}}
\nc{\ldottedsm}[2]{\draw [dotted] (#1+.4,0)--(#2-.4,0);}
\nc{\ldottedsms}[1]{{
\foreach \x/\y in {#1}
{ \ldottedsm{\x}{\y}}
}}

\nc{\stlinest}[2]{\draw(#1,4)--(#2,0);}

\nc{\stlined}[2]{\draw[dotted](#1,2)--(#2,0);}

\nc{\tlab}[2]{\draw(#1,2)node[above]{\tiny $#2$};}
\nc{\tudots}[1]{\draw(#1,2)node{$\cdots$};}
\nc{\tldots}[1]{\draw(#1,0)node{$\cdots$};}

\nc{\uvw}[1]{\fill[white] (#1,2)circle(.1);}
\nc{\huv}[1]{\fill (#1,1)circle(.1);}
\nc{\llv}[1]{\fill (#1,-2)circle(.1);}
\nc{\arcup}[2]{
\draw(#1,2)arc(180:270:.4) (#1+.4,1.6)--(#2-.4,1.6) (#2-.4,1.6) arc(270:360:.4);
}
\nc{\harcup}[2]{
\draw(#1,1)arc(180:270:.4) (#1+.4,.6)--(#2-.4,.6) (#2-.4,.6) arc(270:360:.4);
}
\nc{\arcdn}[2]{
\draw(#1,0)arc(180:90:.4) (#1+.4,.4)--(#2-.4,.4) (#2-.4,.4) arc(90:0:.4);
}
\nc{\cve}[2]{
\draw(#1,2) to [out=270,in=90] (#2,0);
}
\nc{\hcve}[2]{
\draw(#1,1) to [out=270,in=90] (#2,0);
}
\nc{\catarc}[3]{
\draw(#1,2)arc(180:270:#3) (#1+#3,2-#3)--(#2-#3,2-#3) (#2-#3,2-#3) arc(270:360:#3);
}

\nc{\arcr}[2]{
\draw[red](#1,0)arc(180:90:.4) (#1+.4,.4)--(#2-.4,.4) (#2-.4,.4) arc(90:0:.4);
}
\nc{\arcb}[2]{
\draw[blue](#1,2-2)arc(180:270:.4) (#1+.4,1.6-2)--(#2-.4,1.6-2) (#2-.4,1.6-2) arc(270:360:.4);
}
\nc{\loopr}[1]{
\draw[blue](#1,-2) edge [out=130,in=50,loop] ();
}
\nc{\loopb}[1]{
\draw[red](#1,-2) edge [out=180+130,in=180+50,loop] ();
}
\nc{\redto}[2]{\draw[red](#1,0)--(#2,0);}
\nc{\bluto}[2]{\draw[blue](#1,0)--(#2,0);}
\nc{\dotto}[2]{\draw[dotted](#1,0)--(#2,0);}
\nc{\lloopr}[2]{\draw[red](#1,0)arc(0:360:#2);}
\nc{\lloopb}[2]{\draw[blue](#1,0)arc(0:360:#2);}
\nc{\rloopr}[2]{\draw[red](#1,0)arc(-180:180:#2);}
\nc{\rloopb}[2]{\draw[blue](#1,0)arc(-180:180:#2);}
\nc{\uloopr}[2]{\draw[red](#1,0)arc(-270:270:#2);}
\nc{\uloopb}[2]{\draw[blue](#1,0)arc(-270:270:#2);}
\nc{\dloopr}[2]{\draw[red](#1,0)arc(-90:270:#2);}
\nc{\dloopb}[2]{\draw[blue](#1,0)arc(-90:270:#2);}
\nc{\llloopr}[2]{\draw[red](#1,0-2)arc(0:360:#2);}
\nc{\llloopb}[2]{\draw[blue](#1,0-2)arc(0:360:#2);}
\nc{\lrloopr}[2]{\draw[red](#1,0-2)arc(-180:180:#2);}
\nc{\lrloopb}[2]{\draw[blue](#1,0-2)arc(-180:180:#2);}
\nc{\ldloopr}[2]{\draw[red](#1,0-2)arc(-270:270:#2);}
\nc{\ldloopb}[2]{\draw[blue](#1,0-2)arc(-270:270:#2);}
\nc{\luloopr}[2]{\draw[red](#1,0-2)arc(-90:270:#2);}
\nc{\luloopb}[2]{\draw[blue](#1,0-2)arc(-90:270:#2);}

\nc{\larcb}[2]{
\draw[blue](#1,0-2)arc(180:90:.4) (#1+.4,.4-2)--(#2-.4,.4-2) (#2-.4,.4-2) arc(90:0:.4);
}
\nc{\larcr}[2]{
\draw[red](#1,2-2-2)arc(180:270:.4) (#1+.4,1.6-2-2)--(#2-.4,1.6-2-2) (#2-.4,1.6-2-2) arc(270:360:.4);
}

\rnc{\H}{\mathscr H}
\rnc{\L}{\mathscr L}
\nc{\R}{\mathscr R}
\nc{\D}{\mathscr D}
\nc{\J}{\mathscr J}

\nc{\ssim}{\mathrel{\raise0.25 ex\hbox{\oalign{$\approx$\crcr\noalign{\kern-0.84 ex}$\approx$}}}}
\nc{\POI}{\mathcal{O}}
\nc{\wb}{\overline{w}}
\nc{\ub}{\overline{u}}
\nc{\vb}{\overline{v}}
\nc{\fb}{\overline{f}}
\nc{\gb}{\overline{g}}
\nc{\hb}{\overline{h}}
\nc{\pb}{\overline{p}}
\rnc{\sb}{\overline{s}}
\nc{\XR}{\pres{X}{R\,}}
\nc{\YQ}{\pres{Y}{Q}}
\nc{\ZP}{\pres{Z}{P\,}}
\nc{\XRone}{\pres{X_1}{R_1}}
\nc{\XRtwo}{\pres{X_2}{R_2}}
\nc{\XRthree}{\pres{X_1\cup X_2}{R_1\cup R_2\cup R_3}}
\nc{\er}{\eqref}
\nc{\larr}{\mathrel{\hspace{-0.35 ex}>\hspace{-1.1ex}-}\hspace{-0.35 ex}}
\nc{\rarr}{\mathrel{\hspace{-0.35 ex}-\hspace{-0.5ex}-\hspace{-2.3ex}>\hspace{-0.35 ex}}}
\nc{\lrarr}{\mathrel{\hspace{-0.35 ex}>\hspace{-1.1ex}-\hspace{-0.5ex}-\hspace{-2.3ex}>\hspace{-0.35 ex}}}
\nc{\nn}{\tag*{}}
\nc{\epfal}{\tag*{$\Box$}}
\nc{\tagd}[1]{\tag*{(#1)$'$}}
\nc{\ldb}{[\![}
\nc{\rdb}{]\!]}
\nc{\sm}{\setminus}
\nc{\I}{\mathcal I}
\nc{\InSn}{\I_n\setminus\S_n}
\nc{\dom}{\operatorname{dom}} \nc{\codom}{\operatorname{codom}}
\nc{\ojin}{1\leq j<i\leq n}
\nc{\eh}{\widehat{e}}
\nc{\wh}{\widehat{w}}
\nc{\uh}{\widehat{u}}
\nc{\vh}{\widehat{v}}
\nc{\sh}{\widehat{s}}
\nc{\fh}{\widehat{f}}
\nc{\textres}[1]{\text{\emph{#1}}}
\nc{\aand}{\emph{\ and \ }}
\nc{\iif}{\emph{\ if \ }}
\nc{\textlarr}{\ \larr\ }
\nc{\textrarr}{\ \rarr\ }
\nc{\textlrarr}{\ \lrarr\ }

\nc{\comma}{,\ }

\nc{\COMMA}{,\quad}
\nc{\TnSn}{\T_n\setminus\S_n} 
\nc{\TmSm}{\T_m\setminus\S_m} 
\nc{\TXSX}{\T_X\setminus\S_X} 
\rnc{\S}{\mathcal S}

\nc{\T}{\mathcal T} 
\nc{\A}{\mathscr A} 
\nc{\B}{\mathcal B} 
\rnc{\P}{\mathcal P} 
\nc{\K}{\mathcal K}
\nc{\PB}{\mathcal{PB}} 
\nc{\rank}{\operatorname{rank}}

\nc{\mtt}{\!\!\!\mt\!\!\!}

\nc{\sub}{\subseteq}
\nc{\la}{\langle}
\nc{\ra}{\rangle}
\nc{\mt}{\mapsto}
\nc{\im}{\mathrm{im}}
\nc{\id}{\mathrm{id}}
\nc{\al}{\alpha}
\nc{\be}{\beta}
\nc{\ga}{\gamma}
\nc{\Ga}{\Gamma}
\nc{\de}{\delta}
\nc{\ka}{\kappa}
\nc{\lam}{\lambda}
\nc{\Lam}{\Lambda}
\nc{\si}{\sigma}
\nc{\Si}{\Sigma}
\nc{\oijn}{1\leq i<j\leq n}
\nc{\oijm}{1\leq i<j\leq m}

\nc{\comm}{\rightleftharpoons}
\nc{\AND}{\qquad\text{and}\qquad}

\nc{\bit}{\vspace{-3 truemm}\begin{itemize}}
\nc{\bitbmc}{\begin{itemize}\begin{multicols}}
\nc{\bmc}{\begin{itemize}\begin{multicols}}
\nc{\emc}{\end{multicols}\end{itemize}\vspace{-3 truemm}}
\nc{\eit}{\end{itemize}\vspace{-3 truemm}}
\nc{\ben}{\vspace{-3 truemm}\begin{enumerate}}
\nc{\een}{\end{enumerate}\vspace{-3 truemm}}
\nc{\eitres}{\end{itemize}}

\nc{\set}[2]{\{ {#1} : {#2} \}} 
\nc{\bigset}[2]{\big\{ {#1}: {#2} \big\}} 
\nc{\Bigset}[2]{\left\{ \,{#1} :{#2}\, \right\}}

\nc{\pres}[2]{\la {#1} \,|\, {#2} \ra}
\nc{\bigpres}[2]{\big\la {#1} \,\big|\, {#2} \big\ra}
\nc{\Bigpres}[2]{\Big\la \,{#1}\, \,\Big|\, \,{#2}\, \Big\ra}
\nc{\Biggpres}[2]{\Bigg\la {#1} \,\Bigg|\, {#2} \Bigg\ra}

\nc{\pf}{\noindent{\bf Proof.}  }
\nc{\epf}{\hfill$\Box$\bigskip}
\nc{\epfres}{\hfill$\Box$}
\nc{\pfnb}{\pf}
\nc{\epfnb}{\bigskip}
\nc{\pfthm}[1]{\bigskip \noindent{\bf Proof of Theorem \ref{#1}}\,\,  } 
\nc{\pfprop}[1]{\bigskip \noindent{\bf Proof of Proposition \ref{#1}}\,\,  } 
\nc{\epfreseq}{\tag*{$\Box$}}

\nc{\uvert}[1]{\fill (#1,2)circle(.2);}
\rnc{\lvert}[1]{\fill (#1,0)circle(.2);}
\nc{\guvert}[1]{\fill[lightgray] (#1,2)circle(.2);}
\nc{\glvert}[1]{\fill[lightgray] (#1,0)circle(.2);}
\nc{\uvertx}[2]{\fill (#1,#2)circle(.2);}
\nc{\guvertx}[2]{\fill[lightgray] (#1,#2)circle(.2);}
\nc{\uvertxs}[2]{
\foreach \x in {#1}
{ \uvertx{\x}{#2}  }
}
\nc{\guvertxs}[2]{
\foreach \x in {#1}
{ \guvertx{\x}{#2}  }
}

\nc{\uvertth}[2]{\fill (#1,2)circle(#2);}
\nc{\lvertth}[2]{\fill (#1,0)circle(#2);}
\nc{\uvertths}[2]{
\foreach \x in {#1}
{ \uvertth{\x}{#2}  }
}
\nc{\lvertths}[2]{
\foreach \x in {#1}
{ \lvertth{\x}{#2}  }
}

\nc{\vertlabel}[2]{\draw(#1,2+.3)node{{\tiny $#2$}};}
\nc{\vertlabelh}[2]{\draw(#1,2+.4)node{{\tiny $#2$}};}
\nc{\vertlabelhh}[2]{\draw(#1,2+.6)node{{\tiny $#2$}};}
\nc{\vertlabelhhh}[2]{\draw(#1,2+.64)node{{\tiny $#2$}};}
\nc{\vertlabelup}[2]{\draw(#1,4+.6)node{{\tiny $#2$}};}
\nc{\vertlabels}[1]{
{\foreach \x/\y in {#1}
{ \vertlabel{\x}{\y} }
}
}

\nc{\dvertlabel}[2]{\draw(#1,-.4)node{{\tiny $#2$}};}
\nc{\dvertlabels}[1]{
{\foreach \x/\y in {#1}
{ \dvertlabel{\x}{\y} }
}
}
\nc{\vertlabelsh}[1]{
{\foreach \x/\y in {#1}
{ \vertlabelh{\x}{\y} }
}
}
\nc{\vertlabelshh}[1]{
{\foreach \x/\y in {#1}
{ \vertlabelhh{\x}{\y} }
}
}
\nc{\vertlabelshhh}[1]{
{\foreach \x/\y in {#1}
{ \vertlabelhhh{\x}{\y} }
}
}

\nc{\vertlabelx}[3]{\draw(#1,2+#3+.6)node{{\tiny $#2$}};}
\nc{\vertlabelxs}[2]{
{\foreach \x/\y in {#1}
{ \vertlabelx{\x}{\y}{#2} }
}
}

\nc{\vertlabelupdash}[2]{\draw(#1,2.7)node{{\tiny $\phantom{'}#2'$}};}
\nc{\vertlabelupdashess}[1]{
{\foreach \x/\y in {#1}
{\vertlabelupdash{\x}{\y}}
}
}

\nc{\vertlabeldn}[2]{\draw(#1,0-.6)node{{\tiny $\phantom{'}#2'$}};}
\nc{\vertlabeldnph}[2]{\draw(#1,0-.6)node{{\tiny $\phantom{'#2'}$}};}

\nc{\vertlabelups}[1]{
{\foreach \x in {#1}
{\vertlabel{\x}{\x}}
}
}
\nc{\vertlabeldns}[1]{
{\foreach \x in {#1}
{\vertlabeldn{\x}{\x}}
}
}
\nc{\vertlabeldnsph}[1]{
{\foreach \x in {#1}
{\vertlabeldnph{\x}{\x}}
}
}

\nc{\dotsup}[2]{\draw [dotted] (#1+.6,2)--(#2-.6,2);}
\nc{\dotsupx}[3]{\draw [dotted] (#1+.6,#3)--(#2-.6,#3);}
\nc{\dotsdn}[2]{\draw [dotted] (#1+.6,0)--(#2-.6,0);}
\nc{\dotsups}[1]{\foreach \x/\y in {#1}
{ \dotsup{\x}{\y} }
}
\nc{\dotsupxs}[2]{\foreach \x/\y in {#1}
{ \dotsupx{\x}{\y}{#2} }
}
\nc{\dotsdns}[1]{\foreach \x/\y in {#1}
{ \dotsdn{\x}{\y} }
}

\nc{\nodropcustpartn}[3]{
\begin{tikzpicture}[scale=.3]
\foreach \x in {#1}
{ \uvert{\x}  }
\foreach \x in {#2}
{ \lvert{\x}  }
#3 \end{tikzpicture}
}

\nc{\custpartn}[3]{{\lower1.4 ex\hbox{
\begin{tikzpicture}[scale=.3]
\foreach \x in {#1}
{ \uvert{\x}  }
\foreach \x in {#2}
{ \lvert{\x}  }
#3 \end{tikzpicture}
}}}

\nc{\smcustpartn}[3]{{\lower0.7 ex\hbox{
\begin{tikzpicture}[scale=.2]
\foreach \x in {#1}
{ \uvert{\x}  }
\foreach \x in {#2}
{ \lvert{\x}  }
#3 \end{tikzpicture}
}}}

\nc{\dropcustpartn}[3]{{\lower5.2 ex\hbox{
\begin{tikzpicture}[scale=.3]
\foreach \x in {#1}
{ \uvert{\x}  }
\foreach \x in {#2}
{ \lvert{\x}  }
#3 \end{tikzpicture}
}}}

\nc{\dropcustpartnx}[4]{{\lower#4 ex\hbox{
\begin{tikzpicture}[scale=.4]
\foreach \x in {#1}
{ \uvert{\x}  }
\foreach \x in {#2}
{ \lvert{\x}  }
#3 \end{tikzpicture}
}}}

\nc{\dropcustpartnxy}[3]{\dropcustpartnx{#1}{#2}{#3}{4.6}}

\nc{\uvertsm}[1]{\fill (#1,2)circle(.15);}
\nc{\lvertsm}[1]{\fill (#1,0)circle(.15);}
\nc{\vertsm}[2]{\fill (#1,#2)circle(.15);}

\nc{\bigdropcustpartn}[3]{{\lower6.93 ex\hbox{
\begin{tikzpicture}[scale=.6]
\foreach \x in {#1}
{ \uvertsm{\x}  }
\foreach \x in {#2}
{ \lvertsm{\x}  }
#3 \end{tikzpicture}
}}}

\nc{\gcustpartn}[5]{{\lower1.4 ex\hbox{
\begin{tikzpicture}[scale=.3]
\foreach \x in {#1}
{ \uvert{\x}  }
\foreach \x in {#2}
{ \guvert{\x}  }
\foreach \x in {#3}
{ \lvert{\x}  }
\foreach \x in {#4}
{ \glvert{\x}  }
#5 \end{tikzpicture}
}}}

\nc{\gcustpartndash}[5]{{\lower6.97 ex\hbox{
\begin{tikzpicture}[scale=.3]
\foreach \x in {#1}
{ \uvert{\x}  }
\foreach \x in {#2}
{ \guvert{\x}  }
\foreach \x in {#3}
{ \lvert{\x}  }
\foreach \x in {#4}
{ \glvert{\x}  }
#5 \end{tikzpicture}
}}}

\nc{\stline}[2]{\draw(#1,2)--(#2,0);}
\nc{\stlines}[1]{
{\foreach \x/\y in {#1}
{ \stline{\x}{\y} }
}
}

\nc{\uarcs}[1]{
{\foreach \x/\y in {#1}
{ \uarc{\x}{\y} }
}
}

\nc{\darcs}[1]{
{\foreach \x/\y in {#1}
{ \darc{\x}{\y} }
}
}

\nc{\stlinests}[1]{
{\foreach \x/\y in {#1}
{ \stlinest{\x}{\y} }
}
}

\nc{\stlineds}[1]{
{\foreach \x/\y in {#1}
{ \stlined{\x}{\y} }
}
}

\nc{\gstline}[2]{\draw[lightgray](#1,2)--(#2,0);}
\nc{\gstlines}[1]{
{\foreach \x/\y in {#1}
{ \gstline{\x}{\y} }
}
}

\nc{\gstlinex}[3]{\draw[lightgray](#1,2+#3)--(#2,0+#3);}
\nc{\gstlinexs}[2]{
{\foreach \x/\y in {#1}
{ \gstlinex{\x}{\y}{#2} }
}
}

\nc{\stlinex}[3]{\draw(#1,2+#3)--(#2,0+#3);}
\nc{\stlinexs}[2]{
{\foreach \x/\y in {#1}
{ \stlinex{\x}{\y}{#2} }
}
}

\nc{\darcx}[3]{\draw(#1,0)arc(180:90:#3) (#1+#3,#3)--(#2-#3,#3) (#2-#3,#3) arc(90:0:#3);}
\nc{\darc}[2]{\darcx{#1}{#2}{.4}}
\nc{\uarcx}[3]{\draw(#1,2)arc(180:270:#3) (#1+#3,2-#3)--(#2-#3,2-#3) (#2-#3,2-#3) arc(270:360:#3);}
\nc{\uarc}[2]{\uarcx{#1}{#2}{.4}}

\nc{\darcxx}[4]{\draw(#1,0+#4)arc(180:90:#3) (#1+#3,#3+#4)--(#2-#3,#3+#4) (#2-#3,#3+#4) arc(90:0:#3);}
\nc{\uarcxx}[4]{\draw(#1,2+#4)arc(180:270:#3) (#1+#3,2-#3+#4)--(#2-#3,2-#3+#4) (#2-#3,2-#3+#4) arc(270:360:#3);}

\makeatletter
\newcommand\footnoteref[1]{\protected@xdef\@thefnmark{\ref{#1}}\@footnotemark}
\makeatother

\newcounter{theorem}
\numberwithin{theorem}{section}

\newtheorem{thm}[theorem]{Theorem}
\newtheorem{lemma}[theorem]{Lemma}
\newtheorem{cor}[theorem]{Corollary}
\newtheorem{prop}[theorem]{Proposition}

\theoremstyle{definition}

\newtheorem{rem}[theorem]{Remark}
\newtheorem{defn}[theorem]{Definition}
\newtheorem{eg}[theorem]{Example}
\newtheorem{ass}[theorem]{Assumption}

\title{Motzkin monoids and partial Brauer monoids}

\date{}

\author{
Igor Dolinka%
\footnote{Department of Mathematics and Informatics, University of Novi Sad, Trg Dositeja Obradovi\'ca 4, 21101 Novi Sad, Serbia, {\tt dockie\,@\,dmi.uns.ac.rs}.}%
, \
James East%
\footnote{Centre for Research in Mathematics; School of Computing, Engineering and Mathematics, Western Sydney University, Locked Bag 1797, Penrith NSW 2751, Australia, {\tt J.East\,@\,WesternSydney.edu.au}.}%
, \
Robert D.~Gray%
\footnote{School of Mathematics, University of East Anglia, Norwich NR4 7TJ , United Kingdom, {\tt Robert.D.Gray\,@\,uea.ac.uk}.}
}

\maketitle

\vspace{-0.5cm}

\begin{abstract}
We study the partial Brauer monoid and its planar submonoid, the Motzkin monoid.  We conduct a thorough investigation of the structure of both monoids, providing information on normal forms, Green's relations, regularity, ideals, idempotent generation, minimal (idempotent) generating sets, and so on.  We obtain necessary and sufficient conditions under which the ideals of these monoids are idempotent-generated.  We find formulae for the rank (smallest size of a generating set) of each ideal, and for the idempotent rank (smallest size of an idempotent generating set) of the idempotent-generated subsemigroup of each ideal; in particular, when an ideal is idempotent-generated, the rank and idempotent rank are equal.  Along the way, we obtain a number of results of independent interest, and we demonstrate the utility of the semigroup theoretic approach by applying our results to obtain new proofs of some important representation theoretic results concerning the corresponding diagram algebras, the partial (or rook) Brauer algebra and Motzkin algebra.



{\it Keywords}: Partial Brauer monoids, Motzkin monoids, Brauer monoids, idempotents, ideals, rank, idempotent rank, diagram algebras, partial Brauer algebras, Motzkin algebras, cellular algebras.

MSC: 20M20; 20M17, 05A18, 05E15, 16S36.
\end{abstract}




\section{Introduction}\label{sect:intro}

There are many reasons to study the monoids of the title.  On the one hand, \emph{diagram algebras} (including Brauer algebras \cite{Brauer1937}, partial Brauer algebras \cite{HD2014,MarMaz2014}, Temperley-Lieb algebras \cite{TL1971}, Motzkin algebras \cite{BH2014}, rook monoid algebras \cite{HR2001}, partition algebras \cite{Martin1994,Jones1994_2}, etc.)~are ubiquitous in representation theory and statistical mechanics, and several recent approaches to diagram algebras via \emph{diagram monoids} and \emph{twisted semigroup algebras} have proved extremely successful \cite{Wilcox2007,JEgrpm, JEpnsn,DEEFHHL1,DEEFHHL2,HR2005,BDP2002,LF2006}.  On the other hand, diagram monoids are of direct interest to semigroup theorists for a range of reasons.  For one thing, diagram monoids are closely related to several important transformation semigroups; indeed, the partition monoids contain isomorphic copies of the full transformation semigroup and the symmetric inverse monoid \cite{FL2011,EF,JEgrpm}.  For another, diagram monoids provide natural examples of \emph{regular $*$-semigroups} \cite{NS1978}, a variety of semigroups contained strictly between the varieties of regular and inverse semigroups (see Section \ref{sect:prelim} for definitions); thanks to the $*$-regular structure of the partition monoids, it is now known that any semigroup embeds in an idempotent-generated regular $*$-semigroup \cite{EF}, and that any finite semigroup embeds in a $2$-generator regular $*$-semigroup \cite{JEpme}.  Diagram monoids have also played a large role in the development of the theory of pseudovarieties of finite semigroups \cite{ADV2012_2,Auinger2014,Auinger2012,ACHLV}.  Finally, and of particular significance to the current work, algebraic studies of diagram monoids have led to incredibly rich combinatorial structures \cite{DEEFHHL1,DEEFHHL2,EastGray}.

The current article continues in the combinatorial theme of a previous paper of the second and third named authors~\cite{EastGray}, which, in turn, took its inspiration from a number of foundational papers of John Howie on combinatorial aspects of finite full transformation semigroups.  In \cite{Howie1966}, Howie showed that the singular ideal of a finite full transformation semigroup is generated by its idempotents.  In subsequent work, Howie and various collaborators calculated the \emph{rank} (minimal size of a generating set) and \emph{idempotent rank} (minimal size of an idempotent generating set) 
of this singular ideal, classified the minimal (in size) idempotent generating sets, and extended these results to arbitrary ideals~\cite{Howie1978,Gomes1987,HLF,Howie1990}; a prominent role was played by certain well-known number sequences, such as binomial coefficients and Stirling numbers.  Analogous results have been obtained for many other important families of semigroups, such as full linear monoids \cite{Gray2008,Erdos1967,Fountain1991}, endomorphism monoids of various algebraic structures \cite{Fountain1992,Fountain1993,Gray2007,DE1,DEM2015}, sandwich semigroups \cite{DE2,DE3} and, more recently, certain families of diagram monoids \cite{JEpnsn,EF,Maltcev2007,EastGray}. 

The current article mainly concerns the \emph{partial Brauer monoid} $\PB_n$ and the \emph{Motzkin monoid} $\M_n$ (see Section \ref{sect:prelim} for definitions).  There is a growing body of literature on partial Brauer monoids and algebras;%
\footnote{Partial Brauer monoids and algebras are also known in the literature as \emph{rook Brauer} monoids and algebras.} 
see for example \cite{Maz2013,HR2015,Maz2002,HD2014,MarMaz2014,Maz1998,KM2006,DEEFHHL1,ADV2012_2}.  Motzkin monoids and algebras are a more recent phenomenon \cite{Maz2013,HR2015,ADV2012_2,DEEFHHL2,BH2014,PHY2013}, and are \emph{planar} versions of the partial Brauer monoids and algebras (see Section \ref{sect:prelim}); the relationship between $\PB_n$ and $\M_n$ is therefore analogous to that between the \emph{Brauer monoid} $\B_n$ and the \emph{Jones monoid} $\mathcal J_n$.%
\footnote{Jones monoids are occasionally referred to as \emph{Temperley-Lieb} monoids; see \cite{LF2006} for a discussion of naming conventions.  Motzkin monoids and algebras have also been called \emph{partial Jones} (or \emph{partial Temperley-Lieb}) monoids and algebras; we prefer the \emph{Motzkin} nomenclature, since the processes of partialising diagrams and of selecting planar elements do not commute: as noted in \cite{DEEFHHL2}, the ``planar partial Brauer monoid'' \emph{strictly contains} the ``partial planar Brauer monoid''.}  
The planarity built into the definition of the Motzkin and Jones monoids leads to an \emph{aperiodic} structure (all subgroups are trivial), and can have interesting effects when comparing the complexity of certain problems; for example, the enumeration of idempotents is far more difficult in $\M_n$ and $\mathcal J_n$ than in $\PB_n$ and $\B_n$ \cite{DEEFHHL1,DEEFHHL2}, while it is quite the opposite for the enumeration of minimal-size idempotent generating sets of the singular ideals of $\mathcal J_n$ and~$\B_n$~\cite{EastGray}.
The situation in the current article is more akin to the former; although the difficulty is fairly comparable, the combinatorial analysis is more subtle for $\M_n$ than for $\PB_n$, and relies on less-familiar number sequences such as the Motzkin and Riordan numbers.
In the spirit of the articles mentioned in the previous paragraph, we mostly focus on questions concerning ideals, idempotent-generation, and minimal (idempotent) generating sets.  But we also demonstrate the utility of the semigroup theoretic approach by applying our results to the corresponding diagram algebras in order to recover important (known) results on representation theory and cellularity.

Before summarising the results of the article, we first note that some quite general methods exist for studying generating sets and (idempotent) ranks, of so-called (finite) \emph{completely $0$-simple semigroups}. By the much-celebrated Rees theorem~\cite[Theorem 3.2.3]{Howie}, each completely $0$-simple semigroup is isomorphic to a \emph{Rees $0$-matrix semigroup} over a group. Given such a Rees $0$-matrix semigroup, $S = M^0[G;I,\Lambda;P]$, one can associate a finite bipartite graph, with vertex set $I \cup \Lambda$ and with edges (labelled by elements from the group~$G$) corresponding to the non-zero entries in the structure matrix $P$. This is called the \emph{Graham--Houghton graph} of $S$ and is a useful tool for studying questions concerning subsemigroups and generators of~$S$. In particular, this approach was used in \cite{Gray2005} to give a formula for the rank of an arbitrary finite completely $0$-simple semigroup. These techniques were later developed and extended to investigate idempotent ranks in \cite{Gray2008,Gray2007}. These results are relevant for the general problem of determining the rank of a (not necessarily completely $0$-simple) finite semigroup $S$. Indeed, if $S$ is a finite semigroup, then the rank of $S$ is bounded below by the sum of the ranks of the principal factors of the maximal $\gJ$-classes of $S$. (Green's relations, which include the $\gJ$-relation, are defined in Section \ref{sect:prelim}.)  There are in fact many natural examples where this bound is attained, and so the question of the rank of the semigroup reduces to computing the rank of certain finite completely $0$-simple semigroups. This is the case, for example, for the proper two-sided ideals of the full transformation monoid and the full linear monoid, as well as various diagram monoids such as the partition, Brauer and Jones monoids; see \cite{GrayHowieIssuePaper,EastGray,Howie1990} for more details and for additional examples and references. In particular, in all of the diagram monoids considered in \cite{EastGray}, the proper ideals all have a \emph{unique} maximal $\gJ$-class whose elements (indeed, idempotents) generate the whole ideal. For this reason, the questions of rank and idempotent rank in that work reduced to questions about completely $0$-simple semigroups (and certain combinatorial data, including the above-mentioned Graham--Houghton graphs). The current paper represents a significant departure from this pattern: most of the ideals we consider are \emph{not} generated by the elements of a single $\gJ$-class; even though each proper ideal still does have a unique maximal $\gJ$-class, any generating set must contain elements from a non-maximal $\gJ$-class.  As a result, 
the question of computing (idempotent) ranks cannot be resolved simply by applying general results about ranks of completely $0$-simple semigroups, hence the necessity of new methods, such as those we employ here. 

The article is organised as follows.  In Section \ref{sect:prelim}, we gather the preliminary results we will need concerning basic semigroup theory, regular $*$-semigroups and the various diagram monoids that will play a role in our investigations.  The main result of Section \ref{sect:prelim} (Proposition \ref{prop:DClasses}) is that $\PB_n$ and $\M_n$ both have $n+1$ ideals:
\[
I_0(\PB_n)\sub I_1(\PB_n)\sub\cdots\sub I_n(\PB_n)=\PB_n \AND I_0(\M_n)\sub I_1(\M_n)\sub\cdots\sub I_n(\M_n)=\M_n.
\]
We also prove some enumerative results concerning Green's relations in both monoids (Propositions~\ref{prop:DClassesofKn},~\ref{prop:DClassesofMn} and~\ref{prop:m'nr}).  
Section \ref{sect:PBn} is devoted to the partial Brauer monoid $\PB_n$.  
A proper ideal $I_r(\PB_n)\subsetneq\PB_n$ is generated by its top $\D$-class if $r=0$ or $r\equiv n\pmod2$; otherwise, $I_r(\PB_n)$ is generated by its top \emph{two} $\D$-classes (Corollaries \ref{cor:IdDrDr-1} and \ref{cor:r0nmod2} and Proposition \ref{prop:In-1}).  We calculate the rank of each ideal $I_r(\PB_n)$, and the idempotent rank in the case that the ideal is idempotent-generated; this occurs when $0\leq r\leq n-2$, and the rank and idempotent rank are equal for such values of $r$ (Theorem \ref{thm:IrPBn} and Propositions \ref{prop:In-1} and \ref{prop:rankPBn}).  We also describe the idempotent-generated subsemigroup of $\PB_n$, and calculate its rank and idempotent rank, which are again equal (Theorem \ref{thm:EPBn}).
Some of the general structural results we prove in this section are of independent interest; for example, the normal forms we provide in Lemma \ref{lem:normalformPBn} (see also Remark \ref{rem:normalformPBn}) are used in later sections to obtain information about the corresponding diagram algebras.
The situation is somewhat different for the Motzkin monoid~$\M_n$, which is the subject of Section \ref{sect:Mn}.  
Each proper ideal $I_r(\M_n)\subsetneq\M_n$ is still generated by its top two $\D$-classes, but only by a single $\D$-class if $r=0$ (Propositions~\ref{prop:IdDrDr-1M} and~\ref{prop:multinesting}).  While an ideal $I_r(\M_n)$ is only idempotent generated if $r<\floorn$, we are still able to calculate the rank of each ideal, as well as the idempotent rank when applicable; again, the rank and idempotent rank are equal (Proposition \ref{prop:IrIGM} and Theorems~\ref{thm:IrM} and~\ref{thm:rankMn}).  We also describe the idempotent-generated subsemigroup of each ideal, and calculate its rank and idempotent rank (Theorems \ref{thm:IGIrM} and \ref{thm:EMn}); when $r\leq n-2$, these values are again equal and, curiously, also equal to the rank of the corresponding ideal, even in the case that the ideal is not idempotent-generated (i.e., when $\floorn\leq r\leq n-2$).  
In Section \ref{sect:algebras}, we apply results of previous sections to deduce (known) information about the associated diagram algebras, the so-called \emph{partial Brauer algebras} and \emph{Motzkin algebras}; specifically, we show that both algebras are cellular, and calculate the dimensions of the cell modules in the case that the ground ring is a field (these coincide with the irreducible modules in the semisimple case).
Finally, Section \ref{sect:values} contains a number of tables displaying calculated values of various ranks and idempotent ranks.

For non-negative integers $a,b$, we write $[a,b]=\{a,a+1,\ldots,b\}$, which we assume to be empty if $b<a$.  We also write $[n]=[1,n]=\{1,\ldots,n\}$, noting that $[0]=\emptyset$.  If $A\sub\bn$, we write $A^c=\bn\sm A$.  For a non-negative integer $k$, we write $k!!=0$ if $k$ is even or $k!!=k\cdot(k-2)\cdots3\cdot1$ if $k$ is odd.  By convention, we define $(-1)!!=1$.  We also interpret a binomial coefficient ${n\choose k}$ to be $0$ if $k<0$ or $n<k$.  We will sometimes write $X=\{i_1<\cdots<i_k\}$ to indicate that $X=\{i_1,\ldots,i_k\}$ and $i_1<\cdots<i_k$.


\section{Preliminaries}\label{sect:prelim}

In this section, we record some preliminary results concerning regular $*$-semigroups and the various diagram monoids that will play a role in our investigations.
But first we recall some basic notions from semigroup theory.  For more information, the reader is refered to a text such as \cite{Hig,Howie,RSbook}.

Let $S$ be a semigroup.  If $X\sub S$, we write $\la X\ra$ for the subsemigroup of $S$ generated by $X$.  The \emph{rank} of $S$, denoted $\rank(S)$, is the minimal cardinality of a subset $X\sub S$ such that $S=\la X\ra$.  Recall from \cite{HRH} that the \emph{relative rank} of $S$ with respect to a subset $A\sub S$, denoted $\rank(S:A)$, is the minimum cardinality of a subset $X\sub S$ such that $S=\la X\cup A\ra$.  

For any subset $X\sub S$, write $E(X)=\set{x\in X}{x=x^2}$ for the set of all idempotents contained in $X$.  If $S$ is idempotent-generated, the \emph{idempotent rank} of $S$, denoted $\idrank(S)$, is the minimal cardinality of a subset $X\sub E(S)$ such that $S=\la X\ra$; note that $\rank(S)\leq\idrank(S)$.  We write $\bbE(S)=\la E(S)\ra$ for the idempotent-generated subsemigroup of $S$.

Recall that \emph{Green's relations} $\gR,\gL,\gJ,\gH,\gD$ are defined, for $x,y\in S$, by
\begin{gather*}
x\gR y \iff xS^1=yS^1 \COMMA x\gL y \iff S^1x=S^1y \COMMA x\gJ y \iff S^1xS^1 =S^1yS^1  , \\
\gH=\gR\cap\gL \COMMA \gD=\gR\circ\gL=\gL\circ\gR.
\end{gather*}
Here, $S^1$ denotes the monoid obtained by adjoining an identity $1$ to $S$ (if necessary).
If $x\in S$, and if $\gK$ is one of Green's relations, we denote by $K_x$ the $\gK$-class of $x$ in~$S$.  An $\gH$-class contains an idempotent if and only if it is a group, in which case it is a maximal subgroup of $S$.  If $e$ and $f$ are $\D$-related idempotents of $S$, then the subgroups $H_e$ and $H_f$ are isomorphic.  If $S$ is a monoid, then the $\H$-class of the identity element of $S$ is the group of units of $S$.  The $\gJ$-classes of $S$ are partially ordered; we say that $J_x\leq J_y$ if $x\in S^1yS^1$.  If $S$ is finite, then $\gJ=\gD$.  

An element $x\in S$ is \emph{regular} if $x=xyx$ and $y=yxy$ for some $y\in S$ or, equivalently, if $D_x$ contains an idempotent, in which case $R_x$ and $L_x$ do, too.  We say~$S$ is regular if every element of $S$ is regular.  Recall that a semigroup $S$ is \emph{stable} if for any $x,y\in S$:
\[
xy\D x \iff xy\R x \AND xy\D y \iff xy\L y.
\]
In particular, all finite semigroups are stable.  If $T$ is a subsemigroup of $S$, then we write $\gK^T$ for Green's $\gK$-relation on~$T$.  So, for example, if $x,y\in T$, then $x\R^T y \iff xT^1=yT^1$.  If $T$ is a regular subsemigroup of~$S$ (which may or may not be regular itself), and if $\gK$ is one of $\R,\L,\H$, then $\gK^T$ is simply the restriction of $\gK$ to $T$, but this is not necessarily true if $\gK$ is one of $\D$ or $\J$.  We will occasionally write $\gK^S$ for Green's $\gK$-relation on $S$ if it is necessary to distinguish it from the $\gK$-relation on a subsemigroup.

\subsection{Regular $*$-semigroups}

Recall that a semigroup $S$ is a \emph{regular $*$-semigroup} \cite{NS1978} if there is a unary operation ${}^*:S\to S:x\mt x^*$ satisfying
\[
x^{**}=x \COMMA (xy)^*=y^*x^* \COMMA xx^*x=x   \qquad\text{for all $x,y\in S$.}
\]
(The additional identity $x^*xx^*=x^*$ follows as a simple consequence.)  Recall that an element $x$ of a regular $*$-semigroup $S$ is a \emph{projection} if $x^*=x=x^2$.  We denote the set of all projections of $S$ by $P(S)$, and we note that $P(S)\sub E(S)$.  More generally, if $X$ is any subset of $S$, we write $P(X)=P(S)\cap X$.  
The next result is well known; see for example \cite{YS1989}.


\ms
\begin{prop}\label{prop:RSS}
Let $S$ be a regular $*$-semigroup. 
\ben
\itemit{i} We have $P(S)=\{aa^*:a\in S\}=\{a^*a:a\in S\}$.
\itemit{ii} We have $E(S)=P(S)^2$.  In particular, $\la E(S)\ra=\la P(S)\ra$,
and $S$ is idempotent-generated if and only if it is projection-generated.
\itemit{iii} Every $\R$-class of $S$ contains precisely one projection, as does every $\L$-class.  In particular, the number of $\R$-classes (equivalently, $\L$-classes) contained in a $\D$-class $D$ is equal to $|P(D)|$. \epfres
\een
\end{prop}


To further elaborate on part (iii) of the above result, if $x$ is an element of a regular $*$-semigroup $S$, then $xx^*\R x\L x^*x$, with $xx^*,x^*x\in P(D_x)$.  Note also that $x\D x^*$ for all $x\in S$.  



\ms
\begin{lemma}\label{lem:ARSS<D>}
Let $D$ be a $\D$-class of a finite regular $*$-semigroup $S$, and let $T=\la D\ra$ and $P=P(D)$.  
\bit
\itemit{i} We have $\rank(T)\geq|P|$.
\itemit{ii} If the $\H$-classes contained in $D$ have size $1$, then $\rank(T)=|P|$.
\itemit{iii} If $T$ is idempotent-generated, then $T=\la P\ra$ and $\idrank(T)=\rank(T)=|P|$.
\eit
\end{lemma}

\pf Since $T\sm D$ is an ideal of $T$, it is clear that $\rank(T)$ is equal to the minimal cardinality of a subset $A\sub D$ for which $D\sub\la A\ra$.  

For (i), suppose $A\sub D$ is such that $D\sub\la A\ra$.  To show that $\rank(T)\geq|P|$, it suffices to show that $|A|\geq|P|$.  Let $p\in P$, and consider an expression $p=a_1\cdots a_k$ where $a_1,\ldots,a_k\in A$.  In particular, since $p\in D$, we have $a_1 \D p = a_1(a_2\cdots a_k)$.  By stability, since $S$ is finite, it follows that $a_1 \R a_1(a_2\cdots a_k)=p$.  In particular, $A$ contains at least one element from the $\R$-class of $p$.  Since projections belong to distinct $\R$-classes, it follows that $|A|\geq|P|$, as required.  

Next, suppose the $\H$-classes contained in $D$ all have size $1$.  To prove (ii), it suffices to find a generating set for $T$ of size $|P|$.  Now, for any $p,q\in P$, the intersection $R_p\cap L_q$ is an $\H$-class in $D$, so we may write $R_p\cap L_q=\{a_{pq}\}$.  Since $\R$- and $\L$-classes contain unique projections, we have $a_{pq}a_{pq}^*=p$ and $a_{pq}^*a_{pq}=q$.  We first claim that for any $p,q,r,s\in P$, $a_{pq}a_{rs}\in D$ if and only if $qr\in D$, in which case $a_{pq}a_{rs}=a_{ps}$.  Indeed, 
\begin{align*}
a_{pq}a_{rs} \in D \implies a_{pq}a_{rs} \D a_{pq} 
&\implies a_{pq}a_{rs} \R a_{pq} &&\text{by stability}\\
&\implies qa_{rs}=a_{pq}^*a_{pq}a_{rs} \R a_{pq}^*a_{pq}=q &&\text{since $\R$ is a left congruence}\\
&\implies qa_{rs} \D q \D a_{rs} \\
&\implies qa_{rs} \L a_{rs} &&\text{by stability}\\
&\implies qr=qa_{rs}a_{rs}^* \L a_{rs}a_{rs}^*=r &&\text{since $\L$ is a right congruence}\\
&\implies qr\D r
\implies qr\in D.
\end{align*}
The reverse implication is proved in similar fashion.  But then, by stability again,
\[
a_{pq}a_{rs}\in D 
\implies 
\begin{cases}
a_{pq}a_{rs} \D a_{pq} \implies a_{pq}a_{rs} \R a_{pq} \R p\\
a_{pq}a_{rs} \D a_{rs} \implies a_{pq}a_{rs} \L a_{rs} \L s,
\end{cases}
\]
so that $a_{pq}a_{rs}\in R_p\cap L_s=\{a_{ps}\}$, establishing the claim.
Now choose an arbitrary ordering $p_1,\ldots,p_k$ of the elements of $P$, and put
$A=\{a_{p_1p_2},\ldots,a_{p_{k-1}p_k},a_{p_kp_1}\}$, noting that $|A|=|P|$.
Then for any $i,j\in\bk$, we have
\[
a_{p_ip_j} = a_{p_ip_{i+1}}a_{p_{i+1}p_{i+2}}\cdots a_{p_{j-1}p_j},
\]
where we may have to read the subscripts modulo $k$, showing that $D\sub\la A\ra$, and completing the proof of~(ii).

Finally, if $T$ is idempotent-generated, then $T=\la E(D)\ra=\la P\ra$, so that $P$ is an idempotent generating set of size $\rank(T)$, giving $\idrank(T)\leq|P|$.  Since $\rank(T)\leq\idrank(T)$, (iii) is proved, after applying (i). \epf

\subsection{Diagram monoids}

Let $n$ be a non-negative integer, and write $\PB_n$ for the set of all set partitions of $\bn\cup\bn'=\{1,\ldots,n\}\cup\{1',\ldots,n'\}$ into blocks of size at most $2$.  For example, here are two elements of $\PB_6$:
\begin{align*}
\al &= \big\{
\{1,3\}, \{2,3'\}, \{5,6\}, \{4',5'\},
\{4\}, \{1'\}, \{2'\}, \{6'\}
\big\},
\\
\be &= \big\{
\{1,3\}, \{4,3'\}, \{5,6\}, \{4',5'\},
\{2\}, \{1'\}, \{2'\}, \{6'\}
\big\}.
\end{align*}
There is a unique element of $\PB_0$, namely the empty partition.  

If we write $a(m)$ for the number of ways to partition an $m$-element set into blocks of size at most $2$, then one easily establishes the recurrence
\[
a(0)=a(1)=1 \COMMA a(m)=a(m-1)+(m-1)a(m-2) \quad\text{for $m\geq2$.}
\]
In particular, $|\PB_n|=a(2n)$.  An alternative formula, $|\PB_n|=\sum_{k=0}^n{2n\choose2k}(2k-1)!!$, was given in \cite{Maz1998,HD2014}.  The numbers $a(m)$ and $a(2m)$ appear as sequences A000085 and A066223, respectively, on the OEIS \cite{OEIS}.  


An element of $\PB_n$ may be represented (uniquely) by a graph on vertex set $\bn\cup\bn'$; a single edge is included between vertices $u,v\in\bn\cup\bn'$ if and only if $\{u,v\}$ is a block of $\al$ of size $2$.  We typically identify $\al\in\PB_n$ with its corresponding graph.  When drawing such a graph, the vertices $1,\ldots,n$ are arranged in a horizontal line (increasing from left to right) with vertex $i'$ directly below $i$ for each $i\in\bn$, and the edges are always drawn within the interior of the rectangle defined by the vertices.  A graph drawn in such a way is called a \emph{Brauer $n$-diagram}, or simply a \emph{Brauer diagram} if $n$ is understood from context.  So, with $\al,\be\in\PB_6$ as above, we have:
\[
\al = \custpartn{1,2,3,4,5,6}{1,2,3,4,5,6}{\uarcx13{.5}\uarc56\darc45\stline23}
\AND
\be = \custpartn{1,2,3,4,5,6}{1,2,3,4,5,6}{\uarcx13{.5}\uarc56\darc45\stline43}.
\]
Occasionally it will be convenient to order the top and/or bottom vertices differently, but this will always be made clear; see Figure~\ref{fig:DsDr}.

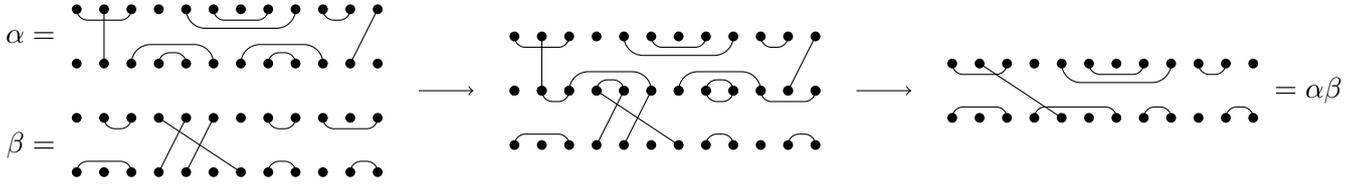
\begin{figure}[H]
   \begin{center}
\begin{tikzpicture}[scale=.36]
\begin{scope}[shift={(0,0)}]	
\buvs{1,...,12}
\blvs{1,...,12}
\stlines{2/2,12/11}
\uarcs{1/3,6/8,10/11}
\uarcx59{.7}
\darcs{4/5,8/9}
\darcx36{.7}
\darcx7{10}{.7}
\draw(0.6,1)node[left]{$\al=$};
\draw[->](13.5,-1)--(15.5,-1);
\end{scope}
\begin{scope}[shift={(0,-4)}]	
\buvs{1,...,12}
\blvs{1,...,12}
\stlines{4/7,5/4,6/5}
\uarcs{2/3,8/9,10/12}
\darcs{1/3,8/9,11/12}
\draw(0.6,1)node[left]{$\be=$};
\end{scope}
\begin{scope}[shift={(16,-1)}]	
\buvs{1,...,12}
\blvs{1,...,12}
\stlines{2/2,12/11}
\uarcs{1/3,6/8,10/11}
\uarcx59{.7}
\darcs{4/5,8/9}
\darcx36{.7}
\darcx7{10}{.7}
\draw[->](13.5,0)--(15.5,0);
\end{scope}
\begin{scope}[shift={(16,-3)}]	
\buvs{1,...,12}
\blvs{1,...,12}
\stlines{4/7,5/4,6/5}
\uarcs{2/3,8/9,10/12}
\darcs{1/3,8/9,11/12}
\end{scope}
\begin{scope}[shift={(32,-2)}]	
\buvs{1,...,12}
\blvs{1,...,12}
\stlines{2/5}
\uarcs{1/3,6/8,10/11}
\uarcx59{.7}
\darcs{1/3,4/7,8/9,11/12}
\draw(12.4,1)node[right]{$=\al\be$};
\end{scope}
\end{tikzpicture}
    \caption{Two Brauer diagrams $\al,\be\in\PB_{12}$ (left), their product $\al\be\in\PB_{12}$ (right), and the product graph $\Ga(\al,\be)$ (centre).}
    \label{fig:multinPB8}
   \end{center}
 \end{figure}

The set $\PB_n$ forms a monoid, known as the \emph{partial Brauer monoid}, under an operation we now describe.  Let $\al,\be\in\PB_n$.  Write $\bn''=\{1'',\ldots,n''\}$.  Let $\alb$ be the graph obtained from $\al$ by changing the label of each lower vertex $i'$ to~$i''$.  Similarly, let $\beb$ be the graph obtained from $\be$ by changing the label of each upper vertex~$i$ to~$i''$.  Consider now the graph $\Ga(\al,\be)$ on the vertex set~$\bn\cup \bn'\cup \bn''$ obtained by joining $\alb$ and~$\beb$ together so that each lower vertex $i''$ of $\alb$ is identified with the corresponding upper vertex $i''$ of $\beb$.  Note that $\Ga(\al,\be)$, which we call the \emph{product graph}, may contain parallel edges.  We define $\al\be\in\PB_n$ to be the Brauer diagram that has an edge $\{x,y\}$ if and only if $x,y\in\bn\cup\bn'$ are connected by a path in $\Ga(\al,\be)$.  An example calculation (with $n=12$) is given in Figure~\ref{fig:multinPB8}.

The partial Brauer monoid is a regular $*$-semigroup.  Indeed, for $\al\in\PB_n$, we write $\al^*$ for the Brauer diagram obtained from $\al$ by interchanging dashed and undashed vertices (i.e., by reflecting $\al$ in a horizontal axis).  For example, with $\al\in\PB_6$ as above, we have
\[
\al = \custpartn{1,2,3,4,5,6}{1,2,3,4,5,6}{\uarcx13{.5}\uarc56\darc45\stline23}
\quad\mt\quad
\al^* = \custpartn{1,2,3,4,5,6}{1,2,3,4,5,6}{\darcx13{.5}\darc56\uarc45\stline32}.
\]
We note that $\PB_n$ is a $*$-submonoid of the larger \emph{partition monoid} $\P_n$; we will not need to discuss $\P_n$ any further, but the reader is referred, for example, to \cite{HR2005,EastGray}.  We now describe a number of important $*$-submonoids of $\PB_n$.
\bit
\item We write $\M_n$ for the \emph{Motzkin monoid}, which is the submonoid of $\PB_n$ consisting of all Brauer diagrams that may be drawn in planar fashion (with vertices arranged as above, and with all edges within the rectangle defined by the vertices); such a planar Brauer diagram is called a \emph{Motzkin diagram}.  
\item We write $\B_n$ for the \emph{Brauer monoid}, which is the submonoid of $\PB_n$ consisting of all Brauer diagrams in which each block is of size $2$.  
\item We call a block of $\al\in\PB_n$ an \emph{upper hook} if it has size $2$ and is contained in $\bn$; \emph{lower hooks} are defined analogously.  We write $\I_n$ for the submonoid of $\PB_n$ consisting of all Brauer diagrams with no hooks.  So $\I_n$ is (isomorphic to) the \emph{symmetric inverse monoid} (also known as the \emph{rook monoid}), which consists of all injective partial maps $\bn\to\bn$.
\eit
For example, with $\al,\be\in\PB_6$ as above, we have $\be\in\M_6$ but $\al\not\in\M_6$.  We also have $\ga\in\B_6$ and $\de\in\I_6$, where:
\[
\ga =  \custpartn{1,2,3,4,5,6}{1,2,3,4,5,6}{\stline13\stline42\uarc23\uarc56\darcx16{.8}\darc45}
\AND
\de = \custpartn{1,2,3,4,5,6}{1,2,3,4,5,6}{\stline13\stline42\stline55}.
\]
Various intersections of the above submonoids are also of importance:
\bit
\item We write $\S_n=\I_n\cap\B_n$ for the \emph{symmetric group}.
\item We write $\mathcal J_n=\B_n\cap\M_n$ for the \emph{Jones monoid} (also known as the \emph{Temperley-Lieb monoid}).
\item We write $\O_n=\M_n\cap\I_n$ for the monoid of all order-preserving injective partial maps $\bn\to\bn$.  (Recall that $\al\in\I_n$ is \emph{order-preserving} if $i\al<j\al$ whenever $i,j$ belong to the domain of $\al$ and $i<j$.  Also, note that $\O_n$ is often denoted $\mathcal{POI}_n$ in the literature.)
\eit
The intersection of all the above monoids is the trivial submonoid, $\{1\}$.  Here,
\[
1=\big\{\{1,1'\},\{2,2'\},\ldots,\{n,n'\}\big\} = \custpartn{1,2,5}{1,2,5}{\dotsups{2/5}\dotsdns{2/5}\stlines{1/1,2/2,5/5}} \in\PB_n
\]
denotes the identity Brauer diagram.  The relevant part of the submonoid lattice of $\PB_n$ is displayed in Figure~\ref{fig:submonoids}.


\begin{figure}[h]
\begin{center}
\begin{tikzpicture}[scale=1]
%
\node (A) at (3,6) {$\PB_n$};
\node (B) at (0,4) {$\B_n$};
\node (C) at (3,4) {$\M_n$};
\node (D) at (6,4) {$\I_n$};
\node (E) at (0,2) {$\mathcal J_n$};
\node (F) at (3,2) {$\S_n$};
\node (G) at (6,2) {$\POI_n$};
\node (H) at (3,0) {$\{1\}$};
\draw (A)--(B) (A)--(C) (A)--(D);
\draw (C)--(E) (C)--(G);
\draw[white,line width=2mm] (B)--(F) (D)--(F);
\draw (B)--(E) (B)--(F);
\draw (D)--(F) (D)--(G);
\draw (E)--(H) (F)--(H) (G)--(H);
\end{tikzpicture}
\qquad\qquad\qquad
\begin{tikzpicture}[scale=1]
%
\node (A) at (3,6) {$\custpartn{1,2,3,4,5,6}{1,2,3,4,5,6}{\uarcx13{.5}\uarc56\darc45\stline23}$};
\node (B) at (0,4) {$\custpartn{1,2,3,4,5,6}{1,2,3,4,5,6}{\stline13\stline42\uarc23\uarc56\darcx16{.8}\darc45}$};
\node (C) at (3,4) {$\custpartn{1,2,3,4,5,6}{1,2,3,4,5,6}{\uarcx13{.5}\uarc56\darc45\stline43}$};
\node (D) at (6,4) {$\custpartn{1,2,3,4,5,6}{1,2,3,4,5,6}{\stline13\stline42\stline55}$};
\node (E) at (0,2) {$\custpartn{1,2,3,4,5,6}{1,2,3,4,5,6}{\uarcx12{.4}\uarc45\darc34\darcx25{.8}\stlines{3/1,6/6}}$};
\node (F) at (3,2) {$\custpartn{1,2,3,4,5,6}{1,2,3,4,5,6}{\stlines{1/2,2/1,3/4,4/6,5/5,6/3}}$};
\node (G) at (6,2) {$\custpartn{1,2,3,4,5,6}{1,2,3,4,5,6}{\stlines{1/1,2/3,6/4}}$};
\node (H) at (3,0) {$\custpartn{1,2,3,4,5,6}{1,2,3,4,5,6}{\stlines{1/1,2/2,3/3,4/4,5/5,6/6}}$};
\draw (A)--(B) (A)--(C) (A)--(D);
\draw (C)--(E) (C)--(G);
\draw[white,line width=2mm] (B)--(F) (D)--(F);
\draw (B)--(E) (B)--(F);
\draw (D)--(F) (D)--(G);
\draw (E)--(H) (F)--(H) (G)--(H);
\end{tikzpicture}
\end{center}
\vspace{-5mm}
\caption{Important submonoids of $\PB_n$ (left) and representative elements from each submonoid (right).  See text for further explanation.}
\label{fig:submonoids}
\end{figure}

Green's relations on $\PB_n$ (and all the submonoids mentioned above) may be conveniently described in terms of a number of parameters we now describe.  With this in mind, let $\al\in\PB_n$.  A block of $\al$ is called a \emph{transversal block} if it has nonempty intersection with both $\bn$ and $\bn'$, and a \emph{nontransversal block} otherwise.  The \emph{rank} of $\al$, denoted $\rank(\al)$, is equal to the number of transversal blocks of $\al$.  For $x\in\bn\cup\bn'$, let $[x]_\al$ denote the block of $\al$ containing $x$.  We define the \emph{domain} and \emph{codomain} of $\al$ to be the sets
\begin{align*}
\dom(\al) = \bigset{ x\in\bn } { [x]_\al\cap\bn'\not=\emptyset} &\AND
\codom(\al) = \bigset{ x\in\bn } { [x']_\al\cap\bn\not=\emptyset}.
\intertext{We also define the \emph{kernel} and \emph{cokernel} of $\al$ to be the equivalences}
\ker(\al) = \bigset{(x,y)\in\bn\times\bn}{[x]_\al=[y]_\al} &\AND
\coker(\al) = \bigset{(x,y)\in\bn\times\bn}{[x']_\al=[y']_\al}.
\end{align*}
To illustrate these ideas, let
\[
\ve = \custpartn{1,2,3,4,5,6,7,8}{1,2,3,4,5,6,7,8}
{
\stlinexs{4/6,7/5,8/8}{0}
\darcxx12{.4}{0}
\darcxx34{.4}{0}
}
\in\PB_8.
\]
Then $\rank(\ve)=3$, $\dom(\ve)=\{4,7,8\}$, $\codom(\ve)=\{5,6,8\}$, and $\ve$ has non-trivial cokernel classes $\{1,2\}$ and $\{3,4\}$, but no non-trivial kernel classes.

It is immediate from the definitions that
\[
\begin{array}{rclcrcl}
\dom(\al\be) \hspace{-.25cm}&\sub&\hspace{-.25cm} \dom(\al), & &
\ker(\al\be)\hspace{-.25cm} &\sp&\hspace{-.25cm} \ker(\al),\\
\codom(\al\be) \hspace{-.25cm}&\sub&\hspace{-.25cm} \codom(\be), & &
\coker(\al\be)\hspace{-.25cm} &\sp&\hspace{-.25cm} \coker(\be),
\end{array}
\]
for all $\al,\be\in\PB_n$.  Also, any upper non-transversal block of $\al$ is an upper non-transversal block of $\al\be$ for any $\al,\be\in\PB_n$, with a similar statement holding for lower non-transversal blocks of $\be$.  

We now recall another way to specify an element of $\PB_n$.  With this in mind, let $\al\in\PB_n$.  We write
\[
\al = \left( \begin{array}{c|c|c|c|c|c} 
i_1 \ & \ \cdots \ & \ i_r \ & \ A_1 \ & \ \cdots \ & \ A_s\ \ \\ \cline{4-6}
j_1 \ & \ \cdots \ & \ j_r \ & \ B_1 \ & \ \cdots \ & \ B_t \ \
\end{array} \!\!\! \right)
=
\left( \begin{array}{c|c} 
i_k \ & \ A_l \ \ \\ \cline{2-2}
j_k \ & \ B_m \ \
\end{array} \!\!\! \right)_{k\in\br,\ l\in\bs,\ m\in\bt}
\]
to indicate that $\al$ has transversal blocks $\{i_1,j_1'\},\ldots,\{i_r,j_r'\}$, upper hooks $A_1,\ldots,A_s$, and lower hooks $B_1',\ldots,B_t'$.  (If $B\sub\bn$, we write $B'=\set{b'}{b\in B}$.)  Note that it is possible for any of $r,s,t$ to be $0$.  In fact, since we do not list the singleton blocks of $\al$, it is possible to have $r=s=t=0$.  For the same reason, it is not possible to determine $n$ from such a tabular representation of $\al\in\PB_n$, so the context will always be made clear.  We will also use variations of this notation.  So, for example, we may write
\[
\dropcustpartn{}{}
{
\fill[white] (1,0)circle(.2);
\stlinexs{3/5}{2}
\stlinexs{5/7}{2}
\uarcxx12{.4}{2}
\uarcxx47{.6}{2}
\darcxx12{.4}{2}
\darcxx34{.4}{2}
\darcxx68{.6}{2}
\uvertxs{1,2,3,4,5,6,7,8}{2}
\uvertxs{1,2,3,4,5,6,7,8}{4}
}
= \left( \begin{array}{c|c|c|c|c} 
3 & 5 & \multicolumn{3}{c}{\begin{array}{c|c} 1,2 \ \ & \ \ 4,7 \ \end{array}} \\ \cline{3-5}
5 & 7 & 1,2 & 3,4 & 6,8  \ \
\end{array} \!\!\! \right).
\]
To gain familiarity with this notation, the reader may like to check that
\[
\left( \begin{array}{c|c|c|c|c|c|c|c|c|c} 
1 & 8 & 20 & 2,7 & 3,4 & 5,6 & 9,12 & 13,18 & 14,17 & 15,16 \ \  \\ \cline{4-10}
7 & 8 & 17 & 1,4 & 2,3 & 9,16 & 10,11 & 12,15 & 13,14 & 18,20  \ \
\end{array} \!\!\! \right) \in E(\M_{20}).
\]
The next result follows quickly from Proposition \ref{prop:RSS}(i) or from \cite[Lemma 4]{EF}.

\ms
\begin{lemma}
Projections of $\PB_n$ are of the form
\[\epfreseq
\left( \begin{array}{c|c|c|c|c|c} 
i_1 & \cdots & i_r & A_1 & \cdots & A_s\ \ \\ \cline{4-6}
i_1 & \cdots & i_r & A_1 & \cdots & A_s \ \
\end{array} \!\!\! \right).
\]
\end{lemma}

Green's relations have been characterized for many of the above semigroups; see for example~\cite{Wilcox2007,FL2011,GMbook,Fernandes2001,EastGray}.  

\ms
\begin{thm}\label{thm:Green}
Let $\al,\be\in\K_n$, where $\K_n$ is one of $\PB_n,\B_n,\M_n,\I_n,\mathcal J_n,\O_n$.  Then
\bit
\itemit{i} $\alpha \R \beta \iff \dom(\alpha)=\dom(\beta)$ and $\ker(\alpha)=\ker(\beta)$,
\itemit{ii} $\alpha \L \beta\iff\codom(\alpha)=\codom(\beta)$ and $\coker(\alpha)=\coker(\beta)$,
\itemit{iii} $\alpha \J \beta\iff\alpha \D \beta\iff\rank(\alpha) = \rank(\beta)$. 
\eit
\end{thm}

\pf Since each $\K_n$ is a regular ($*$-)subsemigroup of the partition monoid $\P_n$, parts (i) and (ii) follow from \cite[Theorem~17]{Wilcox2007}.  For (iii), we note that also \cite[Theorem~17]{Wilcox2007} gives
\[
\al\D^{\K_n}\be \implies \al\D^{\P_n}\be \implies \rank(\al)=\rank(\be).
\]
Conversely, if $\rank(\al)=\rank(\be)$, then we may write
\[
\al = \left( \begin{array}{c|c} 
i_u \ & \ A_v \ \ \\ \cline{2-2}
j_u \ & \ B_w \ \
\end{array} \!\!\! \right)_{u\in U,\ v\in V,\ w\in W}
\AND
\be = \left( \begin{array}{c|c} 
k_u \ & \ C_x \ \ \\ \cline{2-2}
l_u \ & \ D_y \ \
\end{array} \!\!\! \right)_{u\in U,\ x\in X,\ y\in Y}.
\]
But then it is easy to see that $\al\R\ga\L\be$, where
\[\epfreseq
\ga = \left( \begin{array}{c|c} 
i_u \ & \ A_v \ \ \\ \cline{2-2}
l_u \ & \ D_y \ \
\end{array} \!\!\! \right)_{u\in U,\ v\in V,\ y\in Y}\in\K_n.
\]

\ms
\begin{rem}
For some of the above monoids $\K_n$, the descriptions of the $\R$- and $\L$-relations on $\K_n$ may be simplified.  For example, if $\K_n=\B_n$ or $\mathcal J_n$, then $\al\R\be \iff \ker(\al)=\ker(\be)$; or if $\K_n=\I_n$ or $\O_n$, then $\al\R\be \iff \dom(\al)=\dom(\be)$.
\end{rem}

The next result follows quickly from Theorem \ref{thm:Green}.  Again, the statement is well known for some of the semigroups \cite{EastGray,Fernandes2001,Garba1994}.  Recall that a finite semigroup is \emph{aperiodic} if the $\H$-relation is trivial.

\ms
\begin{prop}\label{prop:DClasses}
Let $\K_n$ be one of $\PB_n,\B_n,\M_n,\I_n,\mathcal J_n,\O_n$.  Then the $\D$-classes of $\K_n$ are precisely the sets
\begin{align*}
D_r(\K_n) &= \set{\al\in\K_n}{\rank(\al)=r} \qquad\text{for $0\leq r\leq n$,}
\intertext{where we must have $r\equiv n\ (\operatorname{mod}\, 2)$ if $\K_n=\B_n$ or $\mathcal J_n$.
These $\D$-classes form a chain: $D_r(\K_n)<D_s(\K_n)$ if and only if $r<s$.  The ideals of $\K_n$ are precisely the sets}
I_r(\K_n) &= \set{\al\in\K_n}{\rank(\al)\leq r} \qquad\text{for $0\leq r\leq n$.}
\end{align*}
If $\K_n$ is any of $\PB_n$, $\I_n$ or $\B_n$, then the $\H$-class of an idempotent $\al\in E(\K_n)$ is isomorphic to $\S_r$, where $r=\rank(\al)$.  If $\K_n$ is any of $\M_n$, $\O_n$ or $\mathcal J_n$, then $\K_n$ is aperiodic.  \epfres
\end{prop}



We also need to know the number of $\R$-classes contained in each $\D$-class of the semigroups we consider.  For convenience, we state the values for the Motzkin monoid $\M_n$ separately.
Recall that for a non-negative integer $k$, we write $k!!=0$ if $k$ is even or $k!!=k\cdot(k-2)\cdots3\cdot1$ if $k$ is odd, and that we interpret $(-1)!!=1$.

\ms
\begin{prop}\label{prop:DClassesofKn}
Let $\K_n$ be one of $\PB_n,\B_n,\I_n,\mathcal J_n,\O_n$.  The number of $\R$-classes (equivalently, $\L$-classes) contained in $D_r(\K_n)$ is equal to 
\bit
\itemit{i} ${n\choose r}$ if $\K_n=\I_n$ or $\O_n$,
\itemit{ii} ${n\choose r}\cdot(n-r-1)!!$ if $\K_n=\B_n$,
\itemit{iii} $\frac{r+1}{n+1}{n+1 \choose k}$ if $\K_n=\mathcal J_n$ and $r=n-2k$,
\itemit{iv} ${n \choose r} \cdot a(n-r)$ if $\K_n=\PB_n$,
where the numbers $a(m)$ satisfy 
\[
a(0)=a(1)=1 \COMMA a(m)=a(m-1)+(m-1)a(m-2) \quad\text{for $m\geq2$.}
\]
\eit
\end{prop}

\pf Part (i) is well known; see for example \cite[Proposition 4.6.3]{GMbook} and \cite[Proposition 2.3]{Fernandes2001}.  Parts (ii) and (iii) were proved in \cite[Theorems 8.4 and 9.5]{EastGray}; part (ii) was also proved in \cite[Corollary 3.6(i)]{Larsson}.  A projection $\al$ from $D_r(\PB_n)$ is determined by $\dom(\al)$ and $\ker(\al)$.  There are ${n\choose r}$ ways to choose $\dom(\al)$, and then $\ker(\al)$ is determined by a set partition of $\dom(\al)^c$ into blocks of size at most $2$ (since all elements of $\dom(\al)$ belong to singleton $\ker(\al)$-classes), so there are $a(n-r)$ choices for $\ker(\al)$. \epf



\ms
\begin{prop}\label{prop:DClassesofMn}
Let $m(n,r)$ be the number of $\R$-classes (equivalently, $\L$-classes) contained in $D_r(\M_n)$.  Then
\begin{myalign}
\nonumber m(0,0)=1 \COMMA
m(n,r)=0 &\qquad\qquad\text{if $n<r$ or $r<0$,}\\
\nonumber m(n,r)=m(n-1,r-1)+m(n-1,r)+m(n-1,r+1) &\qquad\qquad\text{if $0\leq r\leq n$ and $n\geq1$.}
\end{myalign}
\end{prop}

\pf It is clear that $m(0,0)=1$, and that $m(n,r)=0$ if $n<r$ or $r<0$.  So suppose $0\leq r\leq n$, where $n\geq1$.  We define a mapping
\[
\bar{}:P(D_r(\M_n)) \to P(D_{r-1}(\M_{n-1})) \cup P(D_r(\M_{n-1})) \cup P(D_{r+1}(\M_{n-1})):\al\mt\bar\al
\]
as follows.  Let $\al\in P(D_r(\M_n))$.  We write $\al^\flat$ for the induced subgraph on vertices $\bnf\cup\bnf'$.  If $\{n,n'\}$ or $\{n\}$ is a block of $\al$, then we write $\bar\al=\al^\flat$, noting that $\bar\al$ then belongs to $P(D_{r-1}(\M_{n-1}))$ or $P(D_r(\M_{n-1}))$, respectively.  If $\{i,n\}$ is a block of $\al$ for some $i\in\bnf$, then we write $\bar\al$ for the Motzkin diagram obtained by adding the edge $\{i,i'\}$ to $\al^\flat$, noting that $\bar\al$ belongs to $P(D_{r+1}(\M_{n-1}))$ in this case, and that $i=\max(\dom(\bar\al))$.  Since the $\bar{\ }$ map is clearly invertible, the result follows. \epf

We will also need to count certain special kinds of projections from $\M_n$.  With this in mind, let $\al\in\M_n$, and suppose $A$ is an upper non-transversal block of $\al$ (of size $1$ or $2$).  We say that $A$ is \emph{nested} if there exists an upper hook $\{x,y\}$ of $\al$ such that $x<\min(A)\leq\max(A)<y$.  (We similarly define nested lower blocks.)  For example, the blocks $\{4,5\}$ and $\{6\}$ are nested, while $\{2\}$ and $\{3,7\}$ are unnested, in
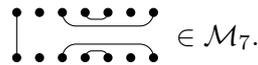
\[
\custpartn{1,2,3,4,5,6,7}{1,2,3,4,5,6,7}{\stline11\darcx37{.6}\darcx45{.3}\uarcx37{.6}\uarcx45{.3}}\in\M_7.
\]

\ms
\begin{prop}\label{prop:m'nr}
Let $m'(n,r)$ denote the number of projections from $D_r(\M_n)$ with no unnested singleton blocks.  Then
\begin{myalign}
\nonumber m'(0,0)=1 \COMMA
m'(n,r)=0 &\qquad\qquad\text{if $n<r$ or $r<0$,}\\
\nonumber m'(n,r)=m'(n-1,r-1)+\sum_{j=r+1}^{n-1}m'(n-1,j) &\qquad\qquad\text{if $0\leq r\leq n$ and $n\geq1$.}
\end{myalign}
\end{prop}

\pf Let $P'(D_r(\M_n))$ denote the set of projections from $D_r(\M_n)$ with no unnested singleton blocks, so that $m'(n,r)=|P'(D_r(\M_n))|$.
Again, the initial values are clear, so suppose $0\leq r\leq n$, where $n\geq1$.  This time, we define a mapping
\[
\bar{}:P'(D_r(\M_n)) \to P'(D_{r-1}(\M_{n-1})) \cup \bigcup_{j=r+1}^{n-1} P'(D_j(\M_{n-1})) :\al\mt\bar\al
\]
as follows.  Let $\al\in P'(D_r(\M_n))$.  If $\{n,n'\}$ is a block of $\al$, then we put $\bar\al=\al^\flat$ (as defined in the previous proof).  If $\{i,n\}$ is a block of $\al$ for some $i\in\bnf$, then we let $\bar\al$ be the Motzkin diagram obtained from~$\al^\flat$ by adding edges $\{k,k'\}$ for each unnested singleton block $\{k\}$ of $\al^\flat$ with $i\leq k\leq n-1$.  (In particular, note that $\{i,i'\}$ is added.)  Again, the $\bar{\ }$ map is invertible, so the result follows. \epf

\ms
\begin{rem}\label{rem:mnr1}
The first few values of $m(n,r)$ and $m'(n,r)$ are given in Table \ref{tab:mnr}.  These are sequences A064189 and A097609 on the OEIS \cite{OEIS}; see also \cite{BPS1992,Bernhart1999,AMartin1998}.  The numbers $m(n)=m(n,0)$ and $m'(0)=m'(n,0)$ are sequences A001006 and A005043 on the OEIS, and are called the \emph{Motzkin} and \emph{Riordan} numbers, respectively.  The numbers $m(n,r)$ and $m'(n,r)$ 
also count various kinds of lattice paths;  bijections between such paths and Motzkin projections lead to different proofs of Propositions~\ref{prop:DClassesofMn}~and~\ref{prop:m'nr}.
\end{rem}




\begin{table}[h]%
\begin{center}
\begin{tabular}{|c|rrrrrrrr|}
\hline
$n$ $\setminus$ $r$	&0	&1	&2	&3	&4	&5	&6	&7    \\
\hline
0&	1&	&&&&&&\\
1&	1&		1&	&&&&&\\
2&	2&		2&		1&	&&&&\\
3&	4&		5&		3&		1&	&&&\\
4&	9&		12&		9&		4&		1&	&&\\
5&	21&		30&		25&		14&		5&		1&	&\\
6&	51&		76&		69&		44&		20&		6&		1&	\\
7&	127&		196&		189&		133&		70&		27&		7&		1\\
\hline
\end{tabular}
\qquad\qquad
\begin{tabular}{|c|rrrrrrrr|}
\hline
$n$ $\setminus$ $r$	&0	&1	&2	&3	&4	&5	&6	&7     \\
\hline
0&	1&	&&&&&&\\
1&	0&		1&	&&&&&\\
2&	1&		0&		1&	&&&&\\
3&	1&		2&		0&		1&	&&&\\
4&	3&		2&		3&		0&		1&	&&\\
5&	6&		7&		3&		4&		0&		1&	&\\
6&	15&		14&		12&		4&		5&		0&		1&	\\
7&	36&		37&		24&		18&		5&		6&		0&		1\\
\hline
\end{tabular}
\end{center}
\vspace{-5mm}
\caption{Values of $m(n,r)$ and $m'(n,r)$ --- left and right, respectively.}
\label{tab:mnr}
\end{table}

\ms
\begin{rem}\label{rem:|IrMn|}
From Propositions \ref{prop:DClasses}, \ref{prop:DClassesofKn} and \ref{prop:DClassesofMn}, we may easily deduce formulae for the sizes of the $\D$-classes and ideals of the monoids we consider.  For example,
\[
|I_r(\PB_n)| = \sum_{s=0}^r{n\choose s}^2a(n-s)^2s! \AND |I_r(\M_n)| = \sum_{s=0}^r m(n,s)^2.
\]
Calculated values of $|I_r(\PB_n)|$ and $|I_r(\M_n)|$ are given in Table \ref{tab:|Ir|}.
\end{rem}

\begin{table}[h]%
\begin{center}
\begin{tabular}{|c|rrrrrrrr|}
\hline
$n$ $\setminus$ $r$	&0	&1	&2	&3	&4	&5	&6	&7 \\
\hline
0 &1 &&&&&&&\\
1 &1& 2&&&&&&\\
2 &4& 8& 10 &&&&&\\
3 &16& 52& 70& 76 &&&&\\
4 &100& 356& 644& 740& 764 &&&\\
5 &676& 3176& 6376& 8776& 9376& 9496 &&\\
6 &5776& 30 112& 75 112& 113 512& 135 112& 139 432& 140 152&\\
7 &53 824& 336 848& 933 080& 1 668 080& 2 138 480& 2 350 160& 2 385 440& 2 390 480\\
\hline
\end{tabular}
\\~\\~\\
\begin{tabular}{|c|rrrrrrrr|}
\hline
$n$ $\setminus$ $r$	&0	&1	&2	&3	&4	&5	&6	 &7    \\
\hline
0 &1 &&&&&&&\\
1 &1& 2&&&&&&\\
2 &4& 8& 9 &&&&&\\
3 &16& 41& 50& 51 &&&&\\
4 &81& 225& 306& 322& 323 &&&\\
5 &441& 1341& 1966& 2162& 2187& 2188 &&\\
6 &2601& \phantom{3}8 377& 13 138& \phantom{3}15 074& \phantom{3}15 474& \phantom{3}15 510& \phantom{3}15 511&\\
7&16 129& \phantom{3}54 545& \phantom{3}90 266& \phantom{3 }107 955& \phantom{3 }112 855& \phantom{3 }113 584& \phantom{3 }113 633& \phantom{3 }113 634\\
\hline
\end{tabular}
\end{center}
\vspace{-5mm}
\caption{Values of $|I_r(\PB_n)|$ and $|I_r(\M_n)|$ --- top and bottom, respectively.}
\label{tab:|Ir|}
\end{table}

\begin{figure}[h]
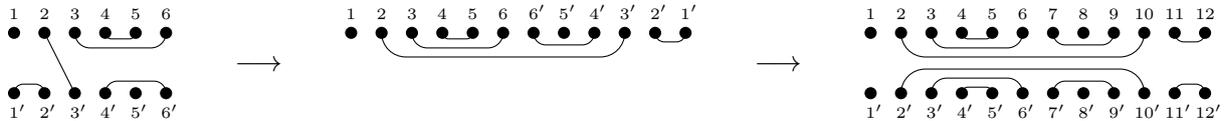

   \begin{center}
\[
\dropcustpartnxy{1,2,3,4,5,6}{1,2,3,4,5,6}
{\stlinexs{2/3}{0}
\uarcxx36{.5}{0}
\uarcxx45{.2}{0}
\darcxx12{.3}{0}
\darcxx46{.4}{0}
\vertlabelshhh{1,...,6}
\vertlabeldns{1,...,6}}
\quad\longrightarrow\quad
\dropcustpartnxy{1,...,12}{}
{
\uarcxx36{.5}{0}
\uarcxx45{.2}{0}
\uarcxx79{.4}{0}
\uarcxx{11}{12}{.3}{0}
\uarcxx2{10}{.8}{0}
\vertlabelshhh{1,...,6}
\vertlabelupdashess{7/6,8/5,9/4,10/3,11/2,12/1}
\vertlabeldnsph{1,...,6}}
\quad\longrightarrow\quad
\dropcustpartnxy{1,...,12}{1,...,12}
{
\uarcxx36{.5}{0}
\uarcxx45{.2}{0}
\uarcxx79{.4}{0}
\uarcxx{11}{12}{.3}{0}
\uarcxx2{10}{.8}{0}
\darcxx36{.5}{0}
\darcxx45{.2}{0}
\darcxx79{.4}{0}
\darcxx{11}{12}{.3}{0}
\darcxx2{10}{.8}{0}
\vertlabelshhh{1,...,12}
\vertlabeldns{1,...,12}
}
\]
    \vspace{-5mm}
\caption{An example illustrating the bijection $\M_n\to P(D_0(\M_{2n}))$; see Remark \ref{rem:mnr2}.}
    \label{fig:|M_n|=m(2n)}
   \end{center}
 \end{figure}

\ms
\begin{rem}\label{rem:mnr2}
There is a bijection between $\M_n$ and $P(D_0(\M_{2n}))$ that may be described as follows.  For $\al\in\M_n$, first ``unfold'' the graph so that all vertices appear on a straight line in the order $1,\ldots,n,n',\ldots,1'$, with all edges ``hanging'' below the vertices.  Then relabel the vertices $1,\ldots,2n$ and add a reflected copy directly below, with vertices labelled $1',\ldots,(2n)'$.  An example is given in Figure~\ref{fig:|M_n|=m(2n)}.  In particular, it follows that $|\M_n|=|P(D_0(\M_{2n}))|=m(2n,0)=m(2n)$, as observed in \cite{BH2014}.  Combined with the formula for $|\M_n|=|I_n(\M_n)|$ from Remark \ref{rem:|IrMn|}, we obtain the well-known identity $m(2n)=\sum_{r=0}^nm(n,r)^2$, as also observed in \cite{BH2014}.  Similar considerations show that $|\PB_n|=|P(D_0(\PB_{2n}))|=a(2n)=\sum_{r=0}^n{n\choose r}^2a(n-r)^2r!$.
\end{rem}

We conclude this section by stating a result that will be used frequently in what follows; for a proof, see \cite[Lemmas~4.1 and~4.7]{ZF2015}.  Similar results were proved for the diagram semigroups $\P_n,\B_n,\mathcal J_n$ in \cite{EastGray}.

\ms
\begin{prop}\label{prop:ZF}
If $0\leq r\leq n-1$, then $I_r(\I_n)=\la D_r(\I_n)\ra$ and $I_r(\O_n)=\la D_r(\O_n)\ra$. \epfres
\end{prop}

\section{The partial Brauer monoid $\PB_n$}\label{sect:PBn}

In this section, we focus on the partial Brauer monoid $\PB_n$.  We calculate the rank of each ideal $I_r(\PB_n)$.  Additionally, we show that $I_r(\PB_n)$ is idempotent-generated if $r\leq n-2$ and we show that $\idrank(I_r(\PB_n))=\rank(I_r(\PB_n))$ for such an $r$; see Theorem \ref{thm:IrPBn}.  We also characterise the elements of the idempotent-generated subsemigroup $\bbE(\PB_n)=\la E(\PB_n)\ra$, and calculate its rank and idempotent rank (which, again, are equal); see Theorem \ref{thm:EPBn}.



Throughout this section, we will use the abbreviations
\begin{align*}
D_r=D_r(\PB_n) = \set{\al\in\PB_n}{\rank(\al)=r} \AND
 I_r=I_r(\PB_n) = \set{\al\in\PB_n}{\rank(\al)\leq r}
\end{align*}
for each $0\leq r\leq n$.  Note that $I_r=D_0\cup\cdots \cup D_r$, and that $D_n=\S_n$, $I_n=\PB_n$ and $I_{n-1}=\PB_n\sm\S_n$.

For each $\oijn$, and each $1\leq k\leq n$, let 
\[
\si_{ij} =  \custpartn{1,3,4,5,7,8,9,11}{1,3,4,5,7,8,9,11}{\dotsups{1/3,5/7,9/11}\dotsdns{1/3,5/7,9/11}\stlines{1/1,3/3,4/8,5/5,7/7,8/4,9/9,11/11}\vertlabelshhh{1/1,4/i,8/j,11/n}} \COMMA
\tau_{ij} =  \custpartn{1,3,4,5,7,8,9,11}{1,3,4,5,7,8,9,11}{\dotsups{1/3,5/7,9/11}\dotsdns{1/3,5/7,9/11}\stlines{1/1,3/3,5/5,7/7,9/9,11/11}\uarc48\darc48\vertlabelshhh{1/1,4/i,8/j,11/n}} \COMMA
\ve_k = \custpartn{1,3,4,5,7}{1,3,4,5,7}{\dotsups{1/3,5/7}\dotsdns{1/3,5/7}\stlines{1/1,3/3,5/5,7/7}\vertlabelshhh{1/1,4/k,7/n}}.
\]
So $\si_{ij}\in\S_n$, $\tau_{ij}\in E(\B_n)$ and $\ve_k\in E(\I_n)$.

For $A=\{a_1<\cdots<a_r\}\sub\bn$, we define
\[
\lam_A = \left( \begin{array}{c|c|c} 
a_1 \ & \ \cdots \ & \ a_r \ \ \\ 
1 \ & \ \cdots \ & \ r \ \
\end{array} \!\!\! \right)
\AND
\rho_A = \left( \begin{array}{c|c|c} 
1 \ & \ \cdots \ & \ r \ \ \\ 
a_1 \ & \ \cdots \ & \ a_r \ \
\end{array} \!\!\! \right).
\]
So $\lam_A$ is the unique element of $\POI_n$ with domain $A$ and codomain $\br$, and $\rho_A=\lam_A^*$.  For example, with $n=7$ and $A=\{1,3,4,6\}$, we have
\[
\lam_A = \custpartn{1,2,3,4,5,6,7}{1,2,3,4,5,6,7}{\stlines{1/1,3/2,4/3,6/4}}
\AND
\rho_A = \custpartn{1,2,3,4,5,6,7}{1,2,3,4,5,6,7}{\stlines{1/1,2/3,3/4,4/6}}.
\]
These mappings were introduced in \cite{JEinsn}, where they were used to derive a presentation for $\O_n$ and the singular ideal $I_{n-1}(\I_n)=\InSn$ of the symmetric inverse monoid $\I_n$.  Specifically, it was shown that $\O_n$ is generated (as a semigroup) by the set $\{1\}\cup\set{\lam_{\{i\}^c}}{i\in\bn}\cup\set{\rho_{\{i\}^c}}{i\in\bn}$.  Note that $\al=\lam_{\dom(\al)}\cdot\rho_{\codom(\al)}$ for any $\al\in\O_n$.

For $A\sub\bn$, we write
\[
\S_A = \set{\al\in\I_n}{\dom(\al)=\codom(\al)=A}.
\]
So clearly $\S_A$ is isomorphic to the symmetric group $\S_r$, where $r=|A|$.  In fact, $\S_A$ is the $\H$-class (in $\PB_n$ or $\I_n$) of the idempotent 
\[
\id_A 
= \left( \begin{array}{c|c} 
a\ & \ \phantom{\emptyset} \ \ \\ \cline{2-2}
a \ & \ \phantom{\emptyset}  \ \
\end{array} \!\!\! \right)_{a\in A}.
\]
The next result gives a useful normal form for the elements of $\PB_n$.



\ms
\begin{lemma}\label{lem:normalformPBn}
If $\al\in D_r$, then
$
\al = \be\cdot \lam_{\dom(\al)} \cdot \ga \cdot \rho_{\codom(\al)} \cdot \de
$
for some $\be,\de\in E(D_r)$ and $\ga\in\S_{\br}$.
\end{lemma}

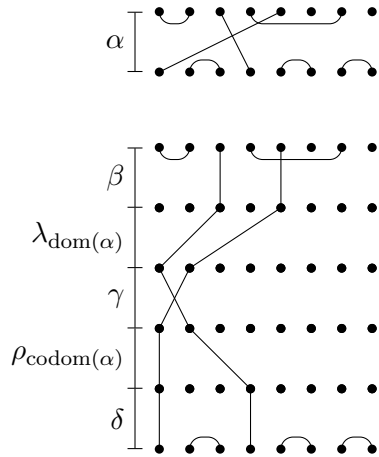
\begin{figure}[h]
\begin{center}
\begin{tikzpicture}[xscale=.4,yscale=0.4]
\uvertths{1,...,8}{.15}
\lvertths{1,...,8}{.15}
\uarc12
\uarc47
\darc23
\darc56
\darc78
\stline34
\stline51
\draw[|-|] (.2,0)--(.2,2);
\draw(.2,1)node[left]{$\al$};
\begin{scope}[shift={(0,-4.5)}]	
\uvertths{1,...,8}{.15}
\lvertths{1,...,8}{.15}
\uarc12
\uarc47
\stline33
\stline55
\draw[|-|] (.2,0)--(.2,2);
\draw(.2,1)node[left]{$\be$};
\end{scope}
\begin{scope}[shift={(0,-6.5)}]	
\uvertths{1,...,8}{.15}
\lvertths{1,...,8}{.15}
\stline31
\stline52
\draw[|-] (.2,0)--(.2,2);
\draw(.2,1)node[left]{$\lam_{\dom(\al)}$};
\end{scope}
\begin{scope}[shift={(0,-8.5)}]	
\uvertths{1,...,8}{.15}
\lvertths{1,...,8}{.15}
\stline12
\stline21
\draw[|-] (.2,0)--(.2,2);
\draw(.2,1)node[left]{$\ga$};
\end{scope}
\begin{scope}[shift={(0,-10.5)}]	
\uvertths{1,...,8}{.15}
\lvertths{1,...,8}{.15}
\stline11
\stline24
\draw[|-] (.2,0)--(.2,2);
\draw(.2,1)node[left]{$\rho_{\codom(\al)}$};
\end{scope}
\begin{scope}[shift={(0,-12.5)}]	
\uvertths{1,...,8}{.15}
\lvertths{1,...,8}{.15}
\darc23
\darc56
\darc78
\stline11
\stline44
\draw[|-] (.2,0)--(.2,2);
\draw(.2,1)node[left]{$\de$};
\end{scope}
\end{tikzpicture}
\end{center}
\vspace{-5mm}
\caption{An example illustrating the factorization $\al = \be\cdot \lam_{\dom(\al)} \cdot \ga \cdot \rho_{\codom(\al)} \cdot \de$ from Lemma \ref{lem:normalformPBn}.}
\label{fig:normalformPBn}
\end{figure}

\pf Suppose $\dom(\al)=\{i_1<\cdots<i_r\}$ and $\codom(\al)=\{j_1<\cdots<j_r\}$.  Then we may write
\[
\al = \left( \begin{array}{c|c|c|c|c|c} 
i_1 \ & \ \cdots \ & \ i_r \ & \ A_1 \ & \ \cdots \ & \ A_s\ \ \\ \cline{4-6}
j_{1\pi} \ & \ \cdots \ & \ j_{r\pi} \ & \ B_1 \ & \ \cdots \ & \ B_t \ \
\end{array} \!\!\! \right)
\]
for some $\pi\in\S_r$.  We then extend $\pi\in\S_r$ to $\ga\in\S_{\br}$, and we see that the result holds with
\[
\be = \left( \begin{array}{c|c|c|c|c|c} 
i_1 \ & \ \cdots \ & \ i_r \ & \ A_1 \ & \ \cdots \ & \ A_s\ \ \\ \cline{4-6}
i_1 \ & \ \cdots \ & \ i_r \ & \multicolumn{3}{|c}{} \ \
\end{array} \!\!\! \right)
\AND
\de = \left( \begin{array}{c|c|c|c|c|c} 
j_1 \ & \ \cdots \ & \ j_r \ & \multicolumn{3}{|c}{} \\ \cline{4-6}
j_1 \ & \ \cdots \ & \ j_r \ & \ B_1 \ & \ \cdots \ & \ B_t \ \
\end{array} \!\!\! \right),
\]
both of which are clearly idempotents or rank $r$. (An example calculation is given in Figure \ref{fig:normalformPBn}.)  \epf

\ms
\begin{rem}\label{rem:normalformPBn}
If $\al\in\M_n$, then $\pi\in\S_r$ from the proof of Lemma \ref{lem:normalformPBn} is the identity permutation, so in fact we have $\al = \be\cdot \lam_{\dom(\al)} \cdot \rho_{\codom(\al)} \cdot \de$, and it is clear that the idempotents $\be,\de$ from the proof are planar, so that $\be,\de\in E(D_r(\M_n))$.
\end{rem}

\ms
\begin{rem}
Lemma \ref{lem:normalformPBn} and Remark \ref{rem:normalformPBn} are of broad general interest, beyond our current investigations (which mostly concern ideals, idempotent generation, and minimal size (idempotent) generating sets).  For example, in Section \ref{sect:algebras}, we show how they may be used to describe cellular structures on the corresponding diagram algebras.  We also note that similar normal form results have been extremely useful in derivations of presentations (generators and relations) for a number of diagram semigroups and transformation semigroups; see for example \cite{JEinsn,JEgrpm,JEpnsn} and references therein.  We therefore expect that Lemma \ref{lem:normalformPBn} and Remark \ref{rem:normalformPBn} will provide useful approaches to presentations for $\PB_n$ and $\M_n$; see \cite{PHY2013,KM2006} for existing presentations.
\end{rem}

Armed with Lemma \ref{lem:normalformPBn}, our first task is to show that the ideals $I_r$ with $r\leq n-2$ are idempotent-generated (see Corollary \ref{cor:IdDrDr-1}).  To do this, we will need several intermediate results.

\ms
\begin{lemma}\label{lem:HrEDr}
If $0\leq r\leq n-2$, then $\S_{\br}\sub\la E(D_r)\ra$.
\end{lemma}

\pf  If $r\leq1$, then $\S_{\br}=\{\id_\br\}\sub E(D_r)$, so suppose $r\geq2$.  Since $\S_{\br}\cong\S_r$, and since $\S_r$ is generated by transpositions,  it follows that $\S_{\br}$ is generated by the set $\set{\si_{ij}\cdot\id_\br}{\oijr}$.  So it suffices to show that each element of this generating set belongs to $\la E(D_r)\ra$.  With this in mind, let $\oijr$.  Then
\[
\si_{ij}\cdot\id_\br \ \ = \gcustpartn{4,8,12,13,1,3,5,7,9,11,14,16}{}{4,8,12,13,1,3,5,7,9,11,14,16}{}
{\dotsups{1/3,5/7,9/11,14/16}\dotsdns{1/3,5/7,9/11,14/16}\gstlines{}\stlines{1/1,3/3,5/5,7/7,9/9,11/11,4/8,8/4}\vertlabelshh{1/1,4/i,8/j,11/r,16/n}} = 
\gcustpartndash{4,8,12,13,1,3,5,7,9,11,14,16}{}{4,8,12,13,1,3,5,7,9,11,14,16}{}
{
\dotsups{1/3,5/7,9/11,14/16}\dotsdns{1/3,5/7,9/11,14/16}
\uvertxs{1,3,5,7,9,11,14,16}{4}\uvertxs{4,8,12,13}{4}\dotsupxs{1/3,5/7,9/11,14/16}{4}
\uvertxs{1,3,5,7,9,11,14,16}{6}\uvertxs{4,8,12,13}{6}\dotsupxs{1/3,5/7,9/11,14/16}{6}
\uvertxs{1,3,5,7,9,11,14,16}{8}\uvertxs{4,8,12,13}{8}\dotsupxs{1/3,5/7,9/11,14/16}{8}
\stlinexs{1/1,3/3,5/5,7/7,9/9,11/11}{0}
\stlinexs{1/1,3/3,5/5,7/7,9/9,11/11}{2}
\stlinexs{1/1,3/3,5/5,7/7,9/9,11/11}{4}
\stlinexs{1/1,3/3,5/5,7/7,9/9,11/11}{6}
\stlinexs{4/4,8/8}{0}
\stlinexs{4/4,8/8}{6}
\darcxx{12}{13}{.4}{6}
\uarcxx{12}{13}{.4}{0}
\uarcxx{8}{12}{.4}{4}
\darcxx{8}{12}{.4}{4}
\stlinexs{4/4,13/13}{4}
\uarcxx{4}{12}{.6}{2}
\darcxx{4}{12}{.6}{2}
\stlinexs{8/8,13/13}{2}
\vertlabelxs{1/1,4/i,8/j,11/r,16/n}{6}
},
\]
with each of the four diagrams in the product on the right hand side belonging to $E(D_r)$.   \epf

\ms
\begin{lemma}\label{lem:lamAEDn-2}
If $A\sub\bn$ with $|A|= n-2$, then $\lam_A,\rho_A\in\la E(D_{n-2})\ra$.
\end{lemma}

\pf By symmetry, we need only prove the statement concerning $\lam_A$.  So suppose $A=\{k<l\}^c$.  We prove that $\lam_A\in\la E(D_{n-2})\ra$ by descending induction on $k+l$.  The maximum value of $k+l$ occurs when $k=n-1$ and $l=n$, in which case $\lam_A=\id_A\in E(D_{n-2})$, so suppose $k+l\leq 2n-2$.  We consider two cases.

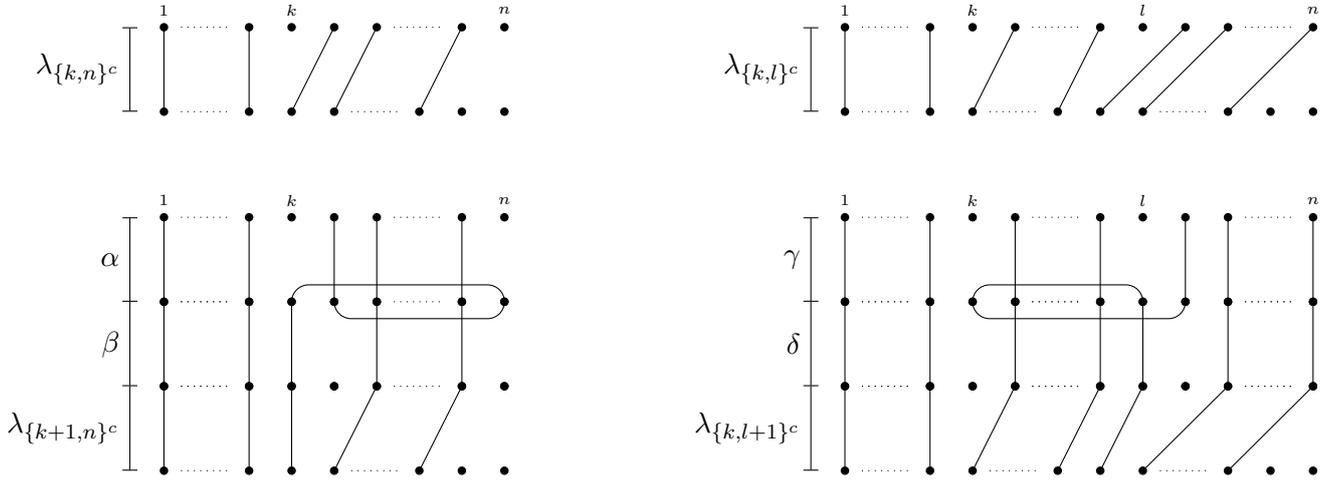
\begin{figure}[h]
\begin{center}
\begin{tikzpicture}[xscale=.56,yscale=0.56]
\uvs{1,3,4,5,6,8,9}
\lvs{1,3,4,5,7,8,9}
\stline11
\stline33
\stline54
\stline65
\stline87
\vertlabelsh{1/1,4/k,9/n}
\ldottedsms{1/3,5/7}
\udottedsms{1/3,6/8}
\draw[|-|] (.2,0)--(.2,2);
\draw(.2,1)node[left]{$\lam_{\{k,n\}^c}$};
\begin{scope}[shift={(16,0)}]	
\uvs{1,3,4,5,7,8,9,10,12}
\lvs{1,3,4,6,7,8,10,11,12}
\stlines{1/1,3/3,5/4,7/6,9/7,10/8,12/10}
\vertlabelsh{1/1,4/k,8/l,12/n}
\udottedsms{1/3,5/7,10/12}
\ldottedsms{1/3,4/6,8/10}
\draw[|-|] (.2,0)--(.2,2);
\draw(.2,1)node[left]{$\lam_{\{k,l\}^c}$};
\end{scope}
\begin{scope}[shift={(0,-4.5)}]	
\uvs{1,3,4,5,6,8,9}
\lvs{1,3,4,5,6,8,9}
\stline11
\stline33
\stline55
\stline66
\stline88
\darc49
\vertlabelsh{1/1,4/k,9/n}
\ldottedsms{1/3}
\udottedsms{1/3,6/8}
\draw[|-|] (.2,0)--(.2,2);
\draw(.2,1)node[left]{$\al$};
\end{scope}
\begin{scope}[shift={(16,-4.5)}]	
\uvs{1,3,4,5,7,8,9,10,12}
\lvs{1,3,4,5,7,8,9,10,12}
\stlines{1/1,3/3,5/5,7/7,9/9,10/10,12/12}
\vertlabelsh{1/1,4/k,8/l,12/n}
\darc48
\udottedsms{1/3,5/7,10/12}
\ldottedsms{1/3,5/7,10/12}
\draw[|-|] (.2,0)--(.2,2);
\draw(.2,1)node[left]{$\ga$};
\end{scope}
\begin{scope}[shift={(0,-6.5)}]	
\uvs{1,3,4,5,6,8,9}
\lvs{1,3,4,5,6,8,9}
\stline11
\stline33
\stline44
\stline66
\stline88
\uarc59
\ldottedsms{1/3}
\udottedsms{1/3,6/8}
\draw[|-] (.2,0)--(.2,2);
\draw(.2,1)node[left]{$\be$};
\end{scope}
\begin{scope}[shift={(16,-6.5)}]	
\uvs{1,3,4,5,7,8,9,10,12}
\lvs{1,3,4,5,7,8,9,10,12}
\stlines{1/1,3/3,5/5,7/7,8/8,10/10,12/12}
\uarc49
\udottedsms{1/3,5/7,10/12}
\ldottedsms{1/3,5/7,10/12}
\draw[|-] (.2,0)--(.2,2);
\draw(.2,1)node[left]{$\de$};
\end{scope}
\begin{scope}[shift={(0,-8.5)}]	
\uvs{1,3,4,5,6,8,9}
\lvs{1,3,4,5,7,8,9}
\stline11
\stline33
\stline44
\stline65
\stline87
\ldottedsms{1/3,5/7}
\udottedsms{1/3,6/8}
\draw[|-] (.2,0)--(.2,2);
\draw(.2,1)node[left]{$\lam_{\{k+1,n\}^c}$};
\end{scope}
\begin{scope}[shift={(16,-8.5)}]	
\uvs{1,3,4,5,7,8,9,10,12}
\lvs{1,3,4,6,7,8,10,11,12}
\stlines{1/1,3/3,5/4,7/6,8/7,10/8,12/10}
\udottedsms{1/3,5/7,10/12}
\ldottedsms{1/3,4/6,8/10}
\draw[|-] (.2,0)--(.2,2);
\draw(.2,1)node[left]{$\lam_{\{k,l+1\}^c}$};
\end{scope}
\end{tikzpicture}
\end{center}
\vspace{-5mm}
\caption{Diagrammatic verification of the equations $\lam_{\{k,n\}^c} = \al\be\cdot\lam_{\{k+1,n\}^c}$ and $\lam_{\{k,l\}^c} = \ga\de\cdot\lam_{\{k,l+1\}^c}$ from the proof of Lemma~\ref{lem:lamAEDn-2} --- left and right, respectively.  See the text for a full explanation.}
\label{fig:lamAEDn-2}
\end{figure}

{\bf Case 1.}  If $l=n$, then $k\leq n-2$, and $\lam_A=\al\be\cdot\lam_{\{k+1,n\}^c}$ (see Figure \ref{fig:lamAEDn-2}), where
\[
\al = \left( \begin{array}{c|c} 
x \ &  \\ \cline{2-2}
x \ & \ k,n\ \ 
\end{array} \!\!\! \right)_{x\in\{k,n\}^c}
\AND
\be = \left( \begin{array}{c|c} 
x \ & \ k+1,n\ \ \\ \cline{2-2}
x \ & 
\end{array} \!\!\! \right)_{x\in\{k+1,n\}^c}.
\]
Clearly $\al,\be\in E(D_{n-2})$, while an inductive hypothesis gives $\lam_{\{k+1,n\}^c}\in\la E(D_{n-2})\ra$.

{\bf Case 2.}  If $l<n$, then $\lam_A=\ga\de\cdot\lam_{\{k,l+1\}^c}$ (see Figure \ref{fig:lamAEDn-2}), where
\[
\ga = \left( \begin{array}{c|c} 
x \ &  \\ \cline{2-2}
x \ & \ k,l\ \ 
\end{array} \!\!\! \right)_{x\in\{k,l\}^c}
\AND
\de = \left( \begin{array}{c|c} 
x \ & \ k,l+1\ \ \\ \cline{2-2}
x \ &  \ \ 
\end{array} \!\!\! \right)_{x\in\{k,l+1\}^c}.
\]
Again, $\ga,\de\in E(D_{n-2})$ and $\lam_{\{k,l+1\}^c}\in\la E(D_{n-2})\ra$. \epf

For the proof of the next result, and in other places, it will be convenient to define some mappings.  With this in mind, we define
\[
{}^\sharp: \PB_{n-1}\to\PB_n:\al\mt\al^\sharp \COMMA \qquad
{}^\natural: \PB_{n-1}\to\PB_n:\al\mt\al^\natural \COMMA \qquad
{}^\flat: \PB_n\to\PB_{n-1}:\al\mt\al^\flat 
\]
as follows.  For $\al\in\PB_{n-1}$, $\al^\sharp\in\PB_n$ denotes the Brauer $n$-diagram obtained by adding the block $\{n,n'\}$ to $\al$.  We similarly define $\al^\natural\in\PB_n$ by the addition of the blocks $\{n\}$ and $\{n'\}$.  For $\al\in\PB_n$, $\al^\flat\in\PB_{n-1}$ denotes the induced subgraph on $\bnf\cup\bnf'$.  Note that, for example, if $\{i,n\}$ is a block of $\al\in\PB_n$, then $\{i\}$ will be a block of $\al^\flat$.  Also, note that $(\al^\sharp)^\flat=(\al^\natural)^\flat=\al$ for any $\al\in\PB_{n-1}$, but that $(\al^\flat)^\sharp=\al$ or $(\al^\flat)^\natural=\al$ will only hold when certain (easily identified) conditions are true of $\al\in\PB_n$.  Also note that ${}^\sharp$ is a monoid monomorphism, ${}^\natural$ is a semigroup monomorphism, but ${}^\flat$ is not a homomorphism (unless $n=1$).


\ms
\begin{lemma}\label{lem:lamAEDr}
If $A\sub\bn$ with $|A|=r\leq n-2$, then $\lam_A,\rho_A\in\la E(D_r)\ra$.
\end{lemma}

\pf Again, it suffices to show that $\lam_A\in\la E(D_r)\ra$, and we do this by induction on $r+n$.  First note that if $r=0$, then $\lam_A=\lam_\emptyset=\id_\emptyset\in E(D_r)$; this includes the base case of the induction, in which $r=0$ and $n=2$.  We now suppose that $1\leq r\leq n-3$, noting that the case in which $r=n-2$ is covered by Lemma~\ref{lem:lamAEDn-2}.  We consider two cases, according to whether or not $n$ belongs to $A$.

{\bf Case 1.}  Suppose first that $n\not\in A$.  
By an induction hypothesis,
we have $\lam_A^\flat=\al_1\cdots\al_k$ for some $\al_1,\ldots,\al_k\in E(D_r(\PB_{n-1}))$.  
Then $\lam_A=\al_1^\natural\cdots\al_k^\natural\in\la E(D_r(\PB_n))\ra$, completing the inductive step in this case.

{\bf Case 2.}  Now suppose $n\in A$.
By an induction hypothesis, $\lam_A^\flat \in \la E(D_{r-1}(\PB_{n-1}))\ra$, and it quickly follows that $(\lam_A^\flat)^\sharp \in \la E(D_r(\PB_n))\ra$.  But then 
$\lam_A = (\lam_A^\flat)^\sharp \cdot \be\ga$ (see Figure \ref{fig:lamAEDr})\footnote{In some figures, shading is used to convey the general ``shape'' of a Brauer diagram.  These should be used as a guide, but the explicit definitions given in the proof should always be consulted.  For example, in Figure \ref{fig:lamAEDr}, the shaded part of $\lam_A$ is understood to be identical to the shaded part of $(\lam_A^\flat)^\sharp$, but the diagram is not supposed to imply that $1$ or $n-1$ belongs to $A$.},
where
\begin{align*}
\be = \left( \begin{array}{c|c|c|c|c} 
1 \ & \ \cdots \ & \ r-1 \ & \ n \ & \\ \cline{5-5}
1 \ & \ \cdots \ & \ r-1 \ & \ n \ & \ r,n-1  \ \ 
\end{array} \!\!\! \right)
\AND
\ga = \left( \begin{array}{c|c|c|c} 
1 \ & \ \cdots \ & \ r \ & \ n-1,n \ \ \\ \cline{4-4}
1 \ & \ \cdots \ & \ r \ & 
\end{array} \!\!\! \right).
\end{align*}
Since $\be,\ga\in E(D_r(\PB_n))$, the inductive step is complete in this case. \epf

\begin{figure}[h]
\begin{center}
\begin{tikzpicture}[xscale=.56,yscale=0.56]
\bluetrap10409212
\uvs{1,9,10}
\lvs{1,4,5,6,8,9,10}
\stlineds{1/1,9/4}
\stlines{10/5}
\dvertlabels{1/1,5/r,10/n}
\ldottedsms{6/8}
\udotteds{1/9}
\ldotteds{1/4}
\draw[|-|] (.2,0)--(.2,2);
\draw(.2,1)node[left]{$\lam_A$};
\begin{scope}[shift={(0,-4.5)}]	
\bluetrap10409212
\uvs{1,9,10}
\lvs{1,4,5,6,8,9,10}
\stlineds{1/1,9/4}
\stlines{10/10}
\ldottedsms{6/8}
\udotteds{1/9}
\ldotteds{1/4}
\draw[|-|] (.2,0)--(.2,2);
\draw(.2,1)node[left]{$(\lam_A^\flat)^\sharp$};
\end{scope}
\begin{scope}[shift={(0,-6.5)}]	
\lvs{1,4,5,6,8,9,10}
\stlines{1/1,4/4,10/10}
\darc59
\ldottedsms{6/8}
\ldottedsms{1/4}
\draw[|-] (.2,0)--(.2,2);
\draw(.2,1)node[left]{$\be$};
\end{scope}
\begin{scope}[shift={(0,-8.5)}]	
\lvs{1,4,5,6,8,9,10}
\stlines{1/1,4/4,5/5}
\uarc9{10}
\dvertlabels{1/1,5/r,10/n}
\ldottedsms{6/8}
\ldottedsms{1/4}
\draw[|-] (.2,0)--(.2,2);
\draw(.2,1)node[left]{$\ga$};
\end{scope}
\end{tikzpicture}
\end{center}
\vspace{-5mm}
\caption{Diagrammatic verification of the equation $\lam_A = (\lam_A^\flat)^\sharp \cdot \be\ga$ from the proof of Lemma~\ref{lem:lamAEDr}.  
See the text for a full explanation.}
\label{fig:lamAEDr}
\end{figure}


The next result now follows from Lemmas \ref{lem:normalformPBn}, \ref{lem:HrEDr} and \ref{lem:lamAEDr}.

\ms
\begin{cor}\label{cor:DrEDr}
If $0\leq r\leq n-2$, then $D_r\sub\la E(D_r)\ra$. \epfres
\end{cor}


Of course it follows quickly from Corollary \ref{cor:DrEDr} that for $0\leq r\leq n-2$, $I_r=D_0\cup\cdots\cup D_r$ is generated by its idempotents.  But we will not state this formally yet, since we will soon prove a stronger result (see Corollary \ref{cor:IdDrDr-1}).

\ms
\begin{lemma}\label{lem:DsDr}
Let $0\leq s\leq r\leq n-2$ with $r\equiv s\pmod2$.  Then $D_s\sub\la D_r\ra$.
\end{lemma}

\begin{figure}[h]
\begin{center}
\begin{tikzpicture}[xscale=.56,yscale=0.56]
\bluebox9{1.6}{13}2
\stlines{1/1,3/3,6/4,8/6}
\uarc45
\vertlabelsh{1/{\phantom{{}_1}i_1},3/{\phantom{{}_s}i_s},4/u,5/v,6/{\phantom{{}_{s+1}}i_{s+1}},8/{\phantom{{}_{r-2}}i_{r-2}}}
\dvertlabels{1/1,3/s,8/r,13/n}
\udottedsms{1/3,6/8}
\ldottedsms{1/3,4/6,9/11}
\udotteds{9/13}
\draw[|-|] (.2,0)--(.2,2);
\draw(.2,1)node[left]{$\be$};
\lvreds{1,3,4,6,7,8,9,11,12,13}
\uvs{1,3,4,5,6,8,9,13}
\begin{scope}[shift={(19,0)}]	
\bluebox5{1.6}82
\bluebox508{.4}
\stlines{1/1,3/3}
\vertlabelsh{1/{\phantom{{}_1}i_1},3/{\phantom{}i_{r-1}},4/u,5/x}
\dvertlabels{1/{\phantom{{}_1}j_1},3/{\phantom{}j_{r-1}},4/v,5/y}
\udottedsms{1/3}
\ldottedsms{1/3,6/8}
\udotteds{5/8}
\ldotteds{5/8}
\draw[|-|] (.2,0)--(.2,2);
\draw(.2,1)node[left]{$\al$};
\uvs{1,3,4,5,6,8}
\lvreds{1,3,4,5,6,8}
\end{scope}
\begin{scope}[shift={(0,-5)}]	
\bluebox9{1.6}{13}2
\stlines{1/1,3/3,4/4,5/5,6/6,8/8}
\vertlabelsh{1/{\phantom{{}_1}i_1},3/{\phantom{{}_s}i_s},4/u,5/v,6/{\phantom{{}_{s+1}}i_{s+1}},8/{\phantom{{}_{r-2}}i_{r-2}}}
\udottedsms{1/3,6/8}
\ldottedsms{1/3,6/8,9/11}
\udotteds{9/13}
\draw[|-|] (.2,0)--(.2,2);
\draw(.2,1)node[left]{$\be_1$};
\lvreds{1,3,4,6,5,8,9,11,12,13}
\uvs{1,3,4,5,6,8,9,13}
\end{scope}
\begin{scope}[shift={(19,-5)}]	
\bluebox5{1.6}82
\stlines{1/1,3/3,4/4}
\vertlabelsh{1/{\phantom{{}_1}i_1},3/{\phantom{}i_{r-1}},4/u,5/x}
\udottedsms{1/3}
\ldottedsms{1/3,6/8}
\udotteds{5/8}
\draw[|-|] (.2,0)--(.2,2);
\draw(.2,1)node[left]{$\be$};
\uvs{1,3,4,5,6,8}
\lvs{1,3,4,5,6,8}
\end{scope}
\begin{scope}[shift={(0,-7)}]	
\stlines{1/1,3/3,6/4,8/6,12/12,13/13}
\uarc45
\dvertlabels{1/1,3/s,8/r,13/n}
\udottedsms{1/3,6/8}
\ldottedsms{1/3,4/6,9/11}
\draw[|-] (.2,0)--(.2,2);
\draw(.2,1)node[left]{$\be_2$};
\lvreds{1,3,4,6,7,8,9,11,12,13}
\uvreds{1,3,4,6,5,8,9,11,12,13}
\end{scope}
\begin{scope}[shift={(19,-7)}]	
\stlines{1/1,3/3,5/5}
\udottedsms{1/3}
\ldottedsms{1/3,6/8}
\draw[|-] (.2,0)--(.2,2);
\draw(.2,1)node[left]{$\ga$};
\draw(.2,1)node[left]{$\phantom{\be_2}$};
\uvs{1,3,4,5,6,8}
\lvreds{1,3,4,5,6,8}
\end{scope}
\begin{scope}[shift={(19,-9)}]	
\bluebox508{.4}
\stlines{1/1,3/3,4/4}
\dvertlabels{1/{\phantom{{}_1}j_1},3/{\phantom{}j_{r-1}},4/v,5/y}
\udottedsms{1/3}
\ldottedsms{1/3}
\ldotteds{5/8}
\draw[|-] (.2,0)--(.2,2);
\draw(.2,1)node[left]{$\de$};
\uvreds{1,3,4,5,6,8}
\lvreds{1,3,4,5,6,8}
\end{scope}
\end{tikzpicture}
\end{center}
\vspace{-5mm}
\caption{Diagrammatic verification of the equations $\be=\be_1\be_2$ and $\al=\be\ga\de$ from the proofs of Lemmas~\ref{lem:DsDr} and \ref{lem:singleton} --- left and right, respectively.  Note the ordering on vertices.  In the left diagram, red vertices are arranged in the order $1,\ldots,n$, and black vertices in the order $i_1,\ldots,i_s,u,v,i_{s+1},\ldots,i_{r-2}$ (followed by an arbitrary but fixed ordering on $\{i_1,\ldots,i_{r-2},u,v\}^c$).  In the right diagram, black vertices are arranged in the order $i_1,\ldots,i_{r-1},u,x$ (followed by an arbitrary but fixed ordering on $\{i_1,\ldots,i_{r-1},u,x\}^c$), and red vertices in the order $j_1,\ldots,j_{r-1},v,y$ (followed by an arbitrary but fixed ordering on $\{j_1,\ldots,j_{r-1},v,y\}^c$).  See the text for a full explanation.}
\label{fig:DsDr}
\end{figure}
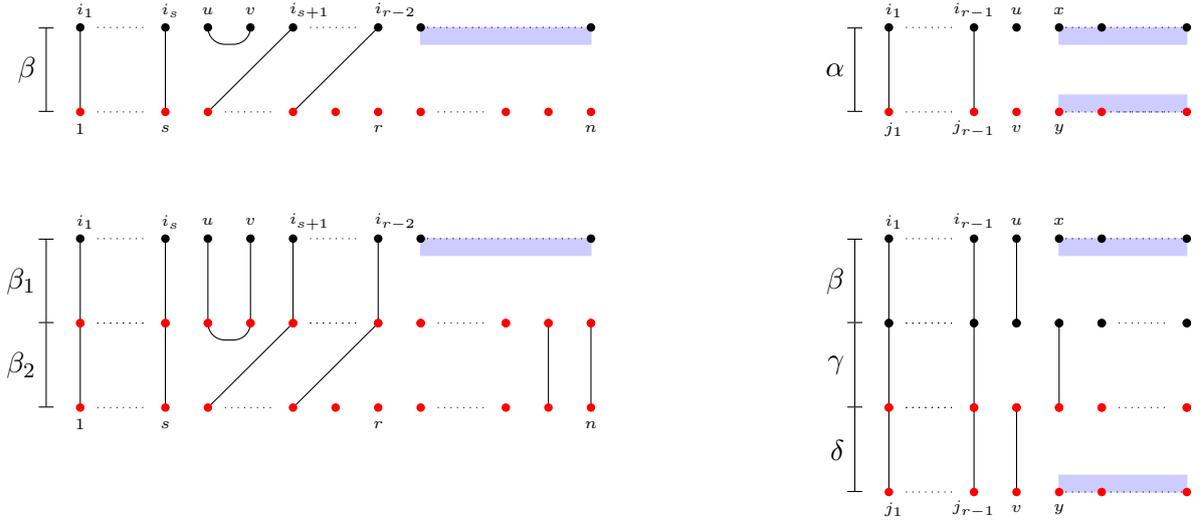

\pf By induction, it suffices to show that $D_{r-2}\sub\la D_r\ra$ for all $2\leq r\leq n-2$.  So suppose $2\leq r\leq n-2$, and let $\al\in D_{r-2}$ be arbitrary.  Let $\dom(\al)=\{i_1<\cdots<i_{r-2}\}$ and $\codom(\al)=\{j_1<\cdots<j_{r-2}\}$.  So then we may write
\[
\al = \left( \begin{array}{c|c|c|c} 
i_1 \ & \ \cdots \ & \ i_{r-2} \ & \ A_k\ \ \\ \cline{4-4}
j_{1\pi} \ & \ \cdots \ & \ j_{(r-2)\pi} \ & \ B_l \ \
\end{array} \!\!\! \right)_{k\in K,\ l\in L},
\]
for some $\pi\in\S_{r-2}$.  As in the proof of Lemma \ref{lem:normalformPBn}, we have $\al=\be\ga\de$ for some $\ga\in\S_{[r-2]}$, where
\[
\be = \left( \begin{array}{c|c|c|c} 
i_1 \ & \ \cdots \ & \ i_{r-2} \ & \ A_k\ \ \\ \cline{4-4}
1 \ & \ \cdots \ & \ r-2 \ & 
\end{array} \!\!\! \right)_{k\in K}
\AND
\de = \left( \begin{array}{c|c|c|c} 
1 \ & \ \cdots \ & \ r-2 \ &  \\ \cline{4-4}
j_1 \ & \ \cdots \ & \ j_{r-2} \ & \ B_l \ \
\end{array} \!\!\! \right)_{l\in L}.
\]
Now, $\ga\in D_{r-2}(\I_n)\sub\la D_r(\I_n)\ra\sub\la D_r\ra$, by Proposition \ref{prop:ZF}, so it remains to show that $\be,\de\in\la D_r\ra$.  By symmetry, we may just do this for $\be$.  Now, if $\be$ has no upper hooks, then $\be\in D_{r-2}(\POI_n)\sub\la D_r(\POI_n)\ra\sub\la D_r\ra$, by Proposition \ref{prop:ZF}, so suppose $\be$ has an upper hook $\{u,v\}$, where $u\in\bn$ is minimal with the property of belonging to an upper hook.  Let $0\leq s\leq r-2$ be such that $i_s<u<i_{s+1}$, where we also define $i_0=0$ and $i_{r-1}=n+1$.  Let $h\in K$ be such that $A_h=\{u,v\}$.  We see then that $\be=\be_1\be_2$ (see Figure \ref{fig:DsDr}), where
\begin{align*}
\be_1 &= \left( \begin{array}{c|c|c|c|c|c|c|c|c} 
i_1 \ & \ \cdots \ & \ i_s \ & \ u \ & \ v \ & \ i_{s+1} \ & \ \cdots \ & \ i_{r-2} \ & \  A_k\ \ \\ \cline{9-9}
1 \ & \ \cdots \ & \ s \ & \ s+1 \ & \ s+2 \ & \ s+3 \ & \ \cdots \ & r \ & 
\end{array} \!\!\! \right)_{k\in K\sm\{h\}},
\\
\be_2 &= \left( \begin{array}{c|c|c|c|c|c|c|c|c} 
1 \ & \ \cdots \ & \ s \ & \ s+3 \ & \ \cdots \ & \ r \ & \ n-1 \ & \ n \ & \  s+1,s+2 \ \ \\ \cline{9-9}
1 \ & \ \cdots \ & \ s \ & \ s+1 \ & \ \cdots \ & \ r-2 \ & \ n-1 \ & \ n \ & 
\end{array} \!\!\! \right).
\end{align*}
Since $\be_1,\be_2\in D_r$, the proof is complete. \epf

\ms
\begin{rem}\label{rem:DsDrM}
The previous proof also shows that $D_s(\M_n)\sub\la D_r(\M_n)\ra$ for $0\leq s\leq r\leq n-2$ with $r\equiv s\pmod2$.  Indeed, if $\al\in D_{r-2}(\M_n)$, then $\al=\be\de$, where $\be,\de$ are as constructed in the above proof, and we see that the planarity of $\al$ entails that $\be,\de,\be_1,\be_2$ are all planar.  (In the case of $\be_1$, planarity of $\al$ gives $i_s<u<v<i_{s+1}$.)
\end{rem}

\ms
\begin{cor}\label{cor:IdDrDr-1}
If $1\leq r\leq n-2$, then $I_r=\la D_r\cup D_{r-1}\ra = \la E(D_r\cup D_{r-1}) \ra$.
\end{cor}

\pf Lemma \ref{lem:DsDr} gives $\la D_r\ra\supseteq D_r\cup D_{r-2}\cup\cdots$ and $\la D_{r-1}\ra\supseteq D_{r-1}\cup D_{r-3}\cup\cdots$, from which it follows that $\la D_r\cup D_{r-1}\ra\supseteq I_r$.  The reverse containment is obvious.  The rest follows from Corollary~\ref{cor:DrEDr}.~\epf

We have shown that for $0\leq r\leq n-2$, the ideal $I_r$ is generated by the idempotents in its top two $\D$-classes.  We will soon calculate the rank and idempotent rank for these ideals (see Theorem \ref{thm:IrPBn}); in particular, we will see that if $r=0$ or $r\equiv n\pmod2$, only the idempotents of rank $r$ are needed to generate $I_r$.  But first we prove some intermediate results.


\ms
\begin{lemma}\label{lem:singleton}
Let $1\leq r\leq n-1$, and suppose $\al\in D_{r-1}$ has an upper singleton block and a lower singleton block.  Then $\al\in\la D_r\ra$.
\end{lemma}

\pf Write 
\[
\al = \left( \begin{array}{c|c|c|c} 
i_1 \ & \ \cdots \ & \ i_{r-1} \ & \ A_k \ \ \\ \cline{4-4}
j_1 \ & \ \cdots \ & \ j_{r-1} \ & \ B_l \ \
\end{array} \!\!\! \right)_{k\in K,\ l\in L},
\]
and suppose $\{u\}$ and $\{v'\}$ are blocks of $\al$.  Since $r\leq n-1$, we may choose further elements $x\in\{i_1,\ldots,i_{r-1},u\}^c$ and $y\in\{j_1,\ldots,j_{r-1},v\}^c$.  Then $\al=\be\ga\de$ (see Figure \ref{fig:DsDr}), where
\begin{gather*}
\be = \left( \begin{array}{c|c|c|c|c} 
i_1 \ & \ \cdots \ & \ i_{r-1} \ & \ u \ & \ A_k \ \ \\ \cline{5-5}
i_1 \ & \ \cdots \ & \ i_{r-1} \ & \ u \ & \  \ \
\end{array} \!\!\! \right)_{k\in K}
\COMMA
\ga = \left( \begin{array}{c|c|c|c} 
i_1 \ & \ \cdots \ & \ i_{r-1} \ & \ x \ \ \\ 
j_1 \ & \ \cdots \ & \ j_{r-1} \ & \ y \ \
\end{array} \!\!\! \right)
,
\\
\epfreseq
\de = \left( \begin{array}{c|c|c|c|c} 
j_1 \ & \ \cdots \ & \ j_{r-1} \ & \ v \ & \  \ \ \\ \cline{5-5}
j_1 \ & \ \cdots \ & \ j_{r-1} \ & \ v \ & \ B_l \ \
\end{array} \!\!\! \right)_{l\in L}
.
\end{gather*}

\ms
\begin{cor}\label{cor:r0nmod2}
If $r=0$ or $1\leq r\leq n-2$ satisfies $r\equiv n\pmod2$, then $I_r = \la D_r\ra = \la E(D_r) \ra$.
\end{cor}

\pf If $r=0$, then $I_0=D_0=\la D_0\ra$.
Next, suppose $1\leq r\leq n-2$ satisfies $r\equiv n\pmod2$.
Since $n-(r-1)$ is odd, every element of $D_{r-1}$ has an upper singleton block and a lower singleton block, so $D_{r-1}\sub\la D_r\ra$, by Lemma~\ref{lem:singleton}.  The result now follows from Corollary \ref{cor:IdDrDr-1}. \epf

\ms
\begin{lemma}\label{lem:upperandlowersingleton}
Suppose $0\leq r\leq n-2$ satisfies $r\not\equiv n\pmod2$, and let $\al\in D_s$ where $s\leq r$.  Then $\al\in\la D_r\ra$ if and only if $r\equiv s\pmod2$ or $\al$ has an upper singleton block and a lower singleton block.
\end{lemma}

\pf We have already seen that $\la D_r\ra\supseteq D_r\cup D_{r-2}\cup\cdots$.  Since $n-r$ is odd, every element of $D_r$, and hence every element of $\la D_r\ra$, has at least one upper singleton block and at least one lower singleton block.  Conversely, if $r\not\equiv s\pmod2$, then any element of $D_s$ with an upper singleton block and a lower singleton block belongs to $\la D_{s+1}\ra\sub\la D_r\ra$, by Lemmas \ref{lem:singleton} and \ref{lem:DsDr}. \epf


For the statement of the next result, recall that $k!!=0$ if $k$ is even, and that ${n\choose-1}=0$.

\ms
\begin{thm}\label{thm:IrPBn}
If $0\leq r\leq n-2$, then $I_r=I_r(\PB_n)$ is idempotent-generated, and
\[
\rank(I_r)=\idrank(I_r) = {n\choose r-1}\cdot(n-r)!!+{n \choose r} \cdot a(n-r),
\]
where the numbers $a(m)$ satisfy 
\[
a(0)=a(1)=1 \COMMA a(m)=a(m-1)+(m-1)a(m-2) \quad\text{for $m\geq2$.}
\]
\end{thm}

\pf If $r=0$ or $1\leq r\leq n-2$ satisfies $r\equiv n\pmod2$, then the first term in the sum is equal to $0$, and the result follows immediately from Corollary \ref{cor:r0nmod2}, Lemma \ref{lem:ARSS<D>}(iii) and Proposition \ref{prop:DClassesofKn}(iv).
%
%
For the remainder of the proof, suppose $1\leq r\leq n-2$ satisfies $r\not\equiv n\pmod2$.  Let $\Si$ be an arbitrary generating set for $I_r$.  Clearly $\Si\cap D_r$ is a generating set for $\la D_r\ra$, so it follows from Lemma \ref{lem:ARSS<D>}(i) and Proposition~\ref{prop:DClassesofKn}(iv) that 
\begin{equation}\label{eq1}
|\Si\cap D_r| \geq |P(D_r)| = {n \choose r} \cdot a(n-r).
\end{equation}
%
Now suppose $\al\in P(D_{r-1}(\B_n))$, noting that $\al$ has no upper singleton blocks.  (Recall that $\B_n\sub\PB_n$ denotes the Brauer monoid.)  
Consider an expression $\al=\be_1\cdots\be_k$, where $\be_1,\ldots,\be_k\in \Si$.  Since each element of $D_r$ has an upper singleton block, it follows that $\be_1\in D_{r-1}$.  In particular, $\be_1\D\al=\be_1(\be_2\cdots\be_k)$.  By stability, it follows that $\be_1\R\be_1(\be_2\cdots\be_k)=\al$.  So $\Si$ contains at least one element from the $\R$-class of every such $\al\in P(D_{r-1}(\B_n))$.  Since projections belong to distinct $\R$-classes, it follows from Proposition~\ref{prop:DClassesofKn}(ii) that
\begin{equation}\label{eq2}
|\Si\cap D_{r-1}| \geq |P(D_{r-1}(\B_n))|={n\choose r-1}\cdot(n-r)!!.
\end{equation}
Adding equations \eqref{eq1} and \eqref{eq2} shows that the stated value of $\rank(I_r)$ is a lower bound.
%
%

To complete the proof, it suffices to give an idempotent generating set for $I_r$ of the stated size.  With this in mind, let $\Ga=P(D_r)\cup P(D_{r-1}(\B_n))$.  
By Lemma \ref{lem:singleton},
$\la D_r\ra$ contains every projection from $D_{r-1}$ having at least one singleton block.  Any other projection from $D_{r-1}$ belongs to $\B_n$, so it follows that $\la\Ga\ra$ contains $P(D_r\cup D_{r-1})$, and we have already seen that $I_r=\la E(D_r\cup D_{r-1})\ra=\la P(D_r\cup D_{r-1})\ra$, so it follows that $I_r=\la\Ga\ra$.  It follows from Proposition \ref{prop:DClassesofKn}(ii) and (iv) that $\Ga$ has the required size. \epf

We now consider the top two ideals $I_{n-1}=\PB_n\sm\S_n$ and $I_n=\PB_n$, which were not covered in the previous result.  Unlike the other ideals, these are not idempotent-generated (unless $n=1$).

\ms
\begin{prop}\label{prop:In-1}
The singular ideal $I_{n-1}(\PB_n)=\PB_n\sm\S_n$ is not idempotent-generated for $n\geq2$.  We have 
\[
\PB_n\sm\S_n = \la D_{n-1}\cup D_{n-2}\ra = \la D_{n-1}(\I_n)\cup D_{n-2}(\B_n)\ra
\AND
\rank(\PB_n\sm\S_n)=\begin{cases}
{n+1\choose2} &\text{if $n\leq3$}\\
1+{n+1\choose2} &\text{if $n\geq4$.}
\end{cases}
\]
\end{prop}

\pf It is obvious that $I_{n-1}=\PB_n\sm\S_n$ is not idempotent-generated if $n\geq2$.  Indeed, the idempotents from $D_{n-1}(\PB_n)=D_{n-1}(\I_n)$ all commute, so no non-idempotent from $D_{n-1}$ can be a product of idempotents.  It is also obvious that $\rank(\PB_0\sm\S_0)=0$ and $\rank(\PB_1\sm\S_1)=1$, so suppose $n\geq2$ for the remainder of the proof.

Suppose $I_{n-1}=\la\Si\ra$.  Since $D_{n-1}(\I_n)=D_{n-1}(\PB_n)\sub I_{n-1}$, and since $\InSn=\la D_{n-1}(\I_n)\ra$, it follows that $\Si$ contains a generating set $\Si_1$ for $\InSn$.  Now,
\[
\rank(\InSn)=\rho+n \qquad\text{where}\qquad
\rho=\begin{cases}
0 &\text{if $n\leq3$}\\
1 &\text{if $n\geq4$.}
\end{cases}
\]
Indeed, this follows quickly for $n=2,3$ from the easily checked\footnote{by hand or by GAP \cite{GAP}} facts that 
\[
\I_2\sm\S_2 = \Big\la \smcustpartn{1,2}{1,2}{\stlines{1/2}}, \smcustpartn{1,2}{1,2}{\stlines{2/1}} \Big\ra
\AND
\I_3\sm\S_3 = \Big\la \smcustpartn{1,2,3}{1,2,3}{\stlines{1/2,3/3}}, \smcustpartn{1,2,3}{1,2,3}{\stlines{1/1,2/3}}, \smcustpartn{1,2,3}{1,2,3}{\stlines{2/2,3/1}} \Big\ra,
\]
and from \cite[Theorem 3.7]{Gomes1987} for $n\geq4$.  (Note that \cite[Theorem 3.7]{Gomes1987} incorrectly states that $\rank(\I_3\sm\S_3)=4$.)  In particular, $|\Si_1|\geq\rho+1$.
Next, we claim that $|\Si\sm\Si_1|\geq{n\choose2}$.  
With this in mind, let $\oijn$ and consider an expression $\tau_{ij}=\al_1\cdots\al_k$, where $\al_1,\ldots,\al_k\in\Si$.  Now $n-2=\rank(\tau_{ij})\leq\rank(\al_1)\leq n-1$, but we could not have $\rank(\al_1)=n-1$ or else $\al_1$ (and hence also $\tau_{ij}$) would then have a singleton upper block, so it follows that $\rank(\al_1)=n-2$.  So $\al_1\D\tau_{ij}$, and the usual stability argument allows us to deduce that $\al_1\R\tau_{ij}$.  
%
In particular, $\Si\sm\Si_1$ contains an element from the $\R$-class of $\tau_{ij}$ for each $\oijn$,
and it follows that $|\Si\sm\Si_1|\geq{n\choose2}$, as claimed.  We then have $|\Si|=|\Si_1|+|\Si\sm\Si_1|\geq (\rho+n)+{n\choose2}=\rho+{n+1\choose2}$.
Since $\Si$ was an arbitrary generating set, it follows that $\rank(I_{n-1})\geq \rho+{n+1\choose 2}$.

Conversely, suppose $\InSn=\la\Ga\ra$ with $|\Ga|=\rho+n=\rank(\InSn)$, and let $\Si=\Ga\cup T$, where $T=P(D_{r-2}(\B_n))=\set{\tau_{ij}}{\oijn}$.  The proof will be complete if we can show that $I_{n-1}=\la\Si\ra$.  Now $D_{n-1}\sub\la\Ga\ra\sub\la\Si\ra$.  Also, $P(D_{n-2})=T\cup\set{\ve_i\ve_j}{\oijn}$ is clearly contained in $\la\Si\ra$, so 
$\la\Si\ra\supseteq\la P(D_{n-2})\ra=I_{n-2}$ by Corollary \ref{cor:r0nmod2}.  Since $I_{n-1}=D_{n-1}\cup I_{n-2}$, the proof is complete. \epf




\ms
\begin{prop}\label{prop:rankPBn}
The semigroup $\PB_n$ is not idempotent-generated for $n\geq2$.  We have 
\[
\rank(\PB_n)=\begin{cases}
n+1 &\text{if $n\leq2$}\\
4 &\text{if $n\geq3$.}
\end{cases}
\]
\end{prop}

\pf The statement is obvious for $n\leq1$, 
so suppose $n\geq2$.  Again, it is clear that $\PB_n$ is not idempotent-generated.  Since $\rank(\S_2)=1$ and $\rank(\S_n)=2$ for $n\geq3$, it suffices to show that $\rank(\PB_n)=2+\rank(\S_n)$.  Since $\PB_n\sm\S_n$ is an ideal, any generating set must contain a generating set for $\S_n$, so it is enough to show that $\rank(\PB_n:\S_n)=2$.

Now suppose $\PB_n=\la\S_n\cup\Si\ra$.  By considering an expression $\ve_1=\al_1\cdots\al_k$, where $\al_1,\ldots,\al_k\in\S_n\cup\Si$, it is easy to see that the first factor not from $\S_n$, say $\al_i$, must have rank $n-1$.  But then $\la\S_n\cup\{\al_i\}\ra=\I_n\not=\PB_n$ (see \cite[Theorem 3.1]{Gomes1987}).
It follows that $|\Si|\geq2$.  The proof will be complete if we can show that $\PB_n$ is generated by $\S_n\cup\{\ve_1,\tau_{12}\}$.  But $\la\S_n\cup\{\ve_1\}\ra=\I_n$ and $\la\S_n\cup\{\tau_{12}\}\ra=\B_n$ (see \cite[Lemma 3.5]{KM2006}), so Proposition~\ref{prop:In-1} gives $\PB_n=\la\S_n\cup\{\ve_1,\tau_{12}\}\ra$. \epf

\ms
\begin{rem}
The fact that $\PB_n=\la\S_n\cup\{\ve_1,\tau_{12}\}\ra$ was also noted in \cite{KM2006}, where a presentation for $\PB_n$ was obtained (with respect to a different generating set).  See also \cite[Lemma 11]{Maz1998}.
\end{rem}

We conclude this section by describing the idempotent-generated subsemigroup $\bbE(\PB_n)=\la E(\PB_n)\ra$.  If $n\leq1$, then $\bbE(\PB_n)=E(\PB_n)=\PB_n$, so we just consider the case $n\geq2$.

\ms
\begin{thm}\label{thm:EPBn}
If $n\geq2$, then $\bbE(\PB_n) = E(D_n\cup D_{n-1})\cup I_{n-2}$, and 
\[
\rank(\bbE(\PB_n))=\idrank(\bbE(\PB_n)) = 1+{n+1\choose2}.
\]
\end{thm}

\pf Since $I_{n-2}$ is idempotent-generated, it is clear that $E(D_n\cup D_{n-1})\cup I_{n-2}\sub\bbE(\PB_n)$.  Conversely, suppose $\al_1,\ldots,\al_k\in E(\PB_n)$, and put $\be=\al_1\cdots\al_k$.  If $\be\not\in I_{n-2}$, then all of $\al_1,\ldots,\al_k$ belong to $E(D_n\cup D_{n-1})=\{1\}\cup\{\ve_1,\ldots,\ve_n\}$, from which it follows quickly that $\be\in E(D_n\cup E_{n-1})$.  We have shown that $\bbE(\PB_n) = E(D_n\cup D_{n-1})\cup I_{n-2}$.

Next, suppose $\bbE(\PB_n)=\la\Si\ra$, where $\Si\sub\bbE(\PB_n)$ is arbitrary.  Clearly, $\Si$ must contain $E(D_n\cup D_{n-1})$.  As in the proof of Proposition \ref{prop:In-1}, we may show that $\Si$ contains at least ${n\choose2}$ elements of rank $n-2$.
It follows that $|\Si|\geq(1+n)+{n\choose2}=1+{n+1\choose2}$, giving $\rank(\bbE(\PB_n))\geq 1+{n+1\choose2}$.

To complete the proof, it suffices to find an idempotent generating set of size $1+{n+1\choose2}$.  Given that $\bbE(\PB_n)=E(D_n\cup D_{n-1})\cup I_{n-2}=E(D_n\cup D_{n-1})\cup\la P(D_{n-2})\ra$, it is clear that $\bbE(\PB_n)$ is generated by $E(D_n\cup D_{n-1})\cup T$, where $T=\set{\tau_{ij}}{\oijn}$, since $\la E(D_{n-1})\cup T\ra\supseteq P(D_{n-2})=T\cup\set{\ve_i\ve_j}{\oijn}$.  \epf

\ms
\begin{rem}
From the above proof, it quickly follows that any minimal (in size) idempotent generating set of $\bbE(\PB_n)$ is of the form $E(D_n\cup D_{n-1})\cup\Ga$, where $\Ga$ is a minimal idempotent generating set of $\B_n\sm\S_n$.  The minimal idempotent generating sets of $\B_n\sm\S_n$ were classified (in terms of factorizations of certain \emph{Johnson graphs}) in \cite[Proposition 8.7]{EastGray}, but the enumeration of such generating sets for $n\geq6$ remains an open problem.
\end{rem}

\section{The Motzkin monoid $\M_n$}\label{sect:Mn}

The situation with Motzkin monoid ideals is somewhat different to partial Brauer monoid ideals.  It is still the case that each proper ideal $I_r(\M_n)$ is generated by its top two $\D$-classes (Proposition \ref{prop:IdDrDr-1M}), but the ideals are never generated by a single $\D$-class, apart from $I_0(\M_n)=D_0(\M_n)$.  Idempotent generation is also more subtle, with the ideal $I_r(\M_n)$ being idempotent-generated if and only if $0\leq r<\floorn$ (Proposition \ref{prop:IrIGM}).  We are still able to calculate the rank and (when appropriate) idempotent rank of the ideals (Theorem \ref{thm:IrM}).  We also describe the subsemigroups of $\M_n$ generated by a single $\D$-class (Proposition \ref{prop:multinesting}) and the idempotent-generated subsemigroup of each ideal (Theorems \ref{thm:IGIrM} and \ref{thm:EMn}), which includes the idempotent-generated subsemigroup $\bbE(\M_n)$.  In particular, apart from the $r=n$ case, we see that the rank of $\bbE(I_r(\M_n))$ is equal to the rank of $I_r(\M_n)$, even when $r\geq\floorn$.


In this section, we use the abbreviations
\[
D_r=D_r(\M_n) = \set{\al\in\M_n}{\rank(\al)=r} \AND
I_r=I_r(\M_n) = \set{\al\in\M_n}{\rank(\al)\leq r}.
\]
As noted in Remark \ref{rem:DsDrM}, the next result has already been proved.

\ms
\begin{lemma}\label{lem:DsDrM}
Let $0\leq s\leq r\leq n-2$ with $r\equiv s\pmod2$.  Then $D_s\sub\la D_r\ra$. \epfres
\end{lemma}

Recall the maps $\lam_A,\rho_A$ defined near the beginning of Section \ref{sect:PBn}.


\ms
\begin{prop}\label{prop:IdDrDr-1M}
If $1\leq r\leq n-1$, then $I_r=\la D_r\cup D_{r-1}\ra$. 
%
%
\end{prop}

\pf The $r\leq n-2$ case follows from Lemma \ref{lem:DsDrM}.  Since $I_{n-1}=D_{n-1}\cup I_{n-2}=D_{n-1}\cup\la D_{n-2}\cup D_{n-3}\ra$, it remains to show that $D_{n-3}\sub\la D_{n-1}\cup D_{n-2}\ra$.  So let $\al\in D_{n-3}$ be arbitrary, and write
\[
\al = \left( \begin{array}{c|c|c|c} 
i_1 \ & \ \cdots \ & \ i_{n-3} \ & \ A_k \ \ \\ \cline{4-4}
j_1 \ & \ \cdots \ & \ j_{n-3} \ & \ B_l \ \
\end{array} \!\!\! \right)_{k\in K,\ l\in L},
\]
where $i_1<\cdots<i_{n-3}$ and $j_1<\cdots<j_{n-3}$.
Note that $K$ and $L$ each have size at most $1$.  Again, $\al=\be\ga$, where
\[
\be = \left( \begin{array}{c|c|c|c} 
i_1 \ & \ \cdots \ & \ i_{n-3} \ & \ A_k \ \ \\ \cline{4-4}
1 \ & \ \cdots \ & \ n-3 \ & \  \ \
\end{array} \!\!\! \right)_{k\in K}
\text{ \ and \ }
\ga = \left( \begin{array}{c|c|c|c} 
1 \ & \ \cdots \ & \ n-3 \ & \  \ \ \\ \cline{4-4}
j_1 \ & \ \cdots \ & \ j_{n-3} \ & \ B_l \ \
\end{array} \!\!\! \right)_{l\in L}.
\]
By symmetry, we just need to show that $\be\in\la D_{n-1}\cup D_{n-2}\ra$.  If $K$ is empty, then $\be\in D_{n-3}(\POI_n)\sub\la D_{n-1}(\POI_n)\ra\sub\la D_{n-1}\cup D_{n-2}\ra$, so suppose $K=\{k\}$, and write $A_k=\{x<y\}$.  Also, let $\dom(\al)^c=\dom(\be)^c=\{x,y,z\}$.  We now consider separate cases depending on the ordering of the set $\{x,y,z\}$.

\begin{figure}[h]
\begin{center}
\begin{tikzpicture}[xscale=.56,yscale=0.56]
\uvs{1,3,4,5,7,8,9,10,12}
\lvs{1,3,4,6,7,9,10,12}
\stlines{1/1,3/3,5/4,7/6,10/7,12/9}
\uarc89
\vertlabelsh{1/1,4/z,8/x,9/{\lower1.3 ex\hbox{$y$}},12/n}
\udottedsms{1/3,5/7,10/12}
\ldottedsms{1/3,4/6,7/9,10/12}
\draw[|-|] (.2,0)--(.2,2);
\draw(.2,1)node[left]{$\be$};
\begin{scope}[shift={(19,0)}]	
\uvs{1,3,4,5,6,7,9}
\lvs{1,3,4,6,7,8,9}
\stlines{1/1,3/3,7/4,9/6}
\uarc46
\vertlabelsh{1/1,4/x,5/z,6/{\lower1.3 ex\hbox{$y$}},9/n}
\udottedsms{1/3,7/9}
\ldottedsms{1/3,4/6}
\draw[|-|] (.2,0)--(.2,2);
\draw(.2,1)node[left]{$\be$};
\end{scope}
\begin{scope}[shift={(0,-4.5)}]	
\uvs{1,3,4,5,7,8,9,10,12}
\lvs{1,3,4,7,8,9,10,12}
\stlines{1/1,3/3,5/5,7/7,10/10,12/12}
\uarc89
\darc89
\vertlabelsh{1/1,4/z,8/x,9/{\lower1.3 ex\hbox{$y$}},12/n}
\udottedsms{1/3,5/7,10/12}
\ldottedsms{1/3,5/7,10/12}
\draw[|-|] (.2,0)--(.2,2);
\draw(.2,1)node[left]{$\tau_{xy}$};
\end{scope}
\begin{scope}[shift={(19,-4.5)}]	
\uvs{1,3,4,5,6,7,9}
\lvs{1,3,4,5,6,8,9}
\stlines{1/1,3/3,4/4,6/5,7/6,9/8}
\vertlabelsh{1/1,4/x,5/z,6/{\lower1.3 ex\hbox{$y$}},9/n}
\udottedsms{1/3,7/9}
\ldottedsms{1/3,6/8}
\draw[|-|] (.2,0)--(.2,2);
\draw(.2,1)node[left]{$\lam_{\{z\}^c}$};
\end{scope}
\begin{scope}[shift={(0,-6.5)}]	
\uvs{1,3,4,5,7,8,9,10,12}
\lvs{1,3,4,6,7,9,10,12}
\stlines{1/1,3/3,5/4,7/6,10/7,12/9}
\udottedsms{1/3,5/7,10/12}
\ldottedsms{1/3,4/6,7/9,10/12}
\draw[|-] (.2,0)--(.2,2);
\draw(.2,1)node[left]{$\lam_{\{x,y,z\}^c}$};
\end{scope}
\begin{scope}[shift={(19,-6.5)}]	
\uvs{1,3,4,5,6,8,9}
\lvs{1,3,4,5,6,8,9}
\stlines{1/1,3/3,6/6,8/8,9/9}
\uarc45
\darc45
\udottedsms{1/3,6/8}
\ldottedsms{1/3,6/8}
\draw[|-] (.2,0)--(.2,2);
\draw(.2,1)node[left]{$\tau_{xz}$};
\end{scope}
\begin{scope}[shift={(19,-8.5)}]	
\uvs{1,3,4,5,6,8,9}
\lvs{1,3,4,6,7,8,9}
\stlines{1/1,3/3,6/4,8/6,9/7}
\udottedsms{1/3,6/8}
\ldottedsms{1/3,4/6}
\draw[|-] (.2,0)--(.2,2);
\draw(.2,1)node[left]{$\lam_{\{x,z\}^c}$};
\end{scope}
\end{tikzpicture}
\end{center}
\vspace{-5mm}
\caption{Diagrammatic verification of the equations $\be=\tau_{xy}\cdot\lam_{\{x,y,z\}^c}$ and $\be=\lam_{\{z\}^c}\cdot\tau_{xz}\cdot\lam_{\{x,z\}^c}$ from the proof of Proposition~\ref{prop:IdDrDr-1M} --- left and right, respectively.  See the text for a full explanation.}
\label{fig:IdDrDr-1M}
\end{figure}
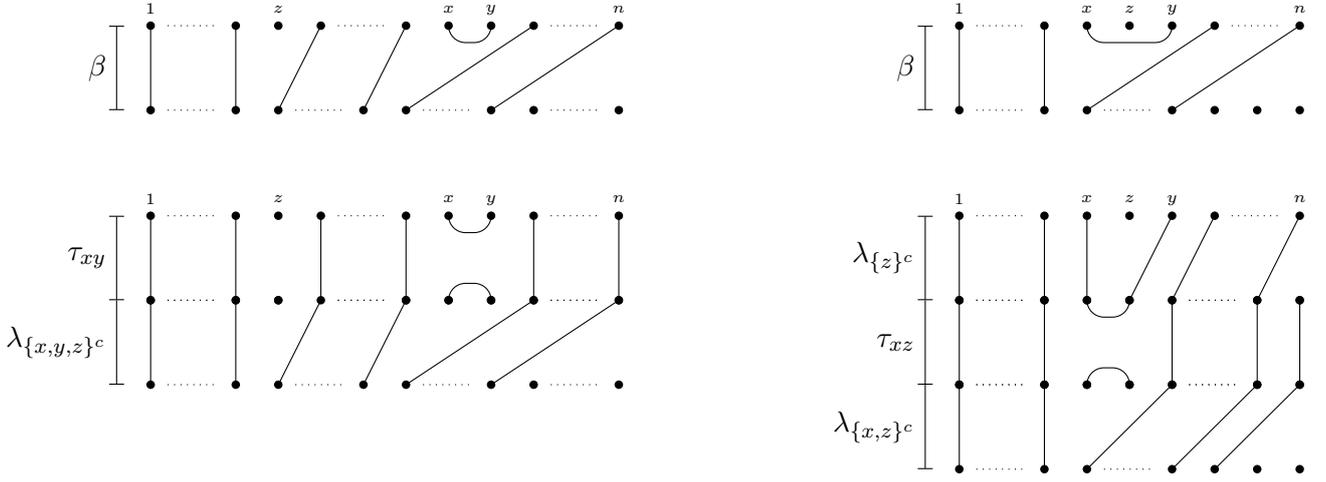

{\bf Case 1.}  Suppose first that $z<x$ or $y<z$.  By planarity, it follows that $y=x+1$.  In either case, we have $\be=\tau_{xy}\cdot\lam_{\{x,y,z\}^c}$ (see Figure \ref{fig:IdDrDr-1M}, which just pictures the $z<x$ case).
Since $\tau_{xy}\in D_{n-2}$ and $\lam_{\{x,y,z\}^c}\in D_{n-3}(\POI_n)\sub\la D_{n-1}\cup D_{n-2}\ra$, the proof is complete in this case.

{\bf Case 2.}  Now suppose that $x<z<y$.  This time, planarity forces $z=x+1$ and $y=x+2$.  In this case, we have $\be=\lam_{\{z\}^c}\cdot\tau_{xz}\cdot\lam_{\{x,z\}^c}$ (again, see Figure \ref{fig:IdDrDr-1M}), with all factors belonging to $D_{n-1}\cup D_{n-2}$. \epf

As noted above, an ideal $I_r=I_r(\M_n)$ of $\M_n$ is never generated by its top $\D$-class $D_r=D_r(\M_n)$ (apart from the trivial case in which $r=0$), as we now describe.  Recall from Lemma \ref{lem:DsDrM} that $\la D_r\ra$ contains all of $D_s$ for any $0\leq s\leq r$ with $r\equiv s\pmod2$.  So, to fully describe $\la D_r\ra$, we need to characterise the elements of $D_s\cap\la D_r\ra$ with $0\leq s\leq r$ and $r\not\equiv s\pmod2$.  We begin with a simple lemma.


\ms
\begin{lemma}\label{lem:<DrMn>proj}
Let $0\leq r\leq n-1$, and let $\al\in \M_n$.  Then $\al\in\la D_r\ra$ if and only if $\al\al^*,\al^*\al\in\la D_r\ra$.
\end{lemma}

\pf If $\al\in\la D_r\ra$, then $\al=\be_1\cdots\be_k$ for some $\be_1,\ldots,\be_k\in D_r$; but then $\al^*=\be_k^*\cdots\be_1^*\in\la D_r\ra$, so that $\al\al^*,\al^*\al\in\la D_r\ra$.  
Conversely, suppose $\al\al^*,\al^*\al\in\la D_r\ra$, and put $s=\rank(\al)=\rank(\al\al^*)$, noting that $s\leq r\leq n-1$.  Let $\be$ be the unique element of $\O_n$ with $\dom(\be)=\dom(\al)$ and $\codom(\be)=\codom(\al)$.  Then $\be\in D_s(\O_n)\sub \la D_r(\O_n)\ra\sub\la D_r\ra$, and $\al = \al\al^*\cdot\be\cdot\al^*\al \in\la D_r\ra$, as we show in Figure \ref{fig:<DrMn>proj}. \epf


\begin{figure}[h]
\begin{center}
\begin{tikzpicture}[xscale=.56,yscale=0.56]
\ublueboxes{1/3,5/6,8/10,15/16,18/20}
\lblueboxes{1/2,4/8,10/11,14/15,17/20}
\uvs{1,3,4,5,6,7,8,10,15,16,17,18,20}
\lvs{1,2,3,4,8,9,10,11,14,15,16,17,20}
\stlines{4/3,7/9,17/16}
\udottedsms{10/15}
\ldottedsms{11/14}
  \draw[|-|] (.2,0)--(.2,2);
  \draw(.2,1)node[left]{$\al$};
\begin{scope}[shift={(0,-4.5)}]	
\ublueboxes{1/3,5/6,8/10,15/16,18/20}
\lblueboxes{1/3,5/6,8/10,15/16,18/20}
\uvs{1,3,4,5,6,7,8,10,15,16,17,18,20}
\lvs{1,3,4,5,6,7,8,10,15,16,17,18,20}
\stlines{4/4,7/7,17/17}
\udottedsms{10/15}
\ldottedsms{10/15}
  \draw[|-|] (.2,0)--(.2,2);
  \draw(.2,1)node[left]{$\al\al^*$};
\end{scope}
\begin{scope}[shift={(0,-6.5)}]	
\uvs{1,3,4,5,6,7,8,10,15,16,17,18,20}
\lvs{1,2,3,4,8,9,10,11,14,15,16,17,20}
\stlines{4/3,7/9,17/16}
\udottedsms{10/15}
\ldottedsms{11/14}
  \draw[|-] (.2,0)--(.2,2);
  \draw(.2,1)node[left]{$\be$};
\end{scope}
\begin{scope}[shift={(0,-8.5)}]	
\ublueboxes{1/2,4/8,10/11,14/15,17/20}
\lblueboxes{1/2,4/8,10/11,14/15,17/20}
\uvs{1,2,3,4,8,9,10,11,14,15,16,17,20}
\lvs{1,2,3,4,8,9,10,11,14,15,16,17,20}
\stlines{3/3,9/9,16/16}
\udottedsms{11/14}
\ldottedsms{11/14}
  \draw[|-] (.2,0)--(.2,2);
  \draw(.2,1)node[left]{$\al^*\al$};
\end{scope}
\end{tikzpicture}
\end{center}
\vspace{-5mm}
\caption{Diagrammatic verification of the equation $\al=\al\al^*\cdot\be\cdot\al^*\al$ from the proof of Lemma \ref{lem:<DrMn>proj}.  See the text for a full explanation.}
\label{fig:<DrMn>proj}
\end{figure}

Let $\al\in\M_n$, and suppose $A,B$ are upper non-transversal blocks of $\al$.  We say that $A$ is \emph{nested by $B$}, and write $A\prec_\al B$, if $\min(B)<\min(A)\leq\max(A)<\max(B)$.  Note that this forces $|B|=2$, but we may have $|A|=1$ or $2$.   (Also note that a block is never nested by itself.)  We define the \emph{nesting depth} of an upper non-transversal block $A$ of $\al$ to be the maximum value of $d$ such that $A\prec_\al B_1\prec_\al\cdots\prec_\al B_d$, where $B_1,\ldots,B_d$ are upper non-transversal blocks of $\al$.  For example, the blocks $\{1\}$, $\{2,10\}$, $\{11,12\}$, $\{3,6\}$, $\{7,9\}$, $\{4,5\}$, $\{8\}$ have nesting depths $0,0,0,1,1,2,2$, respectively, for
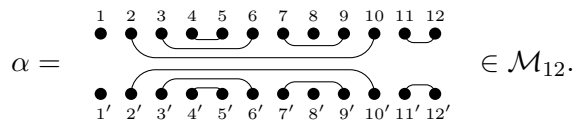
\[
\al=\dropcustpartnxy{1,...,12}{1,...,12}
{
\uarcxx36{.5}{0}
\uarcxx45{.2}{0}
\uarcxx79{.4}{0}
\uarcxx{11}{12}{.3}{0}
\uarcxx2{10}{.8}{0}
\darcxx36{.5}{0}
\darcxx45{.2}{0}
\darcxx79{.4}{0}
\darcxx{11}{12}{.3}{0}
\darcxx2{10}{.8}{0}
\vertlabelshh{1,...,12}
\vertlabeldns{1,...,12}
}\in\M_{12}.
\]
We similarly define nesting and nesting depth for lower blocks.  The next lemma will be useful on a number of occasions.

\ms
\begin{lemma}\label{lem:PDr-1Drunnested}
Let $\al\in P(D_{r-1})$, where $1\leq r\leq n-1$.  Then $\al\in\la D_r\ra$ if and only if $\al$ has an unnested singleton block, in which case $\al\in\la P(D_r)\ra$.
\end{lemma}

\pf Suppose first that $\al$ has an unnested singleton block, $\{x\}$, and let $y\in\dom(\al)^c\sm\{x\}$ be arbitrary.  Then one may easily check that $\al=\be\ga\be$, where $\ga=\id_{\dom(\al)\cup\{y\}}$, and $\be$ is obtained from $\al$ by adding the edge $\{x,x'\}$,
showing that $\al\in\la P(D_r)\ra$.  

Conversely, suppose $\al$ has no unnested singleton blocks, but that $\al=\be_1\cdots\be_k$ for some $\be_1,\ldots,\be_k\in D_r$.  Now, $\dom(\al)\sub\dom(\be_1)$ and, since also $\rank(\be_1)=\rank(\al)+1$, it follows that $\dom(\be_1)=\dom(\al)\cup\{x\}$ for some $x\in\dom(\al)^c$, and hence also $\dom(\al)^c=\dom(\be_1)^c\cup\{x\}$.  We consider separate cases according to the size of the block of $\al$ that~$x$ belongs to, 
showing that either case leads to a contradiction.

{\bf Case 1.}  Suppose first that $\{x\}$ is a block of $\al$.  Then $\{x\}$ is nested by an upper hook $\{y<z\}$ of $\al$.  Since $\dom(\be_1)=\dom(\al)\cup\{x\}$, we see that $y,z\in\dom(\be_1)^c$.  If $\{y\}$ or $\{z\}$ was a block of $\be_1$, then $\{y\}$ or $\{z\}$ would be a block of $\al=\be_1\cdots\be_k$ (respectively), a contradiction.  If $\{y,w\}$ was a block of $\be_1$ for some $w\not=z$, then $\{y,w\}$ would be a block of $\al$, a contradiction.  It follows that $\{y,z\}$ is a block of $\be_1$.  But together with $y<x<z$ and $x\in\dom(\be_1)$, this contradicts the planarity of $\be_1$.

{\bf Case 2.}  Next suppose $\{x,y\}$ is a block of $\al$ for some $y\not=x$.  This time, we would deduce that either $\{y\}$ or $\{y,z\}$ (for some $z\in\{x,y\}^c$) is a block of $\be_1$, and hence also of $\al$, a contradiction. \epf

The next result characterises the elements of $\la D_r(\M_n)\ra$.

\ms
\begin{prop}\label{prop:multinesting}
Let $\al\in D_s$, where $0\leq s\leq r\leq n-2$ and $r\not\equiv s\pmod2$.  Then $\al\in\la D_r\ra$ if and only if $\al$ has singleton upper and lower blocks each with nesting depth at most $(r-s-1)/2$.
\end{prop}

\pf We prove the lemma by induction on 
$k=(r-s-1)/2$.
Suppose first that $k=0$, so that $s=r-1$.  Then
\begin{align*}
\al\in\la D_r\ra &\iff \al\al^*,\al^*\al\in\la D_r\ra &&\text{by Lemma \ref{lem:<DrMn>proj}}\\
&\iff \text{$\al\al^*$ and $\al^*\al$ both have unnested singletons} &&\text{by Lemma \ref{lem:PDr-1Drunnested}}\\
&\iff \text{$\al$ has an unnested upper singleton and an unnested lower singleton,}
\end{align*}
establishing the result for $k=0$.  Suppose now that $k\geq1$.  Note that $s=r-2k-1\leq r-3\leq n-5$.

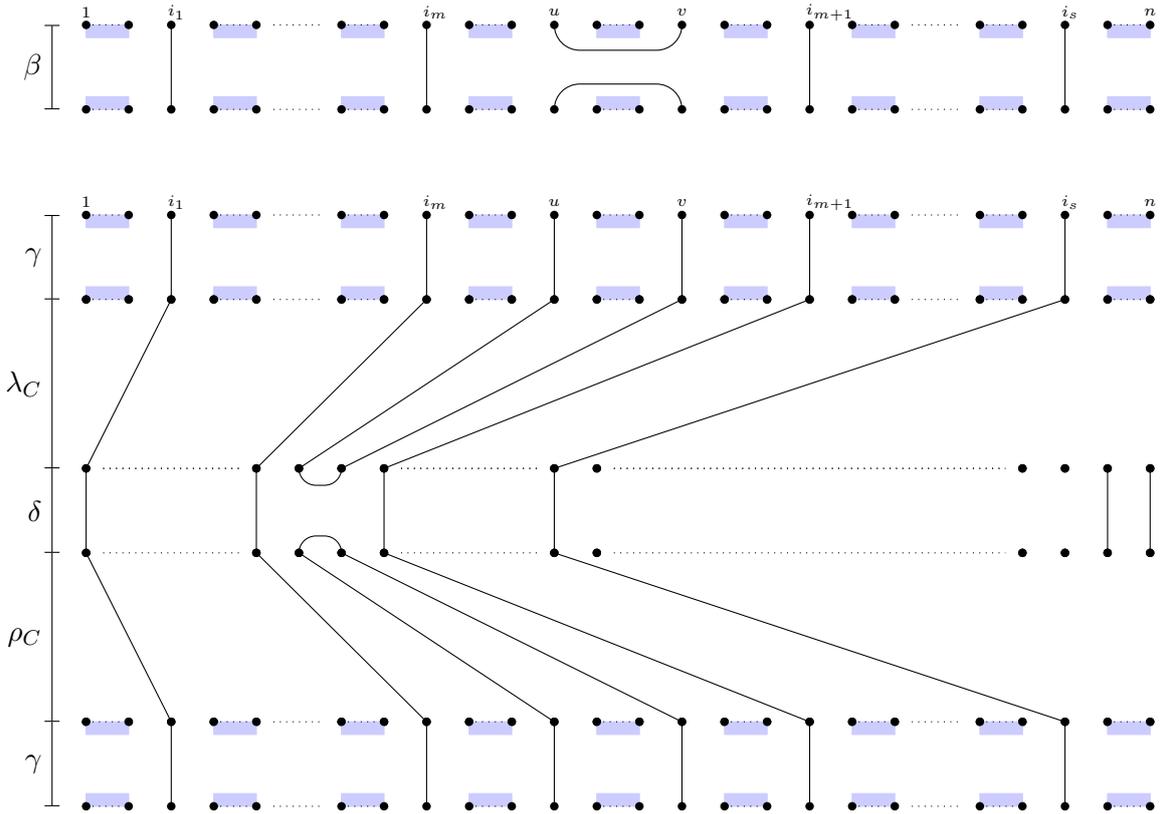
\begin{figure}[h]
\begin{center}
\begin{tikzpicture}[xscale=.56,yscale=0.56]
\ublueboxes{1/2,4/5,7/8,10/11,13/14,16/17,19/20,22/23,25/26}
\lblueboxes{1/2,4/5,7/8,10/11,13/14,16/17,19/20,22/23,25/26}
\uvs{1,2,3,4,5,7,8,9,10,11,12,13,14,15,16,17,18,19,20,22,23,24,25,26}
\lvs{1,2,3,4,5,7,8,9,10,11,12,13,14,15,16,17,18,19,20,22,23,24,25,26}
\stlines{3/3,9/9,18/18,24/24}
\udottedsms{5/7,20/22}
\ldottedsms{5/7,20/22}
\uarcx{12}{15}{.6}
\darcx{12}{15}{.6}
\vertlabels{1/1,3/\phantom{{}_1}i_1,9/\phantom{{}_m}i_m,12/u,15/v,18/\phantom{{}_{m+1}}i_{m+1},24/\phantom{{}_s}i_s,26/n}
  \draw[|-|] (.2,0)--(.2,2);
  \draw(.2,1)node[left]{$\be$};
\begin{scope}[shift={(0,-4.5)}]	
\ublueboxes{1/2,4/5,7/8,10/11,13/14,16/17,19/20,22/23,25/26}
\lblueboxes{1/2,4/5,7/8,10/11,13/14,16/17,19/20,22/23,25/26}
\uvs{1,2,3,4,5,7,8,9,10,11,12,13,14,15,16,17,18,19,20,22,23,24,25,26}
\lvs{1,2,3,4,5,7,8,9,10,11,12,13,14,15,16,17,18,19,20,22,23,24,25,26}
\stlines{3/3,9/9,12/12,15/15,18/18,24/24}
\udottedsms{5/7,20/22}
\ldottedsms{5/7,20/22}
\vertlabels{1/1,3/\phantom{{}_1}i_1,9/\phantom{{}_m}i_m,12/u,15/v,18/\phantom{{}_{m+1}}i_{m+1},24/\phantom{{}_s}i_s,26/n}
  \draw[|-|] (.2,0)--(.2,2);
  \draw(.2,1)node[left]{$\ga$};
\end{scope}
\begin{scope}[shift={(0,-8.5)}]	
\uuvs{1,2,3,4,5,7,8,9,10,11,12,13,14,15,16,17,18,19,20,22,23,24,25,26}
\lvs{1,5,6,7,8,12,13,23,24,25,26}
\stlinests{3/1,9/5,12/6,15/7,18/8,24/12}
\ldottedsms{1/5,8/12,13/23}
  \draw[|-] (.2,0)--(.2,4);
  \draw(.2,2)node[left]{$\lam_C$};
\end{scope}
\begin{scope}[shift={(0,-10.5)}]	
\uvs{1,5,6,7,8,12,13,23,24,25,26}
\lvs{1,5,6,7,8,12,13,23,24,25,26}
\stlines{1/1,5/5,8/8,12/12,25/25,26/26}
\uarc67
\darc67
\ldottedsms{1/5,8/12,13/23}
\udottedsms{1/5,8/12,13/23}
  \draw[|-] (.2,0)--(.2,2);
  \draw(.2,1)node[left]{$\de$};
\end{scope}
\begin{scope}[shift={(0,-14.5)}]	
\lvs{1,2,3,4,5,7,8,9,10,11,12,13,14,15,16,17,18,19,20,22,23,24,25,26}
\uuvs{1,5,6,7,8,12,13,23,24,25,26}
\stlinests{1/3,5/9,6/12,7/15,8/18,12/24}
  \draw[|-] (.2,0)--(.2,4);
  \draw(.2,2)node[left]{$\rho_C$};
\end{scope}
\begin{scope}[shift={(0,-16.5)}]	
\ublueboxes{1/2,4/5,7/8,10/11,13/14,16/17,19/20,22/23,25/26}
\lblueboxes{1/2,4/5,7/8,10/11,13/14,16/17,19/20,22/23,25/26}
\uvs{1,2,3,4,5,7,8,9,10,11,12,13,14,15,16,17,18,19,20,22,23,24,25,26}
\lvs{1,2,3,4,5,7,8,9,10,11,12,13,14,15,16,17,18,19,20,22,23,24,25,26}
\stlines{3/3,9/9,12/12,15/15,18/18,24/24}
\udottedsms{5/7,20/22}
\ldottedsms{5/7,20/22}
  \draw[|-] (.2,0)--(.2,2);
  \draw(.2,1)node[left]{$\ga$};
\end{scope}
\end{tikzpicture}
\end{center}
\vspace{-5mm}
\caption{Diagrammatic verification of the equation $\be = \ga\cdot \lam_C \cdot \de \cdot \rho_C \cdot\ga$ from the proof of Proposition \ref{prop:multinesting}.  See the text for a full explanation.}
\label{fig:multinesting}
\end{figure}

Suppose first that $\al$ has singleton upper and lower blocks each with nesting depth at most $k=(r-s-1)/2$.  We aim to show that $\al\in\la D_r\ra$.  The assumption on the blocks of $\al$ entails that $\al\al^*$ and $\al^*\al$ have singleton blocks of nesting depth at most $k$.  Thus, to prove that $\al\in\la D_r\ra$, it suffices, by Lemma \ref{lem:<DrMn>proj}, to show that any projection $\be\in P(D_s)$ with a singleton block of nesting depth at most $k$ belongs to $\la D_r\ra$, so suppose $\be\in P(D_s)$ is such a projection.  If $\be$ has an unnested singleton, then $\be\in\la D_{s+1}\ra\sub\la D_r\ra$, by Lemmas~\ref{lem:PDr-1Drunnested} and~\ref{lem:DsDrM}, so suppose $\be$ has no unnested singletons.  By the above assumptions, $\be$ has a singleton block $\{x\}$ with nesting depth $1\leq d\leq k$.  Suppose the outer-most block nesting $\{x\}$ is $B=\{u<v\}$.  Let $\ga$ be obtained from $\be$ by replacing the blocks $B$ and $B'$ by $\{u,u'\}$ and $\{v,v'\}$.  Since $\rank(\ga)=s+2\leq r-1$, and since $\{x\}$ is a singleton block of $\ga$ of nesting depth $d-1\leq k-1$, an induction hypothesis gives $\ga\in\la D_r\ra$.  Next, write $\dom(\be)=\{i_1<\cdots<i_s\}$, and suppose $0\leq m\leq s$ is such that $i_m<u<v<i_{m+1}$; here, we define $i_0=0$ and $i_{s+1}=n+1$.  Now let
\[
\de = \left( \begin{array}{c|c|c|c|c|c|c|c|c} 
1 \ & \ \cdots \ & \ m \ & \ m+3 \ & \ \cdots \ & \ s+2 \ & \ n-1 \ & \ n \ & \ m+1,m+2 \ \ \\ \cline{9-9}
1 \ & \ \cdots \ & \ m \ & \ m+3 \ & \ \cdots \ & \ s+2 \ & \ n-1 \ & \ n \ & \ m+1,m+2 \ \
\end{array} \!\!\! \right).
\]
Since $s\leq n-5$, we see that $\de\in P(D_{s+2})$ and that $\{n-2\}$ is an unnested singleton block of $\de$.  So Lemmas~\ref{lem:PDr-1Drunnested} and~\ref{lem:DsDrM} give $\de\in\la D_{s+3}\ra\sub\la D_r\ra$.  But then, writing $C=\dom(\be)\cup\{u,v\}$, we see in Figure \ref{fig:multinesting} that
$\be = \ga\cdot \lam_C \cdot \de \cdot \rho_C \cdot\ga$.
Since also $\lam_C,\rho_C\in D_{s+2}(\O_n)\sub\la D_r(\O_n)\ra$, it follows that $\be\in\la D_r\ra$.  As noted above, this completes the proof that $\al\in\la D_r\ra$.

Conversely, suppose $\al\in\la D_r\ra$.  We will show that $\al$ has an upper singleton block of nesting depth at most $k=(r-s-1)/2$.  Since a symmetrical argument shows that $\al$ also has a lower singleton block of nesting depth at most $k$, the proof will then be complete.  With this goal in mind, suppose to the contrary that $\al$ does \emph{not} have an upper singleton block of nesting depth at most $k$.  Write $\al=\be_1\cdots\be_m$, where $\be_1,\ldots,\be_m\in D_r$.  Since $\dom(\al)\sub\dom(\be_1)$, and since $\rank(\be_1)=r=s+2k+1=\rank(\al)+2k+1$, we have $|\dom(\be_1)\sm\dom(\al)|=2k+1$.  

{\bf Case 1.}  Suppose first that there exists $x\in\dom(\be_1)\sm\dom(\al)$ such that $\{x\}$ is a block of $\al$.  By assumption, $\{x\}$ has nesting depth at least $k+1$, so $\{x\}\prec_\al B_1\prec_\al\cdots\prec_\al B_{k+1}$ for some blocks $B_1,\ldots,B_{k+1}$ of $\al$.  Now, $B_1\cup\cdots\cup B_{k+1}\sub\dom(\al)^c$ and $|B_1\cup\cdots\cup B_{k+1}|=2k+2$, so it follows that at least one element of $B_1\cup\cdots\cup B_{k+1}$ belongs to $\dom(\be_1)^c$.  Suppose $1\leq q\leq k+1$ is such that $B_q=\{z,w\}$ and $z\in\dom(\be_1)^c$.

{\bf Subcase 1.1.}  If $\{z\}$ was a block of $\be_1$, then $\{z\}$ would be a block of $\be_1\cdots\be_m=\al$, contradicting the fact that $B_q=\{z,w\}$ is a (non-singleton) block of $\al$.  

{\bf Subcase 1.2.}  Instead, suppose $\{z,h\}$ is a block of $\be_1$ for some $h$.  Since any upper block of $\be_1$ is a block of $\al=\be_1\cdots\be_m$, it follows that $h=w$, so in fact $B_q=\{z,w\}$ is a block of $\be_1$.  But also $x\in\dom(\be_1)$ and either $w<x<z$ or $z<x<w$, since $\{x\}\prec_\al\{w,z\}$, and this contradicts the planarity of $\be_1$.

{\bf Case 2.}  Finally, suppose that every element of $\dom(\be_1)\sm\dom(\al)$ belongs to an upper hook of $\al$.  In particular, since $\dom(\be_1)\sm\dom(\al)$ has odd size, there must exist $x\in\dom(\be_1)\sm\dom(\al)$ and $y\in\dom(\be_1)^c$ such that $\{x,y\}$ is a block of $\al$.  Since $x\in\dom(\be_1)$ and $y\in\dom(\be_1)^c$, it follows that either $\{y\}$ or $\{y,z\}$ is a block of $\be_1$ for some $y\not=x$.  But then $\{y\}$ or $\{y,z\}$ would be a block of $\be_1\cdots\be_m=\al$, contradicting the fact that $\{x,y\}$ is a block of $\al$. \epf





Another major difference between $\PB_n$ and $\M_n$ concerns idempotent-generation of the ideals.
Unlike the situation with $\PB_n$, where all but the top two of the ideals $I_r(\PB_n)$ are idempotent-generated, more than half of the ideals $I_r=I_r(\M_n)$ are not idempotent-generated.  As such, we must work harder to determine which elements of $\M_n$ are products of idempotents.  To do so, we must first prove a number of technical results.  The proof of the next lemma uses the ${}^\sharp$, ${}^\natural$, ${}^\flat$ mappings defined before Lemma \ref{lem:lamAEDr}.

\ms
\begin{lemma}\label{lem:lamAn-1n}
Suppose $A\sub\bn$ is such that $n-1,n\in A^c$.  Then $\lam_A,\rho_A\in\la E(D_r)\ra$, where $r=|A|$.
\end{lemma}

\pf We just prove the statement for $\lam_A$, and we do this by induction on $n$.  The result is obvious if $n=2$ (in which case $A=\emptyset$), so suppose $n\geq3$.  Let $m=\max(A)$.  We consider two separate cases.

{\bf Case 1.}  Suppose first that $m\leq n-3$.  Since $n-2,n-1\in \bnf\sm A$, an induction hypothesis gives $\lam_A^\flat=\al_1\cdots\al_k$ for some $\al_1,\ldots,\al_k\in E(D_r(\M_{n-1}))$.  But then $\lam_A=\al_1^\natural\cdots\al_k^\natural\in\la E(D_r(\M_n))\ra$.

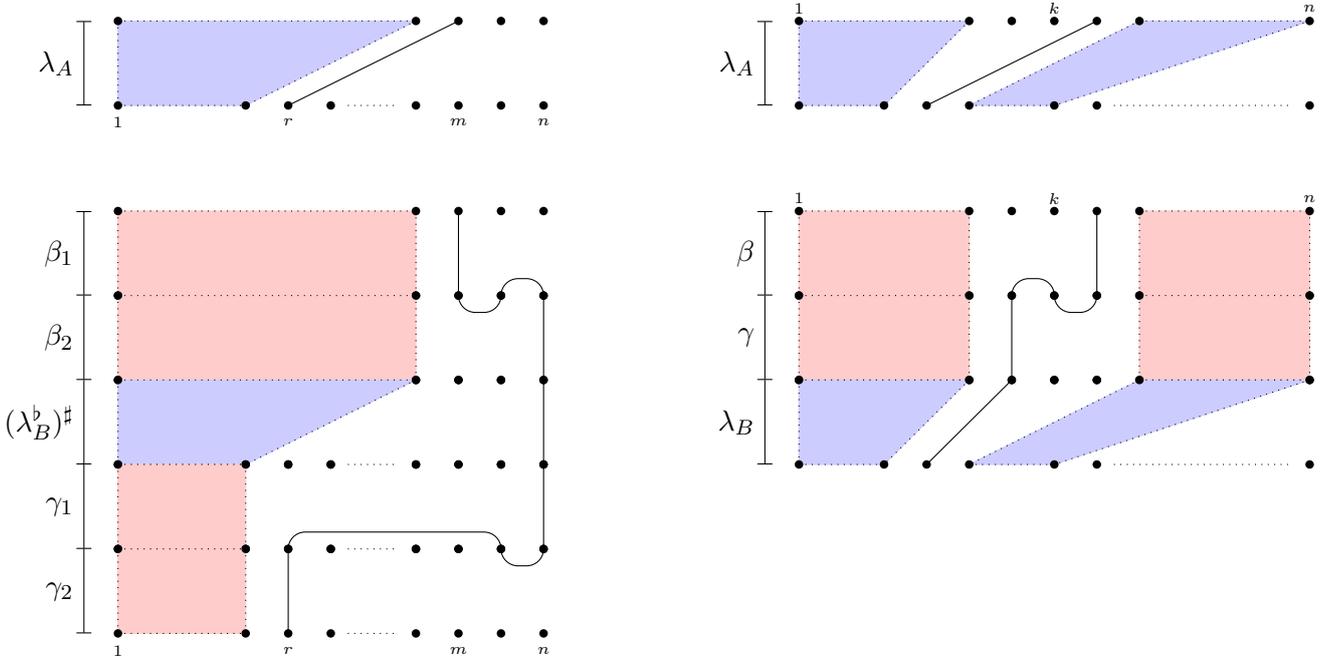
\begin{figure}[h]
\begin{center}
\begin{tikzpicture}[xscale=.56,yscale=0.56]
\bluetrap10408212
\uvs{1,8,9,10,11}
\lvs{1,4,5,6,8,9,10,11}
\stlineds{1/1,8/4}
\stlines{9/5}
\dvertlabels{1/1,5/r,9/m,11/n}
\ldottedsms{6/8}
\udotteds{1/8}
\ldotteds{1/4}
\draw[|-|] (.2,0)--(.2,2);
\draw(.2,1)node[left]{$\lam_A$};
\begin{scope}[shift={(16,0)}]	
\bluetrap10305212
\bluetrap5070{13}292
\uvs{1,5,6,7,8,9,13}
\lvs{1,3,4,5,7,8,13}
\stlineds{1/1,5/3,9/5,13/7}
\stlines{8/4}
\vertlabels{1/1,7/k,13/n}
\ldottedsms{8/13}
\udotteds{1/5,9/13}
\ldotteds{1/3,5/7}
\draw[|-|] (.2,0)--(.2,2);
\draw(.2,1)node[left]{$\lam_A$};
\end{scope}
\begin{scope}[shift={(0,-4.5)}]	
\redbox1082
\uvs{1,8,9,10,11}
\lvs{1,8,9,10,11}
\stlineds{1/1,8/8}
\stlines{9/9}
\darc{10}{11}
\udotteds{1/8}
\draw[|-|] (.2,0)--(.2,2);
\draw(.2,1)node[left]{$\be_1$};
\end{scope}
\begin{scope}[shift={(16,-4.5)}]	
\redbox1052
\redbox90{13}2
\uvs{1,5,6,7,8,9,13}
\lvs{1,5,6,7,8,9,13}
\stlineds{1/1,5/5,9/9,13/13}
\stlines{8/8}
\darc67
\vertlabels{1/1,7/k,13/n}
\udotteds{1/5,9/13}
\ldotteds{1/5,9/13}
\draw[|-|] (.2,0)--(.2,2);
\draw(.2,1)node[left]{$\be$};
\end{scope}
\begin{scope}[shift={(0,-6.5)}]	
\redbox1082
\uvs{1,8,9,10,11}
\lvs{1,8,9,10,11}
\stlineds{1/1,8/8}
\stlines{11/11}
\uarc9{10}
\udotteds{1/8}
\draw[|-] (.2,0)--(.2,2);
\draw(.2,1)node[left]{$\be_2$};
\end{scope}
\begin{scope}[shift={(16,-6.5)}]	
\redbox1052
\redbox90{13}2
\uvs{1,5,6,7,8,9,13}
\lvs{1,5,6,7,8,9,13}
\stlineds{1/1,5/5,9/9,13/13}
\stlines{6/6}
\uarc78
\udotteds{1/5,9/13}
\ldotteds{1/5,9/13}
\draw[|-] (.2,0)--(.2,2);
\draw(.2,1)node[left]{$\ga$};
\end{scope}
\begin{scope}[shift={(0,-8.5)}]	
\bluetrap10408212
\uvs{1,8,9,10,11}
\lvs{1,4,5,6,8,9,10,11}
\stlineds{1/1,8/4}
\stlines{11/11}
\ldottedsms{6/8}
\udotteds{1/8}
\ldotteds{1/4}
\draw[|-] (.2,0)--(.2,2);
\draw(.2,1)node[left]{$(\lam_B^\flat)^\sharp$};
\end{scope}
\begin{scope}[shift={(16,-8.5)}]	
\bluetrap10305212
\bluetrap5070{13}292
\uvs{1,5,6,7,8,9,13}
\lvs{1,3,4,5,7,8,13}
\stlineds{1/1,5/3,9/5,13/7}
\stlines{6/4}
\ldottedsms{8/13}
\udotteds{1/5,9/13}
\ldotteds{1/3,5/7}
\draw[|-] (.2,0)--(.2,2);
\draw(.2,1)node[left]{$\lam_B$};
\end{scope}
\begin{scope}[shift={(0,-10.5)}]	
\redbox1042
\uvs{1,4,5,6,8,9,10,11}
\lvs{1,4,5,6,8,9,10,11}
\stlineds{1/1,4/4}
\stlines{11/11}
\darc5{10}
\ldottedsms{6/8}
\udotteds{1/4}
\ldotteds{1/4}
\draw[|-] (.2,0)--(.2,2);
\draw(.2,1)node[left]{$\ga_1$};
\end{scope}
\begin{scope}[shift={(0,-12.5)}]	
\redbox1042
\uvs{1,4,5,6,8,9,10,11}
\lvs{1,4,5,6,8,9,10,11}
\stlineds{1/1,4/4}
\stlines{5/5}
\uarc{10}{11}
\dvertlabels{1/1,5/r,9/m,11/n}
\ldottedsms{6/8}
\udotteds{1/4}
\ldotteds{1/4}
\draw[|-] (.2,0)--(.2,2);
\draw(.2,1)node[left]{$\ga_2$};
\end{scope}
\end{tikzpicture}
\end{center}
\vspace{-5mm}
\caption{Diagrammatic verification of the equations $\lam_A=\be_1\be_2\cdot(\lam_{B}^\flat)^\sharp\cdot\ga_1\ga_2$ and $\lam_A=\be\ga\cdot\lam_B$ from the proofs of Lemmas \ref{lem:lamAn-1n} and \ref{lem:lamAM} --- left and right, respectively.  Red shaded parts are assumed to be identity maps on the relevant domains.  
See the text for a full explanation.}
\label{fig:lamAn-1n}
\end{figure}

{\bf Case 2.}  Now suppose $m=n-2$.  Put $B=A\sm\{n-2\}$.  Since $n-2,n-1\in\bnf\sm B$, an induction hypothesis gives $\lam_{B}^\flat\in\la E(D_{r-1}(\M_{n-1}))\ra$, and it quickly follows that $(\lam_{B}^\flat)^\sharp\in\la E(D_r(\M_n))\ra$.  In Figure~\ref{fig:lamAn-1n}, we show that $\lam_A=\be_1\be_2\cdot(\lam_{B}^\flat)^\sharp\cdot\ga_1\ga_2$, where
\begin{align*}
\be_1 &= \left( \begin{array}{c|c} 
x \ & \ \ \\ \cline{2-2}
x \ & \ n-1,n \ \
\end{array} \!\!\! \right)_{x\in A}
\COMMA
\be_2 = \left( \begin{array}{c|c} 
x \ & \ n-2,n-1 \ \ \\ \cline{2-2}
x \ & \ \
\end{array} \!\!\! \right)_{x\in B\cup\{n\}},
\\
\ga_1 &= \left( \begin{array}{c|c} 
x \ & \ \ \\ \cline{2-2}
x \ & \ r,n-1 \ \
\end{array} \!\!\! \right)_{x\in [r-1]\cup\{n\}}
\COMMA
\ga_2 = \left( \begin{array}{c|c} 
x \ & \ n-1,n \ \ \\ \cline{2-2}
x \ & \ \
\end{array} \!\!\! \right)_{x\in \br}.
\end{align*}
Since $\be_1,\be_2,\ga_1,\ga_2\in E(D_r(\M_n))$, the proof is complete. \epf

It turns out that for $\al\in\M_n$, being a product of Motzkin idempotents depends on the distribution of elements in $\dom(\al)^c$ and $\codom(\al)^c$.  With this in mind, we say a subset $A\sub\bn$ is \emph{sparse} if for all $i\in \bn$, $i\in A\implies i+1\not\in A$.  We say $A$ is \emph{cosparse} if $A^c$ is sparse.

\ms
\begin{lemma}\label{lem:lamAM}
Suppose $A\sub\bn$ is non-cosparse.  Then $\lam_A,\rho_A\in\la E(D_r)\ra$, where $r=|A|$.
\end{lemma}

\pf For a non-cosparse subset $X\sub\bn$, let $m(X)=\max\set{i\in\bn}{i-1,i\in X^c}$.  We prove the lemma by descending induction on $k=m(A)$.  Note that the case in which $k=n$ was proved in Lemma~\ref{lem:lamAn-1n}, so we assume that $k\leq n-1$.  Note that maximality of $k=m(A)$ gives $k+1\in A$.  Now put $B=(A\sm\{k+1\})\cup\{k-1\}$.  Since $k,k+1\in B^c$, we have $m(B)\geq k+1$ (but we note that it is possible to have $m(B)=k+2$).  In particular, an induction hypothesis gives $\lam_B\in\la E(D_r)\ra$.  But then the proof is complete upon noting that $\lam_A=\be\ga\cdot\lam_B$ (see Figure \ref{fig:lamAn-1n}), where
\[
\be = \left( \begin{array}{c|c} 
x \ & \ \ \\ \cline{2-2}
x \ & \ k-1,k \ \
\end{array} \!\!\! \right)_{x\in A}
\AND
\ga = \left( \begin{array}{c|c} 
x \ & \ k,k+1 \ \ \\ \cline{2-2}
x \ & \ \
\end{array} \!\!\! \right)_{x\in B}
\]
both belong to $E(D_r)$. \epf

\ms
\begin{cor}\label{cor:DrEDrMcosparse}
Let $\al\in D_r$ with $\dom(\al)$ and $\codom(\al)$ both non-cosparse.  Then $\al\in\la E(D_r)\ra$.
\end{cor}

\pf As in 
Remark \ref{rem:normalformPBn}, we have $\al=\be\cdot\lam_{\dom(\al)}\cdot\rho_{\codom(\al)}\cdot\de$, for some $\be,\de\in E(D_r)$.  Lemma \ref{lem:lamAM} gives $\lam_{\dom(\al)},\rho_{\codom(\al)}\in\la E(D_r)\ra$. \epf


We now describe the situation when one or both of $\dom(\al)$ or $\codom(\al)$ is cosparse.

\ms
\begin{lemma}\label{lem:EEMnEmn}
Let $\al\in\M_n$ be such that $\dom(\al)$ or $\codom(\al)$ is cosparse.  Then $\al\in\bbE(\M_n)$ if and only if $\al\in E(\M_n)$, in which case $\al=\id_{\dom(\al)}$.
\end{lemma}

\pf Suppose $\dom(\al)^c$ is sparse, and that $\al=\be_1\cdots\be_k$ for some $\be_1,\ldots,\be_k\in P(\M_n)$.  We claim that each of $\be_1,\ldots,\be_k$ belongs to $E(\I_n)$.  Indeed, suppose this was not the case,
and let $l\in\bk$ be minimal so that $\be_l\not\in E(\I_n)$.  So $\be_1\cdots\be_{l-1}=\id_B$ for some $B\sub\bn$.  Note that $\dom(\al)\sub\dom(\be_1\cdots\be_{l-1})=B$.  In particular, $B^c\sub\dom(\al)^c$ is also sparse.  Since $\be_l\not\in E(\I_n)$, we may write
\[
\be_l = \left( \begin{array}{c|c|c|c} 
i_1 \ & \ \cdots \ & \ i_r \ & \ A_t \ \ \\ \cline{4-4}
i_1 \ & \ \cdots \ & \ i_r \ & \ A_t \ \
\end{array} \!\!\! \right)_{t\in T},
\]
where $T\not=\emptyset$.  Fix some $u\in T$, and write $A_u=\{x<y\}$.  We consider separate cases.

{\bf Case 1.}  Suppose first that $x\in \dom(\al)\sub B$.  Then $\{x,y\}$ or $\{x\}$ would be an upper block of $\id_B\cdot\be_l$, depending on whether $y\in B$ or $y\not\in B$, respectively.  But then $\{x,y\}$ or $\{x\}$ would be an upper block of $\al=\id_B\cdot\be_l\cdots\be_k$, contradicting the fact that $x\in\dom(\al)$.

{\bf Case 2.}  A similar argument shows that we cannot have $y\in\dom(\al)$.

{\bf Case 3.}  Finally, suppose $x,y\in\dom(\al)^c$.  Since $\dom(\al)^c$ is sparse, we have $x+1\in\dom(\al)\sub B$.  If $x+1\not\in\dom(\be_l)$, then we would have $x+1\in\dom(\id_B\cdot\be_l)^c\sub\dom(\id_B\cdot\be_l\cdots\be_k)^c=\dom(\al)^c$, a contradiction.  It follows that $x+1\in\dom(\be_l)$.  But $x<x+1<y$, with $x+1\in\dom(\be_l)$ and $\{x,y\}$ an upper hook of $\be_l$, contradicting the planarity of $\be_l$.

So we have verified the original claim that $\be_1,\ldots,\be_k\in E(\I_n)$.  In particular, $\al=\be_1\cdots\be_k\in E(\I_n)\sub E(\M_n)$, and the proof is complete. \epf

\ms
\begin{prop}\label{prop:IrIGM}
Let $0\leq r\leq n$ where $n\geq2$.  Then $I_r$ is idempotent-generated if and only if $r<\floorn$. 
\end{prop}

\pf Suppose first that $r<\floorn$.  Then there are no cosparse $r$-subsets of $\bn$, and hence no elements of $I_r$ with cosparse domain or codomain.  It then follows from Corollary \ref{cor:DrEDrMcosparse} that every element of $I_r$ is a product of idempotents of rank $\leq r$.
Conversely, suppose $r\geq\floorn$.  Choose arbitrary subsets $A,B\sub\bn$ with $|A|=|B|=\floorn$ such that $A\not=B$ and $A$ is cosparse (there is a unique such set $A$ if $n$ is odd).  Then Lemma~\ref{lem:EEMnEmn} tells us that $\lam_A\rho_B\in I_r$ is not a product of idempotents. \epf

Our next task is to calculate $\rank(I_r)$ (and $\idrank(I_r)$ if appropriate).  It turns out that nested blocks play a major role in this calculation.

\ms
\begin{lemma}\label{lem:IrDrSigma}
Let $1\leq r\leq n-1$.  Then $I_r=\la D_r\cup \Si\ra$, where
\[
\Si=\set{\al\in P(D_{r-1})}{\text{\emph{$\al$ has no unnested singleton block}}}.
\]
\end{lemma}

\pf Since $I_r=\la D_r\cup D_{r-1}\ra$, it suffices to show that $D_{r-1}\sub\la D_r\cup\Si\ra$.  So let $\al\in D_{r-1}$ be arbitrary, and write
\[
\al = \left( \begin{array}{c|c|c|c} 
i_1 \ & \ \cdots \ & \ i_{r-1} \ & \ A_k \ \ \\ \cline{4-4}
j_1 \ & \ \cdots \ & \ j_{r-1} \ & \ B_l \ \
\end{array} \!\!\! \right)_{k\in K,\ l\in L}.
\]
Then $\al = \al\al^* \cdot \be \cdot \al^*\al$, as in Lemma \ref{lem:<DrMn>proj}, where $\be\in\O_n$ satisfies $\dom(\be)=\dom(\al)$ and $\codom(\be)=\codom(\al)$.
Now, $\be\in D_{r-1}(\POI_n)\sub\la D_r(\POI_n)\ra\sub\la D_r\ra$.  Also, $\al\al^*$ belongs to either $\la D_r\ra$ or $\Si$, by Lemma~\ref{lem:PDr-1Drunnested}, depending on whether $\al$ has an unnested upper singleton block or not (respectively).  Similarly, $\al^*\al$ belongs to either $\la D_r\ra$ or $\Si$.  In any case, we have proved that $\al\in\la D_r\cup\Si\ra$. \epf

The statement of the next result uses the numbers $m(n)$, $m(n,r)$ and $m'(n,r)$ defined in Section~\ref{sect:prelim}.

\ms
\begin{thm}\label{thm:IrM}
We have $\rank(I_0(\M_n))=m(n)$.  If $1\leq r\leq n-1$, then
\[
\rank(I_r(\M_n)) = m(n,r) + m'(n,r-1).
\]
Further, $I_r(\M_n)$ is idempotent-generated if and only if $n\leq1$ or $r<\floorn$, in which case $$\idrank(I_r(\M_n))=\rank(I_r(\M_n)).$$
\end{thm}

\pf The formula for $\rank(I_0)=\idrank(I_0)$ follows from Lemma \ref{lem:ARSS<D>}(iii) and Proposition \ref{prop:DClassesofMn}.  For the remainder of the proof, let $1\leq r\leq n-1$.

Since $D_r$ is a maximal $\D$-class of $I_r$, we have $\rank(I_r)=\rank(\la D_r\ra)+\rank(I_r:D_r)$.  By Lemma~\ref{lem:ARSS<D>}(ii) and Proposition \ref{prop:DClassesofMn}, $\rank(\la D_r\ra)=|P(D_r)|=m(n,r)$.  Next, we note that $|\Si|=m'(n,r-1)$, by Proposition~\ref{prop:m'nr}, where $\Si$ is the set in Lemma \ref{lem:IrDrSigma}.  It follows that $\rank(I_r:D_r)\leq m'(n,r-1)$.  

Now suppose $I_r=\la D_r\cup \Ga\ra$.  
To complete the proof that $\rank(I_r)=m(n,r)+m'(n,r-1)$, it suffices to show that $|\Ga|\geq m'(n,r-1)$.  The usual stability argument shows that for each element $\alpha$ of $\Sigma$, $\Gamma$ contains an element $\R$-related
to $\alpha$,
from which $|\Ga|\geq |\Si|=m'(n,r-1)$ quickly follows.


To complete the proof, suppose $r<\floorn$.  By Proposition \ref{prop:IrIGM}, $I_r$ is idempotent-generated.  In particular, $\la D_r\ra=\la E(D_r)\ra=\la P(D_r)\ra$, so $I_r=\la P(D_r)\cup\Si\ra$.  Thus, $P(D_r)\cup\Si$ is an idempotent generating set of size $m(n,r)+m'(n,r-1)$. \epf

Although the previous result does not apply to the entire semigroup $\M_n=I_n$, we may quickly deduce the value of $\rank(\M_n)$.

\ms
\begin{thm}\label{thm:rankMn}
We have $\rank(\M_0)=1$ and $\rank(\M_n)=2n$ for $n\geq1$.
\end{thm}

\pf The $n\leq1$ values being clear, suppose $n\geq2$.  Since $\M_n=\{1\}\cup I_{n-1}$, and since $1$ is an irreducible element of $\M_n$, we see that $I_{n-1}=\la\Ga\ra$ if and only if $\M_n=\la\Ga\cup\{1\}\ra$, so that
\[\epfreseq
\rank(\M_n)=1+\rank(I_{n-1}) = 1+m(n,n-1)+m'(n,n-2)=1+n+(n-1)=2n.
\]

\ms
\begin{rem}
Interpreting the proof of Lemma \ref{lem:ARSS<D>}(ii) in the case that $D=D_{n-1}(\M_n)$, we see that
\[
\{1\}\cup\{\be\}\cup\{\al_1,\ldots,\al_{n-1}\}\cup\{\tau_1,\ldots,\tau_{n-1}\}
\]
is a generating set for $\M_n$ of size $2n=\rank(\M_n)$, where
\[
1 =  \custpartn{1,2,5,6}{1,2,5,6}{\dotsups{2/5}\dotsdns{2/5}\stlines{1/1,2/2,5/5,6/6}\vertlabelshh{1/1,6/n}} \COMMA
\be =  \custpartn{1,2,5,6}{1,2,3,6}{\dotsups{2/5}\dotsdns{3/6}\stlines{1/2,2/3,5/6}\vertlabelshh{1/1,6/n}} \COMMA
\al_j =  \custpartn{1,3,4,5,6,8}{1,3,4,5,6,8}{\dotsups{1/3,6/8}\dotsdns{1/3,6/8}\stlines{1/1,3/3,5/4,6/6,8/8}\vertlabelshh{1/1,4/j,8/n}} \COMMA
\tau_j =  \custpartn{1,3,4,5,6,8}{1,3,4,5,6,8}{\dotsups{1/3,6/8}\dotsdns{1/3,6/8}\stlines{1/1,3/3,6/6,8/8}\uarc45\darc45\vertlabelshh{1/1,4/j,8/n}}.
\]
Note that $\{\be\}\cup\{\al_1,\ldots,\al_{n-1}\}$ generates $\O_n$ as a monoid; see \cite{Fernandes2001}, where a presentation was obtained for $\O_n$ with respect to this generating set.  A monoid presentation for $\M_n$ is given in \cite{PHY2013} in terms of the (monoid) generating set $\{\al_1,\ldots,\al_{n-1}\}\cup\{\al_1^*,\ldots,\al_{n-1}^*\}\cup\{\tau_1,\ldots,\tau_{n-1}\}$.
\end{rem}

We have seen that an ideal $I_r=I_r(\M_n)$ is only idempotent-generated when $r<\floorn$.  In particular, the idempotent-generated subsemigroup $\bbE(I_r)=\la E(I_r)\ra$ is a proper subsemigroup of $I_r$ for $r\geq\floorn$.  Our final task in this section is to describe these idempotent-generated subsemigroups $\bbE(I_r)$.

\ms
\begin{lemma}\label{lem:PDr-2PDrDr-1}
Let $2\leq r\leq n-2$.  Then $P(D_{r-2})\sub\la P(D_r\cup D_{r-1})\ra$.
\end{lemma}

\pf Let $\al\in P(D_{r-2})$ be arbitrary, and write
\[
\al = \left( \begin{array}{c|c|c|c} 
i_1 \ & \ \cdots \ & \ i_{r-2} \ & \ A_k \ \ \\ \cline{4-4}
i_1 \ & \ \cdots \ & \ i_{r-2} \ & \ A_k \ \
\end{array} \!\!\! \right)_{k\in K}.
\]
Again, if $K=\emptyset$, then $\al\in E(D_{r-2}(\O_n)) \sub\la E(D_r(\O_n))\ra\sub\la P(D_r)\ra$, so suppose $K\not=\emptyset$.  Suppose $l\in K$ is such that $A_l=\{x,y\}$, where $x\in\bn$ is minimal such that $x$ belongs to an upper hook of $\al$.  

{\bf Case 1.}  Suppose there exists $z\in\dom(\al)^c$ with $z<x$ or $y<z$.
Then $\al=\be\ga\be$ (see Figure \ref{fig:PDr-2PDrDr-1:1}, which covers the $z<x$ case)), where
\[
\be = \left( \begin{array}{c|c|c|c|c|c} 
i_1 \ & \ \cdots \ & \ i_{r-2} \ & \ x \ & \ y \ & \ A_k \ \ \\ \cline{6-6}
i_1 \ & \ \cdots \ & \ i_{r-2} \ & \ x \ & \ y \ & \ A_k \ \
\end{array} \!\!\! \right)_{k\in K\sm\{l\}}
\AND
\ga = \left( \begin{array}{c|c|c|c|c} 
i_1 \ & \ \cdots \ & \ i_{r-2} \ & \ z \ & \ x,y \ \ \\ \cline{5-5}
i_1 \ & \ \cdots \ & \ i_{r-2} \ & \ z \ & \ x,y \ \
\end{array} \!\!\! \right).
\]
Since $\be\in P(D_r)$ and $\ga\in P(D_{r-1})$, the proof is complete in this case.

\begin{figure}[h]
\begin{center}
\begin{tikzpicture}[xscale=.56,yscale=0.56]
\bluebox6{1.7}{11}2
\bluebox60{11}{.3}
\uvs{1,4,5,6,11,12,13,16}
\lvs{1,4,5,6,11,12,13,16}
\uarcx5{12}{.7}
\darcx5{12}{.7}
\vertlabels{1/1,5/x,12/y,16/n}
\udotteds{6/11}
\ldotteds{6/11}
\udottedsms{1/4,13/16}
\ldottedsms{1/4,13/16}
\stlines{1/1,4/4,13/13,16/16}
\draw[|-|] (.2,0)--(.2,2);
\draw(.2,1)node[left]{$\al$};
\begin{scope}[shift={(0,-4.5)}]	
\bluebox6{1.7}{11}2
\uvs{1,4,5,6,11,12,13,16}
\lvs{1,4,5,6,7,8,11,12,13,16}
\vertlabels{1/1,5/x,12/y,16/n}
\udotteds{6/11}
\udottedsms{1/4,13/16}
\ldottedsms{1/4,13/16,8/11}
\stlines{1/1,4/4,5/5,12/12,13/13,16/16}
\darc67
\draw[|-|] (.2,0)--(.2,2);
\draw(.2,1)node[left]{$\be$};
\end{scope}
\begin{scope}[shift={(0,-6.5)}]	
\uvs{1,4,5,6,11,12,13,16}
\lvs{1,4,5,6,7,8,11,12,13,16}
\udottedsms{1/4,13/16}
\ldottedsms{1/4,13/16,8/11}
\stlines{1/1,4/4,7/7,12/12,13/13,16/16}
\uarc56
\draw[|-] (.2,0)--(.2,2);
\draw(.2,1)node[left]{$\ga$};
\end{scope}
\begin{scope}[shift={(0,-8.5)}]	
\uvs{1,4,5,6,11,12,13,16}
\lvs{1,4,5,6,7,8,11,12,13,16}
\udottedsms{1/4,13/16}
\ldottedsms{1/4,13/16,8/11}
\stlines{1/1,4/4,5/5,6/6,13/13,16/16}
\uarcx7{12}{.6}
\darcx7{12}{.6}
\draw[|-] (.2,0)--(.2,2);
\draw(.2,1)node[left]{$\de$};
\end{scope}
\begin{scope}[shift={(0,-10.5)}]	
\uvs{1,4,5,6,11,12,13,16}
\lvs{1,4,5,6,7,8,11,12,13,16}
\udottedsms{1/4,13/16}
\ldottedsms{1/4,13/16,8/11}
\stlines{1/1,4/4,7/7,12/12,13/13,16/16}
\darc56
\draw[|-] (.2,0)--(.2,2);
\draw(.2,1)node[left]{$\ga^*$};
\end{scope}
\begin{scope}[shift={(0,-12.5)}]	
\bluebox60{11}{.3}
\lvs{1,4,5,6,11,12,13,16}
\uvs{1,4,5,6,7,8,11,12,13,16}
\ldotteds{6/11}
\ldottedsms{1/4,13/16}
\udottedsms{1/4,13/16,8/11}
\stlines{1/1,4/4,5/5,12/12,13/13,16/16}
\uarc67
\draw[|-] (.2,0)--(.2,2);
\draw(.2,1)node[left]{$\be^*$};
\end{scope}
\end{tikzpicture}
\end{center}
\vspace{-5mm}
\caption{Diagrammatic verification of the equation $\al=\be\ga\de\ga^*\be^*$ from Case 2 of the proof of Lemma~\ref{lem:PDr-2PDrDr-1}.  See the text for a full explanation.}
\label{fig:PDr-2PDrDr-1:2}
\end{figure}
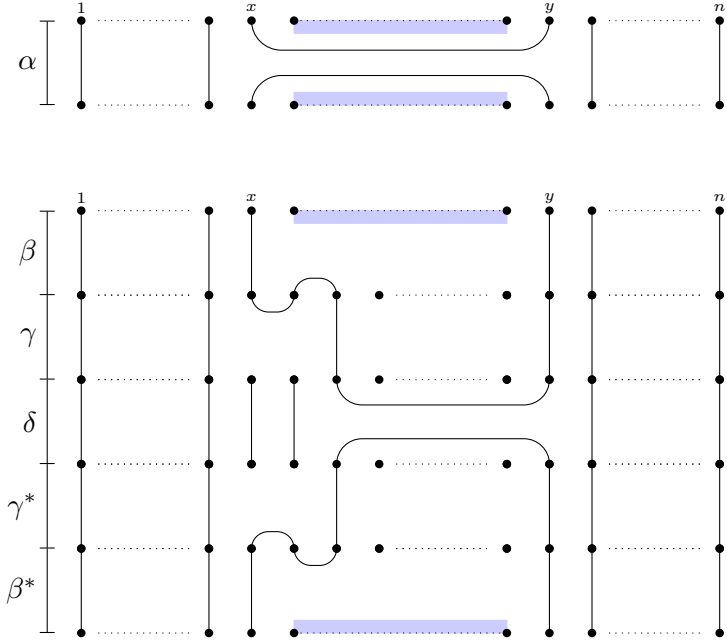

{\bf Case 2.}  Finally, suppose that $\dom(\al)^c=[x,y]$.  Since $r\leq n-2$, we see that $y-x\geq3$.  But then $\al=\be\ga\de\ga^*\be^*$ (see Figure \ref{fig:PDr-2PDrDr-1:2}), where
\begin{align*}
\be &= \left( \begin{array}{c|c|c|c|c|c} 
i_1 \ & \ \cdots \ & \ i_{r-2} \ & \ x \ & \ y \ & \ A_k \ \ \\ \cline{6-6}
i_1 \ & \ \cdots \ & \ i_{r-2} \ & \ x \ & \ y \ & \ x+1,x+2 \ \
\end{array} \!\!\! \right)_{k\in K\sm\{l\}}, \\
\ga &= \left( \begin{array}{c|c|c|c|c|c} 
i_1 \ & \ \cdots \ & \ i_{r-2} \ & \ x+2 \ & \ y \ & \ x,x+1 \ \ \\ \cline{6-6}
i_1 \ & \ \cdots \ & \ i_{r-2} \ & \ x+2 \ & \ y \ & \  \ \
\end{array} \!\!\! \right), \\
\de &= \left( \begin{array}{c|c|c|c|c|c} 
i_1 \ & \ \cdots \ & \ i_{r-2} \ & \ x \ & \ x+1 \ & \ x+2,y \ \ \\ \cline{6-6}
i_1 \ & \ \cdots \ & \ i_{r-2} \ & \ x \ & \ x+1 \ & \ x+2,y \ \
\end{array} \!\!\! \right).
\end{align*}
Since $\be,\ga,\de\in E(D_r)=P(D_r)^2$, the proof is complete in this case. \epf

\begin{figure}[h]
\begin{center}
\begin{tikzpicture}[xscale=.56,yscale=0.56]
\bluebox1032
\bluebox5072
\bluebox{13}0{17}2
\bluebox9{1.7}{11}2
\bluebox90{11}{.3}
\uvs{1,3,4,5,7,8,9,11,12,13,17}
\lvs{1,3,4,5,7,8,9,11,12,13,17}
\stlineds{1/1,3/3,5/5,7/7,13/13,17/17}
\uarcx8{12}{.7}
\darcx8{12}{.7}
\vertlabels{1/1,4/z,8/x,12/y,17/n}
\udotteds{1/3,5/7,9/11,13/17}
\ldotteds{1/3,5/7,9/11,13/17}
\draw[|-|] (.2,0)--(.2,2);
\draw(.2,1)node[left]{$\al$};
\begin{scope}[shift={(0,-4.5)}]	
\bluebox1032
\bluebox5072
\bluebox{13}0{17}2
\bluebox9{1.7}{11}2
\bluebox90{11}{.3}
\uvs{1,3,4,5,7,8,9,11,12,13,17}
\lvs{1,3,4,5,7,8,9,11,12,13,17}
\stlineds{1/1,3/3,5/5,7/7,13/13,17/17}
\stlines{8/8,12/12}
\vertlabels{1/1,4/z,8/x,12/y,17/n}
\udotteds{1/3,5/7,9/11,13/17}
\ldotteds{1/3,5/7,9/11,13/17}
\draw[|-|] (.2,0)--(.2,2);
\draw(.2,1)node[left]{$\be$};
\end{scope}
\begin{scope}[shift={(0,-6.5)}]	
\redbox1032
\redbox5072
\redbox{13}0{17}2
\uvs{1,3,4,5,7,8,9,11,12,13,17}
\lvs{1,3,4,5,7,8,9,11,12,13,17}
\stlineds{1/1,3/3,5/5,7/7,13/13,17/17}
\stlines{4/4}
\uarcx8{12}{.7}
\darcx8{12}{.7}
\udotteds{1/3,5/7,9/11,13/17}
\ldotteds{1/3,5/7,9/11,13/17}
\draw[|-] (.2,0)--(.2,2);
\draw(.2,1)node[left]{$\ga$};
\end{scope}
\begin{scope}[shift={(0,-8.5)}]	
\bluebox1032
\bluebox5072
\bluebox{13}0{17}2
\bluebox9{1.7}{11}2
\bluebox90{11}{.3}
\uvs{1,3,4,5,7,8,9,11,12,13,17}
\lvs{1,3,4,5,7,8,9,11,12,13,17}
\stlineds{1/1,3/3,5/5,7/7,13/13,17/17}
\stlines{8/8,12/12}
\udotteds{1/3,5/7,9/11,13/17}
\ldotteds{1/3,5/7,9/11,13/17}
\draw[|-] (.2,0)--(.2,2);
\draw(.2,1)node[left]{$\be$};
\end{scope}
\end{tikzpicture}
\end{center}
\vspace{-5mm}
\caption{Diagrammatic verification of the equation $\al=\be\ga\be$ from Case 1 of the proof of Lemma~\ref{lem:PDr-2PDrDr-1}.  
Red shaded parts are understood to be parts of $\id_{\dom(\al)}$.  See the text for a full explanation.}
\label{fig:PDr-2PDrDr-1:1}
\end{figure}
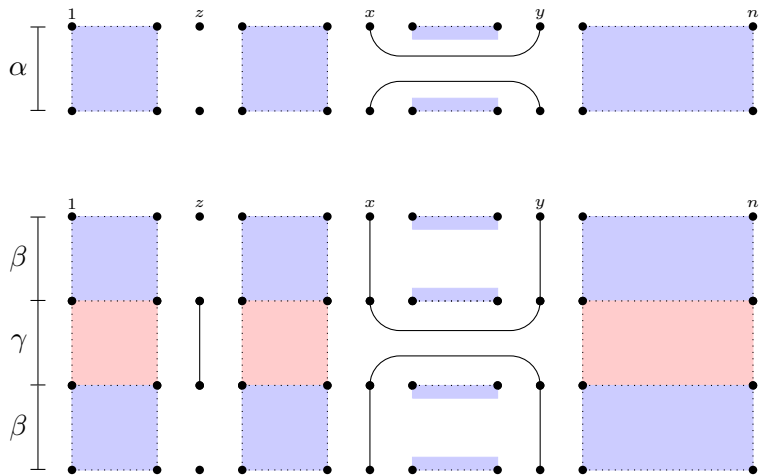

We are now able to describe the idempotent-generated subsemigroups of the ideals $I_r$ with $r\leq n-2$.  Curiously, the rank (and idempotent rank) of $\bbE(I_r)$ is equal to $\rank(I_r)$ for each such $r$, even though we do not necessarily have $I_r=\bbE(I_r)$.  Note that $\bbE(I_0)=I_0$ was covered in Theorem \ref{thm:IrM}.

\ms
\begin{thm}\label{thm:IGIrM}
Let $1\leq r\leq n-2$.  Then
\[
\bbE(I_r) = \set{\id_A}{\text{\emph{$A\sub\bn$ is cosparse and $|A|\leq r$}}} \cup \set{\al\in I_r}{\text{\emph{$\dom(\al)$ and $\codom(\al)$ are non-cosparse}}}.
\]
We have $\rank(\bbE(I_r))=\idrank(\bbE(I_r))=m(n,r)+m'(n,r-1)$.
\end{thm}

\pf The description of the elements of $\bbE(I_r)$ follows quickly from Corollary \ref{cor:DrEDrMcosparse} and Lemma \ref{lem:EEMnEmn}.  Let $\Ga=P(D_r)\cup\Si$, where $\Si\sub P(D_{r-1})$ was defined in Lemma \ref{lem:IrDrSigma}.  By Lemma \ref{lem:PDr-1Drunnested}, $P(D_{r-1})\sm\Si\sub\la P(D_r)\ra$, and it follows that $P(D_r\cup D_{r-1})\sub\la\Ga\ra$.  But then Lemma \ref{lem:PDr-2PDrDr-1} (and a simple induction) shows that $P(I_r)=P(D_r)\cup P(D_{r-1})\cup\cdots\cup P(D_0)\sub\la\Ga\ra$, so that $\bbE(I_r)=\la E(I_r)\ra=\la P(I_r)\ra\sub\la\Ga\ra$.  So $\bbE(I_r)=\la\Ga\ra$ and $\idrank(\bbE(I_r))\leq|\Ga|=m(n,r)+m'(n,r-1)$.

Now suppose $\bbE(I_r)=\la\De\ra$.  To complete the proof, it suffices to show that $|\De|\geq m(n,r)+m'(n,r-1)$.  
The usual stability argument shows that $\De$ contains an element from the $\R$-class of each element of $P(D_r)$, giving $|\De\cap D_r|\geq|P(D_r)|=m(n,r)$.  Now let $\al\in\Si$ be arbitrary, and suppose $\al=\de_1\cdots\de_k$, where $\de_1,\ldots,\de_k\in\De$.  As in the proof of Lemma \ref{lem:PDr-1Drunnested}, we cannot have $\rank(\de_1)=r$, and it quickly follows that
$\rank(\de_1)=r-1$.  As usual, stability then gives $\de_1\R\al$, and we conclude that $|\De\cap D_{r-1}|\geq|\Si|=m'(n,r-1)$.~\epf

The previous result does not apply to the idempotent-generated subsemigroup of the top two ideals $I_n=\M_n$ and $I_{n-1}=\M_n\sm\{1\}$.  We deal with these cases now.  In fact, since $\bbE(I_{n-1})=\bbE(I_n)\sm\{1\}$, it suffices to just consider $\bbE(I_n)=\bbE(\M_n)$.  Things are quite different here: $\bbE(\M_n)$ is not generated by the idempotents in its top two $\D$-classes; instead, we need (some of the) idempotents from the top \emph{four} $\D$-classes.  If $n\leq1$, then $\bbE(\M_n)=E(\M_n)=\M_n$, so we just deal with the $n\geq2$ case.

\ms
\begin{thm}\label{thm:EMn}
If $n\geq2$, then
\[
\bbE(\M_n) = \set{\id_A}{\text{\emph{$A\sub\bn$ is cosparse}}} \cup \set{\al\in\M_n}{\text{\emph{$\dom(\al)$ and $\codom(\al)$ are non-cosparse}}},
\]
and $\rank(\bbE(\M_n))=\idrank(\bbE(\M_n))=3n-2$.
\end{thm}



\pf Again, the description of the elements of $\bbE(\M_n)$ follows from Corollary \ref{cor:DrEDrMcosparse} and Lemma \ref{lem:EEMnEmn}.
Let
\[
\Ga = \{1\} \cup \{\ve_1,\ldots,\ve_n\} \cup \{\tau_1,\ldots,\tau_{n-1}\} \cup \{\mu_1,\ldots,\mu_{n-2}\},
\]
where
\[
\ve_i = \custpartn{1,3,4,5,7}{1,3,4,5,7}{\dotsups{1/3,5/7}\dotsdns{1/3,5/7}\stlines{1/1,3/3,5/5,7/7}\vertlabelshh{1/1,4/i,7/n}} \COMMA
\tau_j =  \custpartn{1,3,4,5,6,8}{1,3,4,5,6,8}{\dotsups{1/3,6/8}\dotsdns{1/3,6/8}\stlines{1/1,3/3,6/6,8/8}\uarc45\darc45\vertlabelshh{1/1,4/j,8/n}} \COMMA
\mu_k =  \custpartn{1,3,4,5,6,7,9}{1,3,4,5,6,7,9}{\dotsups{1/3,7/9}\dotsdns{1/3,7/9}\stlines{1/1,3/3,7/7,9/9}\uarcx46{.6}\darcx46{.6}\vertlabelshh{1/1,4/k,9/n}} .
\]
We first show that $\bbE(\M_n)=\la\Ga\ra$.  It suffices to show that $\la\Ga\ra\supseteq P(\M_n)$.  In fact, since $P(D_n\cup D_{n-1})\sub\Ga$, and since $P(D_{n-2}\cup\cdots\cup D_0)\sub \bbE(I_{n-2})=\la P(D_{n-2}\cup D_{n-3})\ra$, it is enough to show that $P(D_{n-2}\cup D_{n-3})\sub\la\Ga\ra$.  But this is obvious, since
\begin{align*}
P(D_{n-2}) =& \{\tau_1,\ldots,\tau_{n-1}\}\cup\set{\ve_i\ve_j}{\oijn},\\
P(D_{n-3}) =& \{\mu_1,\ldots,\mu_{n-2}\} \cup \set{\ve_i\ve_j\ve_k}{1\leq i<j<k\leq n} \\
& \phantom{ \{\mu_1,\ldots,\mu_{n-2}\}} \cup \set{\tau_i\ve_k}{\text{$1\leq k<i\leq n-1$ or $1\leq i<k-1\leq n-1$}}.
\end{align*}
Next, suppose $\bbE(\M_n)=\la\De\ra$.  To complete the proof, it suffices to show that $|\De|\geq3n-2$.
We claim that:
\bit
\item[(i)] $\De$ contains $E(D_n\cup D_{n-1})=\{1\}\cup\{\ve_1,\ldots,\ve_n\}$,
\item[(ii)] for all $1\leq j\leq n-1$, there exists $\al\in\De$ such that $\al\R\tau_j$, and
\item[(iii)] for all $1\leq k\leq n-2$, there exists $\al\in\De$ such that $\al\R\mu_k$.
\eit
Since the projections $1,\ve_1,\ldots,\ve_n,\tau_1,\ldots,\tau_{n-1},\mu_1,\ldots,\mu_{n-2}$ belong to distinct $\R$-classes of $\M_n$, it would then follow that $|\De|\geq1+n+(n-1)+(n-2)=3n-2$.  Now, (i) is obvious, and (ii) follows by the usual stability argument,
so we just prove (iii).  With this in mind, let $1\leq k\leq n-2$, and suppose $\mu_k=\al_1\cdots\al_t$, where $\al_1,\ldots,\al_t\in\De$.  If $\rank(\al_i)\geq n-1$ for all $i$, then we would have $\al_i\in\{1\}\cup\{\ve_1,\ldots,\ve_n\}$ for all $i$, in which case $\mu_k=\al_1\cdots\al_t\in E(\I_n)$, a contradiction.  Let $s\in\bt$ be minimal so that $\rank(\al_s)\leq n-2$.  Then, as before, we would have $\al_1\cdots\al_{s-1}=\id_A$ for some $A\sub\bn$.  Now, since each element of $\{k+1\}^c$ belongs to a block of $\mu_k$ of size $2$, we could only have $A=\bn$ or $A=\{k+1\}^c$.

{\bf Case 1.}  Suppose first that $A=\bn$, so that $\al_1\cdots\al_{s-1}=1$, and $\mu_k=\al_s\al_{s+1}\cdots\al_t$ with $\rank(\al_s)\leq n-2$.  If we had $\rank(\al_s)=n-2$, then (by planarity) either $\al_s$ would have two singleton upper blocks or else an upper block $\{x,x+1\}$ for some $x\in\bn$.  Since such upper blocks are not in fact blocks of $\mu_k=\al_s\al_{s+1}\cdots\al_t$, it follows that we cannot have $\rank(\al_s)=n-2$.  Since $n-3=\rank(\mu_k)=\rank(\al_s\al_{s+1}\cdots\al_t)\leq\rank(\al_s)$, it follows that $\rank(\al_s)=n-3$.  In particular, $\al_s\D\mu_k=\al_s\al_{s+1}\cdots\al_t$ so, by stability, it follows that $\al_s\R\al_s\al_{s+1}\cdots\al_t=\mu_k$, completing the proof of (iii) in this case.

{\bf Case 2.}  Finally, suppose $A=\{k+1\}^c$, so that $\al_1\cdots\al_{s-1}=\ve_{k+1}$, and $\mu_k=\ve_{k+1}\al_s\cdots\al_t$.  Now, if $\al_s$ had (at least) two upper singleton blocks, then these would be blocks of $\ve_{k+1}\al_s$ and, hence, also of $\mu_k=\ve_{k+1}\al_s\cdots\al_t$, a contradiction.  It follows that $\al_s$ has (at least) one upper hook, say $\{x,y\}$.  
Note that we cannot have $x=k+1$ or else then $\{y\}$ would be a block of $\ve_{k+1}\al_s$ and hence also of~$\mu_k$.  Similarly, $y\not=k+1$.  In particular, $x,y\in\{k+1\}^c=\dom(\ve_{k+1})$, so that $\{x,y\}$ is a block of $\ve_{k+1}\al_s$ and hence also of~$\mu_k$.  It follows that $\{x,y\}=\{k,k+2\}$.  By planarity, $\{k+1\}$ must also be a block of $\al_s$.  It then follows that $\ve_{k+1}\al_s=\al_s$, so that $\mu_k=\al_s\cdots\al_t$, and the proof concludes as in the previous case.  \epf

\ms
\begin{rem}
From the above proof, it quickly follows that any minimal (in size) idempotent generating set of $\bbE(\M_n)$ is of the form $E(D_n\cup D_{n-1})\cup\Ga\cup\De$, where $\Ga$ is a minimal idempotent generating set of $\mathcal J_n\sm\{1\}$, and $\De\sub E(D_{n-3})$ is such that: $|\De|=n-2$, and for all $1\leq k\leq n-2$, $\De$ contains an element $\R$-related to $\mu_k$ and an element $\L$-related to $\mu_k$.  It is easy to check that $R_{\mu_k}\cap L_{\mu_l}$ is non-empty if and only if $k=l$, so that we can only have $\De=\{\mu_1,\ldots,\mu_{n-2}\}$.  It follows that the minimal idempotent generating sets of $\bbE(\M_n)$ are in one-one correspondence with the minimal idempotent generating sets of $\mathcal J_n\sm\{1\}$.  Such generating sets were classified and enumerated (in terms of Fibonacci numbers) in \cite[Theorem 9.8]{EastGray}.
\end{rem}

\ms
\begin{rem}
Presentations for the idempotent-generated subsemigroups $\bbE(\B_n)$ and $\bbE(\P_n)$ of the Brauer monoid $\B_n$ and the partition monoid $\P_n$ were given in \cite{Maltcev2007} and \cite{JEpnsn}, respectively.  These presentations used generating sets of minimal size consisting entirely of projections.  It would therefore be interesting to try and find a presentation for $\bbE(\M_n)$ in terms of the generating set $\Ga$ from the proof of Theorem \ref{thm:EMn}, which is of minimal size and consists entirely of projections.
\end{rem}

\section{Applications to diagram algebras}\label{sect:algebras}

In this section, we apply results of previous sections to obtain important information about the representation theory of the corresponding diagram algebras: namely, the \emph{partial Brauer algebras} \cite{HD2014,MarMaz2014} and \emph{Motzkin algebras}~\cite{BH2014}.  Specifically, we combine general results of Wilcox \cite{Wilcox2007} with our normal forms (Lemma \ref{lem:normalformPBn} and Remark \ref{rem:normalformPBn}) to deduce cellularity of these algebras, and we also use our combinatorial results on regular $*$-semigroups and Green's classes of diagram monoids (Propositions \ref{prop:RSS}, \ref{prop:DClassesofKn} and \ref{prop:DClassesofMn}) to determine the dimensions of the cell representations (which are precisely the irreducible representations in the semisimple case).  Although the results we obtain here are already known \cite{BH2014,HD2014,MarMaz2014}, our techniques are different to existing methods.  We include these derivations to demonstrate the utility of the semigroup theoretic approach to diagram algebras, which we believe will be of use in future investigations.  

\subsection{Cellular algebras}\label{subsect:cellular}

Cellular algebras were introduced in \cite{GL1996} to provide a unified framework for studying the representation theory of many important algebras, including Brauer algebras, Temperley-Lieb algebras, and certain Hecke algebras.  Many other examples of cellular algebras are now known \cite{Wilcox2007,JEcais,Xi2000,Guo2009,Xi1999,AM2000, GW2015,WY2012,Geck2007}.  See also \cite{May2015,May2015_2}, where a broader class of algebras, called \emph{cell algebras}, were introduced in order to capture even more natural examples.
We now recall the definitions and ideas from \cite{GL1996} that we will need.

Let $R$ be a commutative ring with $1$, and $A$ a finite dimensional unital $R$-algebra.  Recall 
that $A$ is a \emph{cellular algebra} with \emph{cell datum} $(\Lam,M,\C,*)$ if:
\bit
\item[(i)] $\Lam$ is a finite partially ordered set,
\item[(ii)] for all $\lam\in\Lam$, $M(\lam)$ is a finite set,
\item[(iii)] $\C=\set{C_{\fs,\ft}^\lam}{\lam\in\Lam,\ \fs,\ft\in M(\lam)}$ is a basis of $A$, and the map $(\lam,\fs,\ft)\mt C_{\fs,\ft}^\lam$ is injective,
\item[(iv)] the map ${}^*:A\to A$ determined on basis elements (and then extended linearly) by $(C_{\fs,\ft}^\lam)^*=C_{\ft,\fs}^\lam$ is an $R$-linear (anti-)involution of $A$,
\item[(v)] for all $\lam\in\Lam$, $\fs,\ft\in M(\lam)$ and $a\in A$,
\begin{equation}\label{eq:cell}
aC_{\fs,\ft}^\lam \equiv \sum_{\fs'\in M(\lam)}r_a(\fs',\fs)C_{\fs',\ft}^\lam \pmod{A(<\lam)},
\end{equation}
where each $r_a(\fs',\fs)\in R$ is independant of $\ft$, and where $A(<\lam)$ is the $R$-submodule of $A$ spanned by $\set{C_{\fu,\fv}^\mu}{\mu\in\Lam,\ \mu<\lam,\ \fu,\fv\in M(\mu)}$.  
\eit
From equation \eqref{eq:cell} (and its $*$-dual), it follows that for all $\lam\in\Lam$ and $\fs_1,\fs_2,\ft_1,\ft_2\in M(\lam)$,
\[
C_{\fs_1,\ft_1}^\lam C_{\fs_2,\ft_2}^\lam \equiv \phi(\ft_1,\fs_2)C_{\fs_1,\ft_2}^\lam \pmod{A(<\lam)}
\]
for some $\phi(\ft_1,\fs_2)\in R$ that depends only on $\ft_1$ and $\fs_2$.  

For each $\lam\in\Lam$, there is a natural left $A$-module $W(\lam)$, called the \emph{cell module} corresponding to $\lam$; $W(\lam)$ has basis $\set{C_{\fs}}{\fs\in M(\lam)}$, and the action of $A$ on $W(\lam)$ is given by
\begin{equation}\label{eq:cell2}
aC_{\fs} = \sum_{\fs'\in M(\lam)}r_a(\fs',\fs)C_{\fs'},
\end{equation}
where the scalars $r_a(\fs',\fs)$ are defined in \eqref{eq:cell}.  
There is a natural bilinear form $\phi_\lam:W(\lam)\times W(\lam)\to R$ defined by $\phi_\lam(C_{\fs},C_{\ft})=\phi(\fs,\ft)$ for each $\fs,\ft\in M(\lam)$.  The \emph{radical} of $\lam$ is then defined to be the $A$-submodule of $W(\lam)$
\[
\rad(\lam) = \set{x\in W(\lam)}{\phi_\lam(x,y)=0\ (\forall y\in W(\lam)},
\]
and we define the quotient modules $L_\lam=W(\lam)/\rad(\lam)$.  
It was shown in \cite[Theorem 3.4(i)]{GL1996} that when $R$ is a field, $\set{L_\lam}{\lam\in\Lam,\ \phi_\lam\not=0}$ is a complete set of (representatives of isomorphism classes of) absolutely irreducible $A$-modules.  
In particular, the dimensions of these irreducible $A$-modules are given by
\[
\dim(L_\lam)=|M(\lam)|-\dim(\rad(\lam)).
\]
In general, calculating the values of $\dim(\rad(\lam))$ is very difficult; see for example \cite[page 48]{Mathas1999}.
It was shown in \cite[Theorem 3.8]{GL1996} that when $R$ is a field, the following are equivalent:
\bit
\item[(i)] $A$ is semisimple,
\item[(ii)] the non-zero cell modules $W(\lam)$ are irreducible and pairwise non-isomorphic,
\item[(iii)] the form $\phi_\lam$ is non-degenerate (i.e., $\rad(\lam)=0$) for each $\lam\in\Lam$.
\eit

\subsection{The partial Brauer algebra}\label{subsect:PBnt}

Let $\al,\be\in\PB_n$.  Recall that when forming the product $\al\be\in\PB_n$, we created the product graph $\Ga(\al,\be)$.  In general, this graph could have several connected components that involve vertices from $\bn''$ (in the middle row); we call these \emph{floating components} (of $\Ga(\al,\be)$).  Floating components come in two kinds: \emph{loops} and \emph{paths}.  Following \cite{MarMaz2014}, we consider singleton floating components to be paths (of length $0$) rather than loops.  We write $l(\al,\be)$ and $p(\al,\be)$ for the number of floating loops and paths in $\Ga(\al,\be)$, respectively.  For example, with $\al,\be\in\PB_{12}$ as in Figure \ref{fig:multinPB8}, we have $l(\al,\be)=1$ and $p(\al,\be)=2$.  
It is easy to check that
\begin{equation}\label{eq:landp}
\left.
\begin{array}{c}
l(\al,\be)+l(\al\be,\ga)=l(\al,\be\ga)+l(\be,\ga) \\
p(\al,\be)+p(\al\be,\ga)=p(\al,\be\ga)+p(\be,\ga) 
\end{array}
\right\}
\qquad\text{for all $\al,\be,\ga\in\PB_n$.}
\end{equation}
Now let $R$ be a commutative ring with identity, and fix some $x,y\in R$.  We define the \emph{twisting map}
\[
\tau:\PB_n\times\PB_n\to R \qquad\text{by}\qquad \tau(\al,\be)=x^{l(\al,\be)}y^{p(\al,\be)}.
\]
For example, with $\al,\be\in\PB_{12}$ as in Figure~\ref{fig:multinPB8}, we have $\tau(\al,\be)=xy^2$.  
The \emph{partial Brauer algebra} $\PBnt$ is the associative $R$-algebra with basis $\PB_n$ and multiplication $\star$ defined on basis elements (and then extended linearly) by 
\[
\al\star\be = \tau(\al,\be) \al\be \qquad\text{for $\al,\be\in\PB_n$.}
\]
(Note that associativity follows from \eqref{eq:landp}.)  In this way, we see that $\PBnt$ is a \emph{twisted semigroup algebra} of $\PB_n$ with respect to the twisting map $\tau$.  As a special case, if $x=y=1$, then $\PBnt=R[\PB_n]$ is the (ordinary) semigroup algebra of $\PB_n$.  It is clear that 
\begin{align}
\label{eq:tw1} \tau(\al,\be)&=\tau(\be^*,\al^*) &&\text{for any $\al,\be\in\PB_n$}\\
\label{eq:tw2} \tau(\al,\be)&=\tau(\al,\ga) &&\text{for any $\al,\be,\ga\in\PB_n$ with $\be\R\ga$.}
\end{align}
%
For $0\leq r\leq n$, recall that the $\H$-class (in $\PB_n$) of the projection
\[
\id_{\br} = \left( \begin{array}{c|c|c|c} 
1 \ & \ \cdots\ & \ r\ & \ \phantom{\emptyset} \ \ \\ \cline{4-4}
1 \ & \ \cdots\ & \ r\ & \ \phantom{\emptyset}  \ \
\end{array}\!\!\! \right) \in P(D_r(\M_n)) , 
\]
was denoted $\S_{\br}$ just before Lemma \ref{lem:normalformPBn}, and is isomorphic to the symmetric group $\S_r$.  It is well known that the (ordinary) symmetric group algebra $R[\S_r]$ is cellular.
%
We fix a cell datum $(\Lam_r,M_r,\C_r,{}^{-1})$ for the group algebra $\RSr$ for each $0\leq r\leq n$,
where $(\Lam_r,\dominates)$ is the set of all integer partitions of $r$ ordered by dominance, $M_r(\lam)$ is the set of all standard tableaux of shape $\lam$ for each $\lam\in\Lam_r$, and $\C_r=\set{C_{\fs,\ft}^\lam}{\lam\in\Lam_r,\ \fs,\ft\in M_r(\lam)}$ is the \emph{Murphy basis}, as described in \cite{Murphy1992} and \cite[Chapter 3]{Mathas1999}.
We now describe a cell datum $(\Lam,M,\C,*)$ for the partial Brauer algebra $\PBnt$.  

First, put $\Lam=\Lam_0\cup\Lam_1\cup\cdots\cup\Lam_n$ with order $\leq$ defined, for $\lam\in\Lam_i$ and $\mu\in\Lam_j$ by
\[
\text{$\lam\leq\mu \iff i<j$ or [$i=j$ and $\lam\dominates\mu$].}
\]
For $0\leq r\leq n$ and $\lam\in\Lam_r$, we define $M(\lam)=P(D_r(\PB_n))\times M_r(\lam)$.  (Recall that $M_r$ is part of the cell datum for $\RSr$.)  

For a projection $\al\in P(D_r(\PB_n))$, write $\xi_\al=\be\cdot\lam_{\dom(\al)}$ and $\zeta_\al=\rho_{\codom(\al)}\cdot\de$, where $\be,\de$ are as in the proof of Lemma \ref{lem:normalformPBn}.  For $0\leq r\leq n$, $\lam\in\Lam_r$, $\al,\be\in P(D_r(\PB_n))$ and $\fs,\ft\in M_r(\lam)$, define
\[
C_{(\al,\fs),(\be,\ft)}^\lam = \xi_\al \cdot C_{\fs,\ft}^\lam \cdot \zeta_\be\in\PBnt.
\]
Since $C_{\fs,\ft}$ is an $R$-linear combination of elements from $\S_{\br}$, we see that $C_{(\al,\fs),(\be,\ft)}^\lam$ is an $R$-linear combination of elements from the $\H$-class of $\xi_\al\cdot\zeta_\be\in D_r(\PB_n)$.  Now put
\[
\C = \bigcup_{r=0}^n \set{C_{(\al,\fs),(\be,\ft)}^\lam}{\lam\in\Lam_r,\ \fs,\ft\in\Lam_r,\ \al,\be\in P(D_r(\PB_n))}.
\]
Then, by the above discussion, \cite[Corollary 7]{Wilcox2007} tells us that $(\Lam,M,\C,*)$, as constructed above, is a cell datum for $\PBnt$.  Note that the anti-involution $*$ on $\PBnt$ is the $R$-linear extension of the anti-involution on $\PB_n$ we have been using throught the article.
The discussion in Section \ref{subsect:cellular} allows us to construct the irreducible representations of $\PBnt$ in the case that $R$ is a field.  In particular, the dimensions of the cell modules $W(\lam)$ are given by
\[
\dim(W(\lam)) = \big|M(\lam)\big| = \big|P(D_r(\PB_n))\big|\cdot\big|M_r(\lam)\big| = {n \choose r} \cdot a(n-r)\cdot\big|M_r(\lam)\big| \qquad\text{for each $\lam\in\Lam_r$,}
\]
by Propositions \ref{prop:DClassesofKn}(iv) and \ref{prop:RSS}(iii), where the numbers $a(m)$ are defined in Proposition \ref{prop:DClassesofKn}.
Formulae for the values $|M_r(\lam)\big|$ are well known
\cite[Exercise 3.25]{Mathas1999}.
In particular, in the case that $\PBnt$ is semisimple,\footnote{When $R=\bbC$ is the field of complex numbers, it is known that $\PBntC$ is generically semisimple \cite[Theorem 1.1 and Corollary 3.6]{MarMaz2014}; i.e., $\PBntC$ is semisimple for all but a finite number of choices of $x,y$.}  
the irreducible $\PBnt$-modules are precisely the cell modules, so $L_\lam=W(\lam)$ for each $\lam\in\Lam$, with dimensions as stated above.

\subsection{Motzkin algebras}\label{subsect:Mnt}

Again let $R$ be a commutative ring with $1$.  
The $R$-submodule of the partial Brauer algebra $\PBnt$ spanned by $\M_n\sub\PB_n$ is a subalgebra of $\PBnt$.  This subalgebra is the \emph{Motzkin algebra} \cite{BH2014}, and we denote it by $\Mnt$.  Again, we may apply results of previous sections to quickly demonstrate the cellularity of $\Mnt$ and calculate the dimensions of the cell modules (which coincide with the irreducible $\Mnt$-modules in the semisimple case).  Since $\M_n$ is aperiodic (i.e., all its subgroups are trivial), the situation is a little easier for $\Mnt$ than for $\PBnt$; this contrasts interestingly with earlier sections, in which we generally had to work harder for $\M_n$ than for $\PB_n$.

This time, we define $\Lam=\{0,1,\ldots,n\}$ ordered by $\leq$.  For $r\in\Lam$, we define $M(r)=P(D_r(\M_n))$ to be the set of all rank $r$ Motzkin projections.  For $r\in\Lam$ and $\al,\be\in M(r)$, we write
\[
C_{\al,\be}^r = \xi_\al\cdot\zeta_\be,
\]
where $\xi_\al$ and $\zeta_\be$ were defined in the previous section (see also Lemma \ref{lem:normalformPBn} and Remark \ref{rem:normalformPBn}).  We then define
\[
\C = \set{C_{\al,\be}^r}{r\in\Lam,\ \al,\be\in M(r)}.
\]
Then by \cite[Corollary 7]{Wilcox2007} again, $\Mnt$ is cellular, with cell datum $(\Lam,M,\C,*)$, where again $*$ denotes the natural anti-involution on $\M_n$.  In fact, by Lemma \ref{lem:normalformPBn} and Remark \ref{rem:normalformPBn}, we see that the cellular basis $\C$ is precisely $\M_n$ itself.  Again, when $R$ is a field, in the semisimple case,\footnote{In \cite[Theorem 5.14]{BH2014}, semisimplicity was characterised in the case $x=y$ in terms of certain expressions involving Chebyshev polynomials.} the irreducible $\Mnt$-modules are precisely the cell modules as defined generally in Section \ref{subsect:cellular}.  In fact, because of the simpler nature of the Motzkin algebra, we can describe these modules quite easily.  For $0\leq r\leq n$, the cell module $W(r)$ has basis $\set{C_\al}{\al\in M(r)}$.  (Recall that $M(r)=P(D_r(\M_n))$.)  The action of $\Mnt$ on $W(r)$ from \eqref{eq:cell2} is given, for $\al\in\M_n$ and $\be\in M(r)$, by
\[
\al\cdot C_\be = \begin{cases}
\tau(\al,\be)C_{\al\be\al^*} &\text{if $\rank(\al\be)=r$}\\
0 &\text{if $\rank(\al\be)<r$.}
\end{cases}
\]
We note that this corresponds to the geometric definition given in \cite[Section 4.1]{BH2014}.  We also note that the dimensions of these cell modules (which are precisely the irreducible modules in the semisimple case) are given by
\[
\dim(W(r)) = \big|M(r)\big|=\big|P(D_r(\M_n))\big|=m(n,r),
\]
by Propositions \ref{prop:DClassesofMn} and \ref{prop:RSS}(iii), where the Motzkin numbers $m(n,r)$ are defined in Proposition \ref{prop:DClassesofMn}.

\section{Calculated values}\label{sect:values}

In this final section, we provide calculated values of the ranks (and idempotent ranks, where appropriate) of the ideals $I_r(\PB_n)$ and $I_r(\M_n)$ as well as the idempotent-generated subsemigroups $\bbE(\PB_n)$ and $\bbE(\M_n)$.

\begin{table}[h]%
\begin{center}
\begin{tabular}{|c|rrrrrrrrrrr|}
\hline
$n\sm r$	&0	&1	&2	&3	&4	&5	&6	&7	&8	&9	&10 \\
\hline
0	&1										&&&&&&&&&&\\
1	&1	&2									&&&&&&&&&\\
2	&2	&3	&3								&&&&&&&&\\
3	&4	&6	&6	&4							&&&&&&&\\
4	&10	&19	&12	&11	&4						&&&&&&\\
5	&26	&50	&55	&20	&16	&4					&&&&&\\
6	&76	&171	&150	&125	&30	&22	&4				&&&&\\
7	&232	&532	&651	&350	&245	&42	&29	&4			&&&\\
8	&764	&1961	&2128	&1876	&700	&434	&56	&37	&4		&&\\
9	&2620	&6876	&9297	&6384	&4536	&1260	&714	&72	&46	&4	&\\
10	&9496	&27 145	&34 380	&32 565	&15 960	&9702	&2100	&1110	&90	&56	&4 \\
\hline
\end{tabular}
\end{center}
\vspace{-5mm}
\caption{Values of $\rank(I_r(\PB_n))$.  These values also give $\idrank(I_r(\PB_n))$ for $r\leq n-2$.}
\label{tab:rank(Ir(PBn))}
\end{table}

\begin{table}[h]%
\begin{center}
\begin{tabular}{|c|rrrrrrrrrrr|}
\hline
$n\sm r$	&0	&1	&2	&3	&4	&5	&6	&7	&8	&9	&10 \\
\hline
0	&1										&&&&&&&&&&\\
1	&1	&2									&&&&&&&&&\\
2	&2	&3	&4								&&&&&&&&\\
3	&4	&6	&5	&6							&&&&&&&\\
4	&9	&15	&11	&7	&8						&&&&&&\\
5	&21	&36	&32	&17	&9	&10					&&&&&\\
6	&51	&91	&83	&56	&24	&11	&12				&&&&\\
7	&127	&232	&226	&157	&88	&32	&13	&14			&&&\\
8	&323	&603	&608	&459	&266	&129	&41	&15	&16		&&\\
9	&835	&1585	&1655	&1305	&832	&419	&180	&51	&17	&18	&\\
10	&2188	&\phantom{3 }4213	&\phantom{3 }4517	&\phantom{3 }3726	&\phantom{3 }2499	&1397	&\phantom{3}626	&\phantom{3}242	&62	&19	&20 \\
\hline
\end{tabular}
\end{center}
\vspace{-5mm}
\caption{Values of $\rank(I_r(\M_n))$.  These values also give $\idrank(\bbE(I_r(\M_n)))$ for $r\leq n-2$.}
\label{tab:rank(Ir(Mn))}
\end{table}

\begin{table}[H]%
\begin{center}
{\footnotesize
\begin{tabular}{|c|rrrrrrrrrrrrrrrrrrrrr|}
\hline
$n$ &0	&1	&2	&3	&4	&5	&6	&7	&8	&9	&10 &11  &12  &13  &14  &15  &16  &17  &18  &19  &20 \\
\hline
$\rank(\bbE(\PB_n))$ &1	&2	&4	&7	&11	&16	&22	&29	&37	&46	&56  &67  &79  &92  &106  &121  &137  &154  &172  &191  &211\\
\hline
$\rank(\bbE(\M_n))$ &1	&2	&4	&7	&10	&13	&16	&19	&22	&25	&28 &31  &34  &37  &40  &43  &46  &49  &52  &55  &58 \\
\hline
\end{tabular}
}
\end{center}
\vspace{-5mm}
\caption{Values of $\rank(\bbE(\PB_n))=\idrank(\bbE(\PB_n))$ and $\rank(\bbE(\M_n))=\idrank(\bbE(\M_n))$.}
\label{tab:rank(E(PBn))}
\end{table}

\section*{Acknowledgements}

The first named author gratefully acknowledges the support of Grant No.~174019 of the Ministry of Education, Science, and Technological Development of the Republic of Serbia, and Grant No.~0851/2015 of the Secretariat of Science and Technological Development of the Autonomous Province of Vojvodina.  This project was initiated while the authors attended the LMS--EPSRC Durham Symposium, \emph{Permutation Groups and Transformation Semigroups}, in July 2015, organised by Peter Cameron, Dugald Macpherson and James Mitchell; we thank the organisers for the opportunity.  The {\sc Semigroups} package of GAP \cite{GAP} was useful during preliminary stages; we also thank James Mitchell and Attila Egri-Nagy for helpful conversations.

\footnotesize
\def\bibspacing{-1.1pt}
\bibliography{biblio}
\bibliographystyle{plain}
\end{document}